\documentclass[10pt]{article}

\usepackage[utf8]{inputenc}

\usepackage{graphicx}
\usepackage{amsmath}
\usepackage{amsthm}
\usepackage{amssymb}
\usepackage{float}
\usepackage{caption}
\usepackage[hidelinks]{hyperref}
\usepackage{xcolor}
\usepackage{color}
\usepackage{mathtools}
\usepackage{dsfont}
\usepackage{comment}
\usepackage{enumitem}

\usepackage{fullpage} 

\usepackage[ruled,vlined,nofillcomment]{algorithm2e}
\usepackage{stmaryrd} 
\usepackage{mathrsfs} 
\usepackage{multirow} 

\usepackage{xspace} 
\usepackage{graphbox} 
\usepackage{float} 

\newif\ifdraft
\drafttrue 

\usepackage{todonotes} 

\def\argmin_#1{\underset{#1}{\mathrm{argmin\, }}}

\def \E{\mathbb{E}}
\newcommand{\be}{\begin{eqnarray}}
\newcommand{\ee}{\end{eqnarray}}
\newcommand{\beno}{\begin{eqnarray*}}
\newcommand{\eeno}{\end{eqnarray*}}
\newcommand{\barr}[1]{\begin{array}{#1}}
\newcommand{\earr}{\end{array}}

\newcommand{\R}{\mathbb{R}}
\newcommand{\eps}{\varepsilon}
\newcommand{\ma}{\alpha}
\newcommand{\mb}{\beta}
\newcommand{\mg}{\gamma}

\newcommand{\mk}{\kappa}
\newcommand{\mO}{\Omega}
\newcommand{\ms}{\sigma}
\newcommand{\mt}{\theta}
\newcommand{\mT}{\Theta}

\newcommand{\cA}{{\mathcal A}}

\newcommand{\converge}{\rightarrow}
\newcommand{\conv}{\converge}

\definecolor{SpecGreen}{cmyk}{0.64,0,0.95,0.40}
\definecolor{SpecGray}{cmyk}{0.64,0,0.05,0.40}
\definecolor{apcolor}{rgb}{0.1,0.5,0.4}

\newcommand{\COMMENTO}[1]{}   

\newcommand{\UN}{\mbox{\rm 1\hspace{-0.6ex}I}}
\newcommand{\ba}{{\bar a}}

\newcommand{\supp}{supp}

\newcommand{\SLTWO}{SL-scheme\xspace}
\newcommand{\SLTHREE}{L-scheme\xspace}
\newcommand{\XA}{X} 
\newcommand{\XcA}{X} 
\newcommand{\YcA}{Y} 
\newcommand{\numA}{{\hat{\mathcal A}}} 
\newcommand{\numV}{{\hat{\mathcal V}}} 

\newtheorem{thm}{Theorem}[section]
\newtheorem{prop}[thm]{Proposition}
\newtheorem{corol}[thm]{Corollary}
\newtheorem{defi}[thm]{Definition}
\newtheorem{lem}[thm]{Lemma}
\newtheorem{rem}[thm]{Remark}

\title{Neural networks for first order HJB equations and application to front propagation with obstacle terms%
\footnote{This research benefited from the support of the FMJH Program PGMO and from the support 
to this program from EDF.
}
}
\author{Olivier Bokanowski,%
\footnote{Université Paris Cité, Laboratoire Jacques-Louis Lions (LJLL), F-75013 Paris, France}
\footnote{Sorbonne Université, CNRS, LJLL, F-75005 Paris, France \texttt{olivier.bokanowski@u-paris.fr}}
\
Averil Prost,%
\footnote{INSA Rouen, LMI (EA 3226 - FR CNRS 3335), 76801 St Etienne du Rouvray, France. \texttt{averil.prost@insa-rouen.fr}}
\
Xavier Warin%
\footnote{EDF R\&D \& FiME, 91120 Palaiseau, France \texttt{xavier.warin@edf.fr}}
}

\begin{document}

\date{\today} 

\maketitle

\begin{abstract}
We consider a deterministic optimal control problem with a maximum running cost functional, in a finite horizon context,
and propose deep neural network approximations for Bellman's dynamic programming principle,
corresponding also to some first order Hamilton-Jacobi-Bellman equation.
This work follows the lines of Huré  et al.
(SIAM J. Numer. Anal., vol. 59 (1), 2021, pp. 525-557) where algorithms are proposed
in a stochastic context.
However we need to develop
a completely new approach in order to deal with the propagation of errors in the
deterministic setting, where no diffusion is present in the dynamics. 
Our analysis gives precise error estimates in an average norm. 
The study is then illustrated on several academic numerical examples related to
front propagations models in presence of obstacle constraints,
showing the relevance of the approach for average dimensions
(e.g. from $2$ to $8$), even for non smooth value functions.
\end{abstract}




\newcommand{\ha}{{\hat a}}
\newcommand{\hb}{{\hat b}}
\newcommand{\hU}{{\hat U}}
\newcommand{\hV}{{\hat V}}
\newcommand{\dt}{{\Delta t}}
\newcommand{\cN}{{\mathcal N}}
\newcommand{\cB}{{\mathcal B}}

\textwidth 15cm


\newcommand{\refexrot}{1}
\newcommand{\refexeik}{2}
\newcommand{\refexadv}{3}
\newcommand{\refeiknc}{4}
\newcommand{\refeikc}{5}



\newcommand{\barre}{\hbox{\rule[1ex]{\textwidth}{.1cm}}}
\newcommand{\barrethin}{\hbox{\rule[1ex]{\textwidth}{0.02cm}}}
\renewcommand{\div}{\mathrm{div}}

\newcommand{\hypfive}{(H4)}
\newcommand{\hypfourA}{(H5)}
\newcommand{\hypfourB}{(H6)}

\section{Introduction}

In this work we are interested by the approximation of
a deterministic optimal control problem with finite horizon involving a maximum running cost, defined as
\begin{align}\label{pb:1}
        v(t,x) = \inf_{a(\cdot) \in \mathbb{A}_{[t,T]}} \max \bigg( \max_{\mt\in[t,T]} g(y_{x}^a(\mt)),\ \varphi(y_{x}^a(T)) \bigg),
\end{align}
where the state $x$ belongs to $\R^d$ and $t\in [0,T]$ for some $T\geq 0$.
Here the trajectory $y(s)=y_{x}^a(s)$ obeys the following dynamics
\begin{align}\label{eq:defysOpenLoop}
        \dot y(s) = f(y(s),a(s)), \quad a.e.\ s\in [t,T],
\end{align}
with initial condition $y(t)=x$, and control $a \in \mathbb{A}_{[t,T]} \coloneqq L^{\infty}\left([t,T],A\right)$. It is assumed that 
$A$ is a non-empty compact subset of $\R^\mk$ ($\mk\geq 1$)
and ($f$, $\varphi$, $g$) are Lipschitz continuous. 
The value $v$ is solution of the following Hamilton-Jacobi-Bellman (HJB) partial differential equation, in the viscosity sense
(see for instance \cite{bok-for-zid-2010-2})
\begin{subequations}\label{eq:hjb-pde}
\be 
  & & \min\bigg(- v_t  + \max_{a\in A} (-f(x,a)\cdot \nabla_x v), \ v- g(x)\bigg) = 0, \quad t\in [0,T] \\
  & & v(T,x)=\max(\varphi(x),g(x)).
\ee
\end{subequations}

Tremendous numerical efforts have been made in order to propose efficient algorithms 
for solving problem  related to \eqref{pb:1}, or the corresponding HJB equation 
\eqref{eq:hjb-pde}. Precise numerical methods have been developed, using approximations on grids, such as 
Markov Chains approximations \cite{kus-dup-01}, finite difference schemes  (monotone schemes \cite{cra-lio-84},  
semi-Lagrangian schemes (see e.g. \cite{deb-jak-14,fal-fer-16})
ENO or WENO higher-order schemes \cite{Osher_Shu_91,Shu99}, 
finite element methods (see \cite{jensen-smears-2013}), 
discontinuous Galerkin methods \cite{Hu_1999_SIAM_DG_FEM, Li_2005_AML_DG_HamiJaco},
and in particular \cite{bok-che-shu-2014,bok-che-shu-2016} for \eqref{eq:hjb-pde}, see also \cite{smears-suli-2016},
or max-plus approaches \cite{aki-gau-lak-08}).
However, grid-based methods are limited to low dimensions because of the well-known curse of dimensionality.
In order to tackle this difficulty, various approaches are studied, such as
spare grids methods \cite{bok-gar-gri-klo-2013,garcke_kroner_2017},
tree structure approximation algorithm (see e.g.\cite{alla_falcone_Saluzzi_tree_2020}),
tensor decomposition methods \cite{dolgov_tensor_2021}, max-plus approaches in \cite{mac-des-gau-08}.


\if{
\begin{rem}[Links with classical problems]\label{rem:links-with-class-pbs}
  The following Bolza problem with state constraints
  \be \label{pb:0-1}
  u(t,x) = 
    \inf_{a(.)\in\mathbb{A}_{[t,T]}} \left\{ \int_t^T \ell(y_x^a(\mt),a(\mt)) d\mt + \varphi(y_x^a(T)), 
     \quad y_x^a(\theta) \in K \ \ \forall \theta\in [t,T] \right\}
  \ee
  may be recast into the form \eqref{pb:1} 
  by augmenting the dynamics (thus leading to a control problem with dynamics in dimension $d+1$), see~\cite{alt-bok-zid-2013}.
\end{rem}
}\fi


In the deterministic context,
problem~\eqref{pb:1} is motivated by deterministic optimal control with state constraints 
(see e.g.~\cite{bok-for-zid-2010-2} and \cite{alt-bok-zid-2013}).
In \cite{ban-tom-21}, the HJB equation \eqref{eq:hjb-pde} is approximated by deep neural networks (DNN)
for solving state constrained reachability control problems 
of dimension up to $d=10$. 
In \cite{bok-gam-zid-22} or in  \cite{ass-bok-des-zid-18}, 
formulation \eqref{pb:1} is used to solve an abort landing problem (using different numerical approaches);
in \cite{bok-bou-des-zid-sso-17}, equations such as \eqref{eq:hjb-pde} are used to solve an aircraft payload optimization problem;
a multi-vehicle safe trajectory planning is considered in \cite{ban-che-tan-tom-21}.

%



On the other hand, for stochastic control, DNN approximations were already used 
for gas storage optimization  
in \cite{barrera2006numerical}, where the control approximated by a neural
network was the amount of gas injected or withdrawn in the storage.
This approach has been adapted and popularized recently for the resolution of BSDE
in~\cite{han-jen-e-18} (deep BSDE algorithm). 
For a convergence study of such algorithms in a more general context, see \cite{han_convergence_2020}.

In this work, we study some neural networks approximations for \eqref{pb:1}.
We are particularly interested for the obtention of a rigorous error analysis of such approximations.
We follow the approach of \cite{hur-pha-21-a} (and its companion paper~\cite{hur-pha-21-b}),
combining neural networks approximations and Bellman's dynamic programming principle.
We obtain precise error estimates in an average norm.

Note that the work of \cite{hur-pha-21-a} is developed in the stochastic context, where an error analysis
is given. 
However this error analysis somehow relies strongly on a diffusion assumption of the model (transition probabilities with 
densities are assumed to exists).
In our case, we would need to assume that the deterministic process admits a density, which is not the case
(see remark~\ref{rem:hur-not-for-deterministic}).
Therefore the proof of \cite{hur-pha-21-a} does not apply to the deterministic context.
Here we propose a new approach for the convergence analysis, leading to new error estimates.
We chose to present the algorithm on a running cost optimal control problem, 
but the approach can be generalized
to Bolza or Mayer problems (see e.g.~\cite{alt-bok-zid-2013,ass-bok-des-zid-18}).

%
%

For sake of completeness, let us notice that the ideas of \cite{hur-pha-21-a} 
are related to methods already proposed in \cite{gobet2005regression} and \cite{bender2007forward} 
for the resolution of Backward Stochastic Differential Equations (BSDE),
where the control function is calculated by regression on a space of some basis functions
(the {\em Hybrid-Now} algorithm is related to \cite{gobet2005regression},
and the {\em performance iteration} algorithm is related to an improved algorithm in \cite{bender2007forward}).
For recent developments, see \cite{gobet2016linear} using classical linear regressions,
and \cite{hure2020deep} and \cite{germain2022approximation} for BSDE approximations using neural networks.





From the numerical point of view, we illustrate our algorithms on some 
academic front-propagation problems with or without obstacles.
We focus on a "Lagrangian scheme"
(a deterministic equivalent of the {\em performance iteration} scheme of \cite{hur-pha-21-a}),
and also compare with other algorithms :
a "semi-Lagrangian algorithm" (similar to the {\em Hybrid-Now} algorithm of \cite{hur-pha-21-a})
and an hybrid algorithm combining the two previous,
involving successive projection of the value function on neural network spaces.

The plan of the paper is the following.
In section~2 we define a semi-discrete value approximation for \eqref{pb:1}, with controlled error with respect to the continuous value, using piecewise constant controls.
In section~3, equivalent reformulations of the problem are given using feedback controls and dynamic programming.
In section~4, an approximation result of the discrete value function by using Lipschitz continuous feedback controls is given.
In section~5 we present three numerical schemes using neural networks approximations (for the approximation of feedback controls and/or for the value),
using general Runge Kutta schemes for the approximation of the controlled dynamics.
Section~6 contains our main convergence result for one of the proposed scheme
(the Lagrangian scheme) which involves
only approximations of the feedback controls, and section~7 focuses on the proof of our main result.
Section~8 is devoted to some numerical academic examples of front-propagation problems 
with or without an obstacle term (state constraints), 
for average dimensions, 
showing the potential of the proposed algorithms in this context, and also giving comparisons between
the different algorithms introduced.
An appendix contains some details for computing reference solutions for some of the considered examples.


{\bf Notations.}
Unless otherwise precised, the norm $|.|$ on $\R^q$ ($q\geq 1$) 
is the max norm $|x|=\|x\|_\infty=\max_{1\leq i\leq q} |x_i|$.
The notation $\llbracket p,q\rrbracket=\{p,p+1,\dots,q\}$ is used, for any integers $p\leq q$.
For any function $\ma:\R^p\conv \R^{q}$ for some $p,q\geq 1$, $[\ma]:=\sup_{y\neq x}\frac{|\ma(y)-\ma(x)|}{|y-x|}$ denotes
the corresponding Lipschitz constant. We also denote $a\vee b :=\max(a,b)$ for any $a,b\in \R$.
The set of "feedback" controls is defined as 
$ \cA:=\{a:\R^d\conv A, \ a(.)\ \mbox{measurable}\}$.


\section{Semi-discrete approximation with piecewise constant controls} 

In this section, we first aim to define a semi-discrete approximation of \eqref{pb:1} in time.

Let the following assumptions hold on the set $A$ and functions $f,g,\varphi$.

\medskip


\newcommand{\hypo}{(H0)\xspace}
\noindent\textbf{\hypo}
{\em
$A$ is a non-empty compact subset of $\R^\mk$ ($\mk\geq 1$), and is a convex set.
}

\medskip

\newcommand{\hypf}{(H1)\xspace}
\noindent {\textbf\hypf}
{\em
$f:\R^d\times A \conv \R^d$ is Lipschitz continuous and we denote 
$[f]_{1},[f]_2 \geq 0$ constants such that
\beno
 |f(x,a)-f(x',a')| \ \leq \ [f]_1|x-x'| + [f]_2|a-a'|, \quad \forall (x,x')\in (\R^d)^2,\ \forall (a,a')\in A^2.
\eeno
}

\newcommand{\hypG}{(H2)\xspace}
\noindent{\textbf \hypG}
{\em
$g:\R^d\conv \R$ is Lipschitz continuous.
}

\medskip

\newcommand{\hypPhi}{(H3)\xspace}
\noindent\textbf{\hypPhi}{\ \em
$\varphi:\R^d\conv \R$ is Lipschitz continuous.
}


\medskip


Let $T>0$ be the horizon, let $N\in\mathbb{N}^*$ be a number of iterations, and $(t_k)_{k\in\llbracket0,N\rrbracket} \subset [0,T]$ 
be a time mesh with $t_0 = 0$ and $t_N = T$. 
To simplify the presentation, we restrict ourselves to the uniform mesh $t_k = k \dt$ with $\dt = \frac{T}{N}$, 
but the arguments would carry over unchanged with a non-uniform time mesh.

Let us consider $F_h:\R^d\times A\conv\R^d$ (for a given $h>0$),
corresponding to some one time step approximation of $y^a_{x}(h)$ (starting from $y^a_x(0)=x$).
For instance, we may consider the Euler scheme $F_h(x,a) = x + h f(x,a)$, or
the Heun scheme  $F_h(x,a) = x + \frac{h}{2} (f(x,a) + f(x+ h f(x,a),a))$, and so on.
General Runge Kutta schemes
are considered later on in section~\ref{sec:RK}. Assumptions on $F_h$ will be made precise when needed.

For a given sequence $a \coloneqq \left(a_n, a_{n+1}, \dots, a_{N-1}\right) \in A^{N-n}$ (with $n\in\llbracket0,N-1\rrbracket$), which corresponds to a 
piecewise constant control approximation, and a given integer $p\geq 1$,
we define two levels of approximations for the trajectories.

The fine approximation involves time step $h=\frac{\dt}{p}$, 
is denoted $Y_{j,x}^{a_k}$ (for a fixed control $a_k$, $\forall 0\leq j\leq p$),
and is defined  recursively by
\begin{subequations}
\label{eq:defYcA}
\be
    & & Y^{a_k}_{0,x} = x \\
    & & Y^{a_k}_{j+1,x} = F_h(Y^{a_k}_{j,x},a_k), \quad j=0,\dots,p-1, 
\ee
\end{subequations}
which also corresponds to $j$ iterates of $y\conv F_h(y,a_k)$, starting from $y=x$, with the same control $a_k$.
This fine level will be used to obtain approximation of the trajectory 
at intermediate time steps $t_k + j h$ which lie into $[t_k,t_{k+1}]$.

The coarse approximation with time step $\dt$ 
is denoted $(X^a_{k,x})_{n\leq k\leq N}$ and is defined  recursively by
\begin{subequations}
\label{eq:defXnOpenLoop}
\be
  & & X_{n,x}^a \coloneqq x  
  \label{eq:defXnOpenLoop-a} \\
  & & X_{k+1,x}^a = Y^{a_k}_{p, X^a_{k,x}}, \quad k=n,\dots,N-1.
  \label{eq:defXnOpenLoop-b}
\ee
\end{subequations}
{\bf Notations.} We will often use the notations, for a given $a\in A$,
 $F(\cdot,a)\equiv F^a(\cdot) := Y_{p,.}^a$ and the fact that 
 \eqref{eq:defXnOpenLoop-b} can also be written
$$
  X^a_{k+1,x} = F(X^a_{k,x},a_k), \quad  k=n,\dots,N-1.
$$

We can now define the following cost functional,
for $x\in \R^d$, $a\in A^{N-n}$ 
and $n \in \llbracket0,N\rrbracket$
\begin{align}\label{eq:defJn} 
    J_n (x,a) \coloneqq 
      \max\limits_{n\leq k < N}  
      \max\limits_{0\leq j<p} g(Y^{a_k}_{j,X_{k,x}^a})
        \bigvee (g\vee \varphi) (\XcA_{N,x}^a), \quad x\in \R^d,\ a \in A^{N-n},
\end{align}
and the following semi-discrete version of \eqref{pb:1}, 
for $x \in \R^d$ and $n \in \llbracket0,N\rrbracket$
\be\label{eq:defVn} 
     V_n(x) := \min_{a\in A^{N-n}} J_n(x,a).
\ee
(for $n=N$, we have $V_N(x)=J_N(x)=g(x)\vee\varphi(x)$).
It will be also useful to introduce the following notation, for $a\in A$ and $x\in \R^d$,
\be\label{eq:Ga}
  G^a(x) := \max\limits_{0\leq j<p} g(Y_{j,x}^a).
\ee

The values $(V_n)_{0\leq n\leq N}$ satisfy also $V_N(x)=g(x)\vee \varphi(x)$, and the following dynamic programming principle (DPP) for $n=0,\dots,N-1$:
\begin{align}\label{eq:DPP-1}
  V_n(x) = \inf\limits_{a\in A}  G^a(x) \vee V_{n+1}(F(x,a)), \quad x\in \R^d.
\end{align}

Let us notice that the case $p=1$ leads to the following simplifications: $h=\dt$, $F(x,a)=F_\dt(x,a)$,  $G^a(x)=g(x)$,
$J_n(x)=\max_{n\leq k\leq N} g(X^a_{k,x}) \vee \varphi(X^a_{N,x})$, as well as $V_N(x)=g(x)\vee\varphi(x)$ and the DPP
$V_n(x)=  \inf\limits_{a\in A}  g(x) \vee V_{n+1}(F_\dt(x,a))$ ($0\leq n\leq N-1$).

The motivation behind the introduction of the finer level of approximation $(Y_{j,x}^{a_k})_{0\leq j\leq p}$ is first numerical.
It enables a better evaluation of the running cost term $g(.)$ along the trajectory, without the computational cost of more intermediate controls.
The numerical improvement is illustrated in the examples of section~\ref{ex:rot}. 
Furthermore, from the theoretical point of view,
the convergence analysis in our main result will strongly 
use the fact that $x\conv F_h(x,a)$ is a change of variable for $h$ sufficiently small 
(i.e., $p$ sufficiently large).

\medskip

We start by showing some uniform Lipschitz bounds. 

\begin{lem}\label{lem:vn-lip-bound}
  Assume (H0)-(H3), and the Lipschitz bound  $\sup_{a\in A} [F_h(.,a)]\leq 1 + ch$ for some constant $c\geq 0$.

  $(i)$
  The function $J_n(.,a)$ is Lipschitz for all $a\in A^{N-n}$, 
  with uniform bound $[J_n(.,a)] \leq [g]\vee[\varphi] e^{cT}$. 

  $(ii)$ 
  In particular, the uniform bound $\max_{0\leq n\leq N}[V_n] \leq [g]\vee[\varphi] e^{cT}$ holds.

\end{lem}

\begin{proof}

  $(i)$
  Notice that $1+ch\leq e^{ch}$. Then for $a\in A$ and for the $j$-th iterate $F_h^{(j)}(.,a)$ of $F_h$, we obtain
  $|Y^a_{j,x}-Y^a_{j,y}|=|F_h^{(j)}(x,a) - F_h^{(j)}(y,a)| \leq e^{jc h} |x-y|\leq e^{c\dt } |x-y|$ for any $0\leq j \leq p$ (by recursion).
  Hence also $\sup_{a\in A} [F(.,a)] \leq e^{c\dt}$,
  from which we deduce for any $a=(a_n,\dots,a_{N-1})\in A^{N-n}$ and $n\leq k\leq N$, 
  $|X^a_{k,x} - X^a_{k,y}| \leq  e^{c (k-n) \dt} |x-y| \leq e^{c T} |x-y|$. 
  The desired result follows from the definition of $J_n$ and repeated use of $\max(a,b)-\max(c,d)\leq \max(a-c,b-d)$.

  $(ii)$
  As a direct consequence of $(i)$ and the definition of $V_n$.
\end{proof}

The following result shows that $V_n(x)$ is a first order approximation of $v(t_n,x)$ in time.

\begin{thm} \label{thm:VnA_close_to_v}
  Assume (H0)-(H3), and that there exists $h_0>0$ such that
  \begin{itemize}
  \item[-]
  $F_h$ is consistent with the dynamics $f$ in the following sense:
  \be\label{eq:F-consistent}
      \exists C\geq 0,\ \forall (x,a,h)\in \R^d\times A\times(0,h_0), \quad
      |F_h(x,a) - (x + h f(x,a))|\leq C h^2 (1+|x|),
  \ee
  \item[-]
    for all $h\in ]0,h_0]$,
    $\sup_{a\in A} [F_h(.,a)]\leq 1 + c h$ for some constant $c\geq 0$ ($c$ may depends on $[f]$),

  \item[-]
  $f(x,A)$ is convex for all $x\in \R^d$. 
  \end{itemize}
  Let $h=\frac{\dt}{p}\leq h_0$ (with $p\geq 1$).
  Then 
   \be \label{eq:error-estim-v}
      \max_{0\leq n\leq N} |V_n(x)-v(t_n,x)| \leq \widetilde{C}\dt(1+|x|)
   \ee
  for some constant $\widetilde{C}\geq 0$ independent of $x,\dt,p$ (and $N$).
\end{thm}
\begin{proof}
  This follows from the arguments of Theorem B.1. in \cite{bok-gam-zid-22}.
\end{proof}

\begin{corol}
  In particular, if $F_h$ is a consistent RK scheme (see definition \ref{def:RK}) 
  and $f(x,A)$ is convex for all $x$,
  then by Lemma~\ref{lem:RK-is-F-consistent},
  \eqref{eq:F-consistent} holds and therefore the error estimate \eqref{eq:error-estim-v} also holds.
\end{corol}

Our aim is now to propose numerical schemes for the approximation of $V_n(.)$.


\section{Reformulation with feedback controls} \label{sec:feedback}





In this section, equivalent definitions for $V_n$ are given using feedback controls in $\cA$
(the set of measurable functions $a:\R^d\conv A$).
These formulations will lead to the numerical schemes. 


First, for a given $a_k\in \cA$, 
the fine approximation $Y_{x,j}^{a_k}$ 
(with time step $h=\frac{\dt}{p}$)
is defined by $Y_{x,j}^{a_k} := Y_{x,j}^{a_k(x)}$
using definition \eqref{eq:defYcA}.
This corresponds also to
\be
  \label{eq:defYcA-TWO}
    & & Y^{a_k}_{0,x} = x \\
    & & Y^{a_k}_{j+1,x} = F_h(Y^{a_k}_{j,x},a_k(x)), \quad j=0,\dots,p-1 
\ee
(that is, $j$ iterates of $y\conv F_h(y,a_k(x))$, starting from $y=x$, with the fixed control $a_k(x)$).
Then, for a given sequence $a \coloneqq \left(a_n, \dots, a_{N-1}\right) \in \cA^{N-n}$,
the coarse approximation is defined by
\begin{subequations} \label{eq:defXnFeedback}
\be
     & &  X_{n,x}^a = x \\
     & &  X_{k+1}^a  = F^{a_k}(X_{k,x}^a) \equiv F\left(\XcA_{k,x}^{a}, a_k(\XcA_{k,x}^a)\right),
     \quad k=n,\dots,N-1.
\ee
\end{subequations}
(with notation $F^{a_k}(x)=F(x,a_k(x))$.
We also extend the definition of $J_n$ to the feedback space, 
for $x\in \R^d$ and $a\in \cA^{N-n}$, as follows 
\begin{align}\label{eq:defJnE}
  J_n (x,a) \coloneqq  \big(\max_{n\leq k \leq N} G^{a_k}(\XcA_{k,x}^{a})  \big)
   \bigvee \varphi(\XcA_{N,x}^a)
\end{align}
where now, for a given control $a\in \cA$, we extend the definition of \eqref{eq:Ga} by
\be\label{eq:Ga-TWO}
  G^a(x) := \max\limits_{0\leq j<p} g(Y_{j,x}^{a(x)})
\ee
(this also corresponds to define $G^a(x)$ as $G^{a(x)}(x)$).
With this definitions, we have the following results.

\begin{prop} \label{prop:dpp-Vn-cA}

  $(i)$ 
  $V_n(x)$ is the minimum of $J_n(x,.)$ over feedback controls: 
  $$ 
    V_n(x) = \min_{a\in \cA^{N-n}} J_n(x,a), \quad x\in \R^d.
  $$

  $(ii)$ For all $0\leq n\leq N-1$, 
  $V_n$ satisfies the following
  dynamic programming principle over feedback controls
  \be\label{eq:DPP-1-vn-feedback} 
    V_n(x) = \min_{a\in \cA} G^a(x) \bigvee V_{n+1}(F^a(x)), \quad x\in \R^d, \ n=0,\dots,N-1
  \ee 
  and in particular, the infimum is reached by some some $\ba_n\in \cA$.
\end{prop}
\begin{proof}
The problem is to show the existence of a measurable feedback control.
By using a measurable selection procedure 
(see for instance Lemmas~2A,~3A p. 161 of~\cite{leeFoundationsOptimalControl1986}), 
  since $A$ is compact and $G$ (such that $G^a(x)=G(x,a(x))$), $V_{n+1}$ and $F$ (such that $F^a(x)=F(x,a(x))$) 
  are continuous, we may choose $\ba_n$ measurable in \eqref{eq:DPP-1-vn-feedback}.
\end{proof} 

The following well-known result links pointwise minimization over open-loop controls $a\in A$  and minimization 
of an averaged value over feedback controls $a\in \cA$.
\begin{lem}
  \label{lem:equiv-minimization}
  Let $X$ be a random variable with values in $\R^d$
  which admits a density $\rho$, and such that $\E[|X|]<\infty$.
  Then for any measurable $\mO\subset\R^d$ such that $\rho(x)>0$ a.e. $x\in \mO$, and $n\in \llbracket 0,N-1\llbracket$,
  \beno
   & & \hspace{-2cm}\ba_n(.) \in \argmin_{a\in \cA } \E\bigg[\UN_\mO(X)\ G^a(X) \vee V_{n+1} (F^a(X))\bigg]  \\
   & & \quad \Longleftrightarrow \quad
   \bigg( \ba_h(x) \in \argmin_{a\in A }\big( G^a(x) \vee V_{n+1} (F^a(x))\big),\ \mbox{a.e. $x\in \mO$} \bigg).
  \eeno 
\end{lem}

We now introduce a new assumption on a sequence of sets $\mO_n$, densities $\rho_n$ (supported in $\bar\mO_n$) and associated random variables $X_n$ 
(with associated probability densities $\rho_n$).

\medskip
\noindent{\textbf \hypfive}
\begin{em}
    The functions $\rho_k\in L^1(\R^d)$ and open sets $\mO_k\subset\R^d$, $0\leq k\leq N$, are such that
\begin{subequations} 
\be
  & &
    \mbox{$\rho_k(x)>0$ on $\mO_k$, and $\supp(\rho_k)\subset \overline{\mO_k}$,  $\forall k=0,\dots,N$}
    \label{eq:H5-a}
    \\[1ex]
  & &
    \mbox{$F(\mO_k,a)\subset \mO_{k+1}$, $\forall a\in A$, 
    $\forall k=0,\dots,N-1$} 
    \label{eq:H5-b}
    \\
  & &
    C_{k,\dt}:= \sup_{x\in \mO_k} \sup_{a\in A} \frac{\rho_k(x)}{\rho_{k+1}(F(x,a))}<\infty,
    \ \forall k=0,\dots,N-1.
    \label{eq:H5-c}
\ee
\end{subequations} 
  Furthermore, 
  we consider random variables $(X_k)_{0\leq k\leq N}$ on some probability space, with values in $\R^d$, 
  absolutely continuous with respect to Lebesgue's measure and admitting 
  $(\rho_k)_{0\leq k\leq N}$ as associated densities.
\end{em}

From the definitions we have $\E[\phi(X_k)] = \int_{\mO_k} \phi(x) \rho_k(x) dx$ for any measurable bounded function $\phi$.

\medskip

The technical assumption \eqref{eq:H5-c} is not important in this section but will be needed for the main result later on.
Before going on, we give some examples where \hypfive\ holds:

\begin{itemize}
\item
  case $\mO_0$ is a bounded subset of $\R^d$, 
 $\mO_k=\mO_0 + \mathscr{B}(0, c k\dt \|f\|_\infty)$ (where $\mathscr{B}(0,r)$ is the ball of radius $r$,
 $\|f\|_\infty$ is a bound for $|f|$ on $\mO_N\times A$, 
    and assuming $|F_h(x)|\leq |x| + c h \|f\|_\infty$ as it will be the case for RK schemes as in \eqref{def:RK}(i)),
 and
 $(\rho_k(.))_k$ is any set of bounded functions such that $\rho_k(x)\geq \eta$, 
 $\forall x\in \mO_{k+1}$, $\forall k=0,\dots,N-1$, for some $\eta>0$.

    A useful example is the case of $\mO_k:=\mathscr{B}(0,L_0 + c k \dt\|f\|_\infty)$, $\forall k\geq 0$,
    with uniform densities $\rho_k$ compactly supported on $\mO_k$.
    In that case we notice that 
    $C_{k,\dt}=\frac{|\mO_{k+1}|}{|\mO_k|}
              = \big(\frac{L_0 + c (k+1) \dt\|f\|_\infty}{L_0 + c k \dt\|f\|_\infty}\big)^d$
    and also the following uniform estimate holds:
    \be\label{eq:useful-estim}
      \max_{0\leq k<N} \prod_{k=n}^{N-1} C_{k,\dt} 
              \leq \bigg(\frac{L_0 + c T\|f\|_\infty}{L_0}\bigg)^d.
    \ee 

 \item
  case $\mO_k=\R^d$ and $\rho_k(x)\stackrel{|x|\conv\infty}{\sim} e^{-q_k |x|}$ with $q_k>0$, $\forall k$.

\item
  case $\mO_k=\mO$ and $\rho_k=\rho$, $\forall k$, where $\mO$ is bounded,
assuming furthermore that $\mO$ is invariant by the dynamics, i.e.,
    $F_h^a(\mO)\subset \mO$ for all $a\in A$ and $h\geq 0$.
\end{itemize}


We can now give the following equivalent properties for $V_n$.

\begin{prop}
  Let $n\in \llbracket 0,N-1\rrbracket$ and $(\mO_n,\rho_n)$ as in \hypfive\ (with associated random variables $X_n$). Then
  $V_n$ satisfies the following dynamic programming principle
  \be\label{eq:DPP-2-vn} 
    V_n(x) =  G^{\ba_n}(x) \bigvee V_{n+1}(F^{\ba_n}(x)), \quad \forall x\in \mO_n 
  \ee
  for any 
  \be\label{eq:DPP-2-an}
    \ba_n(.) \in \argmin_{a\in \cA} \E\bigg[G^a(X_n) \vee V_{n+1} (F^a(X_n))\bigg].
  \ee
  In particular, we have 
  \be \label{eq:dpp-useful}
     \E[V_n(X_n)] =  \E\bigg[ G^{\ba_n}(X_n) \bigvee V_{n+1}(F^{\ba_n}(X_n))\bigg]
                  = \inf_{a\in \cA} \E\bigg[G^a(X_n) \vee V_{n+1} (F^a(X_n))\bigg].
  \ee

\end{prop}

\begin{proof}   
  The proof follows from Lemma~\ref{lem:equiv-minimization}
  and the dynamic programming principle of Proposition~\ref{prop:dpp-Vn-cA}.
\end{proof}

The above reformulation with an averaging criteria is motivated by numerical aspects: 
the problem can then be relaxed with an approximation $\numA$ of the control space $\cA$, 
for instance by neural networks. However, in general, $\ba_n$ is no more than measurable.
To circumvent this difficulty,
we first approximate problem \eqref{eq:dpp-useful} by more regular feedback controls.

\section{Approximation by Lipschitz continuous feedback controls}
\label{sec:LIP-FEEDBACK}

We aim to approximate \eqref{eq:dpp-useful} by using by Lipschitz continuous feedback controls.
Note that in Krylov~\cite{kry-80}, 
some approximations using feedback controls are given, yet in a different
context with stochastic differential equations and for non-degenerate diffusions.

Let $\rho \in \mathcal{C}^1\left(\mathbb{R}^d,\mathbb{R}\right)$ 
be a smooth function such that supp$(\rho) \subset \mathscr{B}(0,1)$, 
and $\int_{\mathbb{R}^d} \rho(x)dx=1$. 
Let $(\rho_\eps)_{\eps>0}$ be the mollifying sequence such that $\rho_\eps(x) \coloneqq \frac{1}{\eps^d}\rho(\frac{x}{\eps})$. 
For any sequence $a \coloneqq \left(a_0,\dots,a_{N-1}\right) \in \cA^{N}$, 
we associate the regularization by convolution
\be\label{eq:a-eps-1}
  a_k^\eps \coloneqq \rho_{\eps} * a_k.
\ee
Therefore $a^\eps_k$ is Lipschitz continuous, 
and $\|\nabla a^\eps_k\|_{L^\infty} = \| (\nabla \rho_\eps) * a_k\|_{L^\infty} \leq  \frac{1}{\eps}\|\nabla \rho\|_{L^1} \|a_k\|_{L^\infty}$. 
By classical arguments,
$\lim_{\eps\conv 0} a^\eps(x)  = a(x)$ a.e $x \in \mathbb{R}^d$.


In this section the following assumptions on $F_h$ will be needed.

\medskip

\begin{em}
\noindent
{\bf \hypfourA} 
The function $F_h$ satisfies:
\begin{itemize}
  \item there exists a constant $C>0$ and $h_0>0$ such that, for all $0<h\leq h_0$:
     \begin{align}\label{eq:F-bound-sec4}
        |F_h (x,a)|\leq |x| + C h (1+|x|), \quad x\in \R^d, a\in A,
     \end{align}
  \item $F_h$ satisfies a continuity property (for all $0<h\leq h_0$):
      \begin{align}\label{eq:F-cont-sec4}
      \forall x \in \R^d, \quad a\in A \conv F_h(x,a)\in \R^d \ \mbox{is continuous.}
      \end{align}
\end{itemize}
\end{em}
Such assumptions are naturally satisfied by Euler or Heun schemes already mentioned, 
and will be satisfied by more general RK schemes.

The following result will be used later on 
in order to obtain a regularized sequence
of controls (for the approximation of the dynamic programming principle) which will be more and more precise as $k$ varies
from $k=n$ to $k=N-1$.

\begin{prop}\label{prop:eps-weak-approx-refined}
  Let $k \in \llbracket0,N\rrbracket$.
  Assume (H0)-(H4) and \eqref{eq:F-bound-sec4}-\eqref{eq:F-cont-sec4}.
  Then 
  \be \label{eq:dpp-error-regularization}
    \lim_{\eps\conv 0^+}  
    \bigg|
    \E [ V_{k}(X_k) ]  - 
    \E [ G^{\ba_k^{\eps}}(X_k) \bigvee V_{k+1}(F(X_k,\ba_k^{\eps})) ] 
    \bigg| = 0
  \ee
  (with $\ba_k$ as in \eqref{eq:DPP-2-an} and $\ba_k^\eps$ as in \eqref{eq:a-eps-1}).
\end{prop}

\begin{lem}\label{lem:prelim}
  Assume \eqref{eq:F-bound-sec4}. There exists constant $\ma,\mb$ independent of $p\geq 1$, such  that 
$$
  \forall (x,a)\in \R^d\times A,\quad  |F(x,a)|\leq \ma |x| + \mb.
$$
\end{lem}
\begin{proof}
  By using the bound $|F_h(x,a)|\leq |x|+Ch(1+|x|) \leq e^{Ch}|x| + Ch$,
  by recursion (discrete Gronwall estimates)
  we obtain $|F(x,a)| = |F_h^{(p)}(.,a)| \leq e^{Chp} (|x|+ pCh) \leq e^{C\dt} (|x|+C\dt)$.
\end{proof}

\begin{proof}[Proof of Proposition~\ref{prop:eps-weak-approx-refined}]
In order to simplify the presentation, we consider the case of $p=1$ ($G^a(x)=g(x)$), the proof being similar 
in the general case $p\geq 1$.
The optimal control $\ba_n(.)$ satisfies
$V_n(x)=g(x) \vee V_{n+1}(F(x,\ba_n(x)))$ a.e. $x\in \mO_n$, hence
$\E[V_n(X_n)] = \E[ g(X_n)\vee V_{n+1}(F(X_n,\ba_n))]$.
  By assumption \eqref{eq:F-cont-sec4} of \hypfourA, and the pointwise convergence of $\ba_n^\eps(x)$,
we deduce that $\lim_{\eps\conv 0} F(x,\ba_n^\eps(x)) = F(x,\ba_n(x))$ a.e.  
  On the other hand, by using Lemma~\ref{lem:prelim} and the fact that $g$ is Lipschitz continuous, 
  there exists constants $\ma',\mb'$ such that $|g(x)\vee F(x,\ba_n^\eps(x))|\leq \ma' |x| + \mb'$ $\forall x\in \R^d$.
Therefore, the result is obtained by Lebesgue's dominated convergence theorem,
  the continuity of $(x,y)\conv g(x)\vee V_{n+1}(y)$ and the integrability assumption on $X_n$.
\end{proof}

We give also an other approximation result of $V_n$ by $J_n$.


\begin{prop}\label{prop:eps-weak-approx}
  Let $N\geq 1$ be given, assume (H0)-(H4), and \eqref{eq:F-bound-sec4}-\eqref{eq:F-cont-sec4}. Then
  \begin{align*}
    \lim_{\eps\conv 0} \max_{0\leq n \leq N-1} \E\Big[|J_n\big(X_n, (\ba_{n}^\eps,\dots,\ba_{N-1}^{\eps}) \big) 
      - V_n(X_n)|\Big] = 0.
  \end{align*}
\end{prop}

\begin{proof}
  It suffices to prove the result for a given $n$.
Let $a \in \cA^{N-n}$ be an arbitrary sequence. 
  By Lemma~\ref{lem:prelim} we have
$$ 
  |X_{k+1,x}^a| \leq e^{C\dt} (|X_{k,x}^a| + C \dt).
$$
Then by similar estimates, we obtain
\begin{align*}
  |X_{n,x}^a| 
     \leq e^{C n\dt} (|x| + Cn\dt) 
     \leq e^{C T} (|x| + C T).
\end{align*}
Hence
\begin{align*}
        \left|J_n(x,a)\right| 
        \leq \max_{n\leq k \leq N} \left[|g(0)| + [g]\left|X_{k,x}^{a}\right|\right] \bigvee \left[|\varphi(0)| + [\varphi]\left|X_{N,x}^a\right|\right]
        \leq K_0 + K_1 (|x| + C T)
\end{align*}
where $K_0 \coloneqq |g(0)|\vee|\varphi(0)|$ and $K_1 \coloneqq ([g]\vee[\varphi])e^{C T}$. 
  In particular
$\E\left(\left|J_n(X,a)\right|\right) <\infty$.
  In the same way, we can also obtain, for any $k\leq n \in[0,N]$,
  a Lipschitz bound of the form 
$$
   |G^a(x) \vee G^a(X^a_{k+1,x}) ... \vee G^a(X^a_{n,x}) \vee V_{n+1}(X^a_{n+1,x}))| \leq \ma' |x| + \mb'.
$$
  (for some constant $\ma',\mb'$).

In order to simplify the presentation, we consider again the case $p=1$ ($G^a(x)=g(x)$), the proof being similar 
in the general case $p\geq 1$.
Let $\eta >0$.
Consider the optimal control sequence $\ba$ and its regularization $\ba^\eps$. 
Let $n\in [0,\dots,N-1]$. By using the optimality of $\ba_n$, 
we have $V_n(x)= g(x) \vee V_{n+1}(F^{\ba_n}(x))$,
and as in Proposition~\ref{prop:eps-weak-approx-refined},
for $\eps>0$ small enough,
$$
  \big| \E[V_n(X_n)]  - \E[ g(X_n)\vee V_{n+1}(F^{\ba_n^\eps}(X_n)] \big|  \leq  \frac{\eta}{N}.
$$
Then we remark that $V_{n+1}(F^{\ba_n^\eps}(x))= g(F^{\ba_n^\eps}(x)) \vee V_{n+2}(F^{\ba_{n+1}}( F^{\ba_{n}^\eps}(x)) )$.
Hence by the same argument as before, for $\eps>0$ small enough, 
$$
  \bigg| \E[V_n(X_n)]  - \E[ g(X_n)\vee g(F^{\ba_n^\eps}(X_n)) \vee V_{n+2}(F^{\ba_{n+1}^\eps}(F^{\ba_{n}^\eps}(X_n)) )]
   \bigg| \leq  2\frac{\eta}{N}.
$$
Iterating this argument we deduce the existence of Lipschitz continuous controls 
$\ba^\eps:=(\ba_n^{\eps},\dots,\ba_{N-1}^\eps)$ 
such that 
\be \label{eqlem:2}
  \bigg| \E[V_n(X_n)]  - \E[ J_n(X_n,\ba^\eps)] \bigg| \leq  (N-n) \frac{\eta}{N} \ \leq \eta.
\ee
This concludes the proof.
\end{proof}



\section{Numerical schemes}

\subsection{Dynamic programming schemes}

%

It is natural to consider approximation schemes that mimic the dynamic programming principle
\eqref{eq:DPP-2-vn} - \eqref{eq:DPP-2-an}.
We see that \eqref{eq:DPP-1} (or \eqref{eq:DPP-1-vn-feedback})
has been relaxed by \eqref{eq:DPP-2-an} by minimizing 
a certain expectation over a set of feedback controls. 

We consider three schemes. Two of them may be seen as deterministic counterparts of 
the "value iteration" scheme (Huré \textit{et al.} \cite{hur-pha-21-a} or the BSDE scheme of \cite{gobet2005regression})
and the "performance iteration" scheme (Huré \textit{et al.} \cite{hur-pha-21-a}
or the BSDE scheme of \cite{bender2007forward}), hereafter denoted
the "SL-scheme" and the "L-scheme", respectively.

The third one is an hybrid combination of both paradigm (hereafter denoted the "H-scheme").

The set of measurable functions $\cA$ will be approximated by finite-dimensional
spaces $(\numA_{n})_n$, with $\numA_n$ typically a neural network space. 
When needed, neural networks will also be used in order to approximate value functions: in this case, 
we will denote by $\numV_n \subset \mathcal{C}\left(\mO, \mathbb{R}\right)$ 
a finite-dimensional space for the approximation of $\hV_n$.

Let $(\rho_n)_{0\leq n\leq N}$ be a sequence of densities supported 
in domains $(\Omega_n)_{0\leq n\leq N}$, as in \hypfive, 
with associated random variables $(X_n)_{0\leq n\leq N}$.
Recall that, for feedback controls $a\in \cA$,
$F^a(x)$ and $G^a(x)$ are defined at the beginning of section~\ref{sec:feedback} in terms of the approximate 
dynamics $F_h(.,.)$ ($F^a(x)$ corresponds to $p$ iterates of $y\conv F_h(y,a(x))$ starting from $y=x$,
and $G^a(x)$ corresponds to the maximum of $g(.)$ taken at the first previous $p$ iterates).

\medskip


\paragraph{{Semi-Lagrangian scheme (or "SL-scheme")}}
Let $(\numA_n)_{n\in\llbracket0,N-1\rrbracket}$ and $(\numV_n)_{n\in\llbracket0,N-1\rrbracket}$ 
be two given sequences of finite-dimensional spaces. 
Set $\hat{V}_N \coloneqq g \vee \varphi$. 
Then, for $n = N-1,\dots,0$: 
\begin{quote}
 - compute a feedback control $\ha_n$  according to
  \begin{subequations}
  \begin{align}\label{eq:SL-scheme-a}
  \ha_n \in \argmin_{a\in\numA_n} \E\bigg[  G^a(X_n) \bigvee \hV_{n+1} (F^a(X_n))  \bigg]
  \end{align}

 - set 
  \begin{align}\label{eq:SL-scheme-b}
    \hV_n := \argmin_{V\in\numV_n}
       \E\bigg[ \left| V(X_n) - G^{\ha_n}(X_n) \bigvee \hV_{n+1}\big(F^{\ha_n}(X_n)\big)\right|^2 \bigg].
  \end{align}
  \end{subequations}
\end{quote}

The approximations $(\hV_n)_{n\in\llbracket0,N-1\rrbracket}$ 
are stored, and only $\hV_{n+1}$ is used at iteration $n$. 
This explains the "semi-Lagrangian" terminology.
Owing to these projections, the computational cost is in $O(N)$, where $N$ is the number of time steps. 

\medskip
\noindent
\paragraph{{Lagrangian scheme (or "L-scheme")}}
Let $(\numA_n)_{n\in\llbracket0,N-1\rrbracket}$ be a given sequence of finite-dimensional spaces.
Set $\hat{V}_N \coloneqq g \vee \varphi$. 
Then, for $n = N-1,\dots,0$:
\begin{quote}
 - compute a feedback control $\ha_n$ according to
  \begin{subequations} \label{eq:L-scheme}
  \begin{align}\label{eq:L-scheme-a}
  \ha_n \in \argmin_{a\in\numA_n} \E\bigg[  G^a(X_n) \bigvee \hV_{n+1} (F^a(X_n))  \bigg]
  \end{align}

 - set 
  \begin{align}\label{eq:L-scheme-b}
  \hat{V}_{n} (x) \coloneqq G^{\ha_n}(x) \bigvee \hat{V}_{n+1}\left(F^{\ha_n}(x)\right)
      \equiv J_n(x,(\ha_n,\dots,\ha_{N-1})).
  \end{align}
  \end{subequations}
\end{quote}

In this algorithm, only the feedback controls $(\ha_k)$ are stored ($\hat{V}_n$ is not stored).
Each evaluation of the value $\hat{V}_{n+1}(x)$ 
uses the previous controls $(\hat{a}_{n+1},\dots,\hat{a}_{N-1})$ to compute the approximated characteristic, 
in a full Lagrangian philosophy. Therefore the overall computational cost is then of order $O(N^2)$. 
This scheme completely avoids projections of the value on functional subspace.

\begin{rem}
    From a computational point of view, an approximation of the minimum \eqref{eq:L-scheme-a}
  is obtained by using a stochastic gradient algorithm
  (see numerical section for details). 
  Hence  
  the optimality of $\ha_k$ in \eqref{eq:L-scheme-a} should therefore
  be replaced by some approximation $\ha_k$ such that
  \begin{align}\label{eq:L-scheme-a-approx}
    \E\bigg[  G^{\ha_k}(X_k) \bigvee \hV_{n+1} (F^{\ha_k}(X_k))  \bigg]
     \leq \min_{a\in\numA_k} 
     \E\bigg[  G^a(X_k) \bigvee \hV_{n+1} (F^a(X_k))  \bigg]
     + \mg_k
  \end{align}
  for some $\mg_k\geq 0$ (which takes into account some error on the optimal feedback control).
  Then an error analysis still holds (see Corollary~\ref{cor:a-approx}) showing some 
  robustness of the approach.
\end{rem}

\paragraph{{Hybrid scheme (or "H-scheme")}}
\newcommand{\Vtemp}{\hV^{[tmp]}}
Let $(\numA_n)_{n\in\llbracket0,N-1\rrbracket}$ and $(\numV_n)_{n\in\llbracket0,N-1\rrbracket}$ be two given sequences of finite-dimensional spaces. 
Set $\hat{V}_N \coloneqq g \vee \varphi$ as well as 
$\Vtemp_N \coloneqq g \vee \varphi$.
Then, for $n = N-1,\dots,0$:
\begin{quote}
 - compute a feedback control $\ha_n$  according to
  \begin{subequations}
  \begin{align}\label{eq:H-scheme-a}
    \ha_n \in \argmin_{a\in\numA_n} \E\bigg[  G^a(X_n) \bigvee \Vtemp_{n+1} (F^a(X_n))  \bigg]
  \end{align}

 - if $n\geq 1$, compute 
    $\Vtemp_{n} \in \numV_n$ (in prevision of the computation of $\ha_{n-1}$),
    such that
  \begin{align}\label{eq:H-scheme-c}
    \Vtemp_{n} \coloneqq \argmin_{V\in\numV_n} 
      \E\big[\left(V(X_{n}) - \hV_n(X_{n})\right)^2\big],
  \end{align}
  where $\hV_n$ is such that
  \begin{align}\label{eq:H-scheme-b}
    \hV_{n}(x) \coloneqq G^{\ha_n}(x) \bigvee \hV_{n+1}\big(F^{\ha_n}(x)\big) \ 
      \equiv J_n(x,(\ha_n,\dots,\ha_{N-1})).
  \end{align}   
  \end{subequations}
\end{quote}

The sequence of controls $(\ha_0,\dots,\ha_{N-1})$ is the output of the algorithm,
and $\hV_n(x)$ can be recovered using \eqref{eq:H-scheme-b}.
At each iteration $1\leq n\leq N$, $\hat{V}_n$ is projected on the space $\numV_n$, 
and its projection $\Vtemp_n$ is used to compute $\ha_{n-1}$. 
In this hybrid method, we still avoid some of the projection errors, 
by computing $\hV_n$ from the feedback controls in~\eqref{eq:H-scheme-b}.
Each evaluation of $\hV_n$ costs $N-n$ evaluations of the controls mappings, leading to an overall quadratic cost in $O(N^2)$. 
However, in the minimization procedure for \eqref{eq:H-scheme-a}, 
we can directly access to the values of $\Vtemp_n(.)$, which is less costly that computing $\hV_n(.)$.

In the present work, only the convergence of the L-scheme is analyzed.
However, the three proposed schemes will be compared 
on several examples in the numerical section (see in particular Sec.~\ref{ex:rot}).

\subsection {Runge Kutta schemes}
\label{sec:RK}

In this section, we consider a particular class of Runge-Kutta (RK) schemes for the definition of $F_h(x,a)$
which will corresponds to some approximation of the characteristics for a given control $a$.

For given $c=(c_i)_{1\leq j\leq q}$, $B=\in \R^{q\times q}$, let us denote 
$$ |c|_1:=\sum_{j=1}^q |c_i|, \quad \|B\|_{\infty} := \max_i \sum_{j} |b_{ij}| $$
and let also
$$ 
  C_f:=|f(0,A)| \coloneqq \max_{a\in A} |f(0,a)|.
$$

\begin{defi}[Runge-Kutta scheme]
\label{def:RK}
  \begin{enumerate}[label=(\roman*)]
    \item \label{eq:RKdef} 
      For a given $a\in A$, we say that
    $x \to F_h(x,a)$ is a Runge-Kutta scheme for $\dot{y} = f(y,a)$ 
    with time step $h>0$, if there exists 
    $q\in\mathbb{N}^*$, $\left(b_{ij}\right)_{i,j} \in \mathbb{R}^{q\times q}$ and $(c_i)_i \in \mathbb{R}^q$ such that
    \beno
      y_i(x)  =  x + h \sum_{j=1}^q b_{ij} f(y_j(x),a)\quad \forall i \in \llbracket1,q\rrbracket 
    \eeno
      and
    \beno
      F_h(x,a)  =   x + h \sum_{i=1}^q c_i f(y_i(x),a).
    \eeno
  \item \label{eq:RKconsist} The scheme is said to be \emph{consistent} if $\sum_{i=1}^q c_i = 1$.
  \item The scheme is said to be \emph{explicit} if $b_{ij} = 0$ for all $j \geq i$.

\end{enumerate} 
\end{defi}

\begin{rem}

  (a) Note that a consistent RK scheme 
satisfies $F_h(x,a) = x + h f(x,a) + O(h^2)$ for $h$ sufficiently small. 
This is made precise in Lemma~\ref{lem:RK-is-F-consistent}.

  (b) Also, for any given value $a\in A$, $y_i:= x + h \sum_{i=1}^q b_{ij} f(y_j,a)$ can be solved by a fixed point argument in
$(\R^d)^q$, as soon as $h \|B\|_\infty [f]_1 <1$.
Hence for $h$ small enough such that $h \|B\|_\infty [f]_1 <1$,
the RK scheme is well defined.
Explicit RK schemes are always well defined.

   (c) Note that in the above definition, the value of the control $a_0=a(x)$ is frozen at the foot of the characteristic.
Hence we consider an RK approximation of $\dot y(t) = f(y(t),a(x))$  on $[t_n,t_{n+1}]$ with $y(t_n)=x$ and fixed $a(x)$,
  rather than an approximation of $\dot y(t) = f(y(t),a(y(t)))$.
\end{rem}

We now give some estimates that will be useful later on.

The following lemma gives an estimate between two trajectories led by different controls.
\begin{lem}
  \label{lem:estim-a-abar}
  Assume $h \|B\|_{\infty} [f]_1 \leq \frac{1}{2}$, and $(a,\ba)\in \cA^2$.

  $(i)$ 
  \be\label{eq:estim-a-abar-a}
    |F_h(x,a(x)) - F_h(y,\ba(y))| \leq 
       e^{2 h |c|_1 [f]_1} |x-y| + 2 h |c|_1 [f]_2 |a(x) - \ba(y)|.
  \ee

  $(ii)$ $\forall\ 0\leq j\leq p$,
  \be\label{eq:estim-a-abar-b}
    | (F_h^a)^{(j)}(x) - (F_h^\ba)^{(j)}(x) | \leq C_F \dt |a(x) - \ba(x)|
  \ee
  (where $(F_h^a)^{(j)}(x)$ corresponds to $j$ iterates of $y\conv F_h(y,a(x))$ starting from $y=x$),
  where $C_F$ is a constant independent of $\dt$ and such that 
  \be\label{eq:C_F} 
    2 |c|_1 [f]_2 e^{2\dt |c|_1 [f]_1} \leq C_F.
  \ee
  In particular, 
  denoting $F^a(x)=F(x,a(x))$, we have also
  \be\label{eq:estim-a-abar-c}
    | F^a(x) - F^\ba(x) | \leq C_F \dt |a(x) - \ba(x)|
  \ee
  (recall that $F^a(x)$ corresponds to $p$ iterates of $y\conv F_h(y,a(x))$ starting with $y=x$, 
  and $h=\dt/p$).

  $(iii)$
  More generally,
  \be\label{eq:estim-a-abar-b-general}
    | F^a(x) - F^\ba(y) | \leq e^{C_1\dt} (|x-y| +   C_2 \dt |a(x) - \ba(y)|)
  \ee
  where $C_1:=2 |c|_1 [f]_1$ and $C_2:=2 |c|_1 [f]_2$.

\end{lem}

\begin{proof}
  $(i)$
  From the definitions, denoting $a_0=a(x)$ and $\ba_0=\ba(y)$, we have 
  \be
    \big| F_h(x,a_0)- F_h(y,\ba_0) \big|
       & \leq & |x-y| + h \sum_{1\leq j\leq q} |c_j| |f(y_j^a(x), a_0)-f(y_j^\ba(y),\ba_0)| \nonumber \\
       & \leq & |x-y| + h |c|_1 [f]_1 \|Y^a-Y^{\ba}\|_\infty  + h |c|_1 [f]_2 |a_0 - \ba_0|
       \label{eq:inter1}
  \ee
  where the intermediate values of the RK schemes are denoted $Y^a=(y^a_j(x))_{1\leq j\leq q}$ and
  $Y^\ba=(y^\ba_j(y))_{1\leq j\leq q}$, and satisfy
  $Y^a = X + h B f(Y^a,a_0)$ and $Y^\ba = Y + h B f(Y^\ba,\ba_0)$  
  (with $X=(x,\dots,x)$ and $Y=(y,\dots,y)$).
  Hence we have also
  \beno
     \| Y^a - Y^\ba \|_\infty
       & \leq & |x-y| + h \|B\|_\infty [f]_1 \|Y^a-Y^\ba\|_\infty  + h \|B\|_\infty [f]_2 |a_0 - \ba_0|,
  \eeno
  from which we deduce, using the assumption of the present Lemma,
  \beno
     \| Y^a - Y^\ba \|_\infty
       & \leq & 2 \big( |x-y| + h \|B\|_\infty [f]_2 |a_0 - \ba_0| \big).
  \eeno
  Combining with \eqref{eq:inter1}, we obtain
  \beno
    \big| F_h(x,a_0)- F_h(y,\ba_0) \big|
       & \leq & e^{2 h |c|_1 [f]_1} |x-y|   + 2 h |c|_1 [f]_2 |a_0 - \ba_0|.
  \eeno

  $(ii)$-$(iii)$ From the previous bound, 
  denoting $e_j:= 
    \big| (F_h^{a_0})^{(j)}(x)- (F_h^{\ba_0})^{(j)}(y) \big|
    \equiv
    \big| Y^{a_0}_{j,x} - Y^{\ba_0}_{j,y}\big|$, we have
  \beno 
    e_{j+1} 
       & \leq & e^{2 h |c|_1 [f]_1} e_j + 2 h |c|_1 [f]_2 |a_0 - \ba_0|,
  \eeno 
  and, by recursion,
  $$ 
    e_{j} \leq e^{2 j h |c|_1 [f]_1} (e_0 + 2 j h |c|_1 [f]_2 |a_0 - \ba_0|), \quad 0\leq j \leq p.
  $$
  This concludes to $(iii)$, and also for $(ii)$ using $y=x$, $e_0=0$ and $jh\leq \dt$. 
\end{proof}


\begin{lem}
  \label{lem:RKmajoration}
  Assume $h \|B\|_{\infty} [f]_1 \leq \frac{1}{2}$
  and let $F_h$ be an RK scheme. 

  $(i)$ For $a\in A$ and $x\in \R^d$:
  \beno
     |F_h(x,a)| \leq |x| + 2 h |c|_1 (C_f\vee [f]_1) (|x|+1).
  \eeno

  $(ii)$
  Let $a$ in $\cA$ and denote $F^a_h(x):=F_h(x,a(x))$, then it holds
  \be
    \label{eq:p-large-enough}
    \quad [F_h^a - i_d] \leq  2  |c|_1 ([f]_1+ [f]_2[a])\, h.  
  \ee
\end{lem}
\begin{proof}
  We start by proving $(ii)$.
  As in the proof of Lemma~\ref{lem:estim-a-abar}$(ii)$, with $\ba=a$,
  we obtain 
  $$
    |(F^a(x,a(x))-F^a(y,a(y)))- (x-y)|\leq  h|c|_1 [f]_1 (2|x-y| + 2h\|B\|_\infty [f]_2[a]|x-y|)
     +  h|c|_1 [f]_2[a] |x-y|.
  $$
  Combining with the assumption that $2h \|B\|_{\infty} [f]_1 \leq \frac{1}{2}$,
  we obtain the desired bound.

  For the proof of $(i)$, for any $a$ in $A$ we have $[f^a]=[f]_1$ and
  $|F_h^a(x)|\leq |F_h^a(0)| + [F_h^a]|x| \leq |F_h^a(0)| + (1+2h|c|_1[f]_1])|x|$.
  By direct bounds we have also $|F_h^a(0)|\leq 2h |c|_1 C_f$, 
  from which we deduce the desired bound.
\end{proof}

\begin{lem}
  \label{lem:RK-is-F-consistent}
  Assume $h |B|_{\infty} [f]_1 \leq \frac{1}{2}$
  and let $F_h^a$ be a consistent RK scheme. 
  Then $F_h$ is consistent with $f$ in the following sense:
  \be\label{eq:F-consistent-bis}
      \exists C\geq 0,\ \exists h_0>0, \forall (x,a,h)\in \R^d\times A\times(0,h_0), \quad
      \quad
        |F_h(x,a) - (x + h f(x,a))|\leq C h^2 (1+|x|),
  \ee
\end{lem}
\begin{proof}
We use $|f(y_j,a)|\leq [f]_1 |y_j| + C_f$ where  $C_f= \max_{a\in A}|f(0,a)|$ to obtain in the RK scheme
  $|y_i-x| \leq   h \sum_{j} |b_{ij}| ([f]_1 |y_j| + C_f)$, hence 
  \be\label{eq:majo-inter-un} 
    \max_i |y_i-x| \leq h \|B\|_\infty [f]_1 (\max_i |y_i|) + h \|B\|_\infty C_f.
  \ee 
  By using the assumption we then get 
  $\max_i |y_i|  \leq 2 |x| +  2h \|B\|_\infty C_f$.
  Let $h\in ]0,h_0]$ with $h_0$ such that $h_0 \|B\|_\infty [f]_1 = \frac{1}{2}$.
  Using \eqref{eq:majo-inter-un} (and the fact that $h \|B\|_\infty [f]_1 \leq \frac{1}{2}$), we obtain
  \be\label{eq:majo-inter-deux} 
    \max_i |y_i-x| \leq h \|B\|_\infty [f]_1 |x|  + 2 h \|B\|_\infty C_f = Ch (1+|x|)
  \ee 
  for some constant $C\geq 0$.
  Then
  \beno
    F_h(x,a) & = & x + h \sum_{j=1}^q c_j f(y_j,a) 
     \ = \  x + h \sum_{j=1}^q c_j f(x + O(h(1+|x|)),a) \\
     & = &  x + h \sum_{j=1}^q c_j \big(f(x,a)  + O(h(1+|x|)) \big)\\
     & = &  x + h f(x,a) + O(h^2(1+|x|))
  \eeno
which is the desired result.
\end{proof}

\subsection {Neural network spaces}

Neural networks are functions build by compositions of other "simple" functions. 
They are widely used for their approximating capabilities, 
and are known to be dense in the class of continuous multivariate functions under mild hypotheses
(see or instance Lemma 16.1 of \cite{gyorfiDistributionFreeTheoryNonparametric2002}).
We restrict ourselves to so-called \emph{feedforward} neural networks, in the following sense. 

\begin{defi}[Feedforward neural network]
  Let $L\in\mathbb{N}^*$ (the \emph{number of layers}),
  and $(d_k)_{k\in\llbracket0,L\rrbracket} \subset\mathbb{N}^*$ be a sequence of dimensions. 
  A neural network is a function $\mathcal{R} : \mathbb{R}^{d_0} \mapsto \mathbb{R}^{d_L}$ of the form
    $$
        \mathcal{R}(x) = \sigma_{L} \circ \mathcal{L}_L \circ \cdots \circ \sigma_1 \circ \mathcal{L}_1 (x)
    $$
        where $\mathcal{L}_k : \mathbb{R}^{d_{k-1}} \mapsto \mathbb{R}^{d_k}$ is an affine transformation, 
        $\sigma_k : \mathbb{R}^{d_k} \mapsto \mathbb{R}^{d_k}$ is a nonlinear \emph{activation} function
        which acts coordinate by coordinate:
        $$ \ms_k(x) = (\bar \ms_k(x_1), \bar \ms_k(x_2),\dots, \bar\ms_{k}(x_{d_k})) $$
        for $x\in \R^{d_k}$ and for a certain $\bar\sigma_k : \mathbb{R} \mapsto \mathbb{R}$. 
\end{defi}
Each affine transformation is represented by a \emph{weight}
matrix $\omega_{k} \in \mathbb{R}^{d_k \times (d_{k-1}+1)}$,
with $\mathcal{L}_k (x) = \omega_k {\begin{pmatrix} x \\ 1 \end{pmatrix}}$.

Classical examples of activation functions include the sigmoid function $\bar\ms(x)=\frac{1}{1+e^{-x}}$,
the rectified linear unit (ReLU) $\bar\ms(x)=\max(0,x)$, etc... 
The last activation function, $\bar\ms_L$, may be set to the identity function, 
or $\ms_L$ may be a more complex function (so that $\ms_L(.)\in A$ for the control approximation)
depending on the example.

We may now define the sets $\numA$ and $\numV$ as (recall that $A\subset \R^\kappa$)
\begin{align*}
        \numA &\coloneqq \left\{\text{feedforward neural networks with }d_0 = d \text{ and }d_L = \kappa\right\}, \\
        \numV &\coloneqq \left\{\text{feedforward neural networks with }d_0 = d \text{ and }d_L = 1\right\}. 
\end{align*}
In our numerical examples, 
we will always choose ReLU for the inner layer activation functions. 
The last activation function $\ms_L$
will vary to fit the definition of $A$ for the given problem (see the numerical section). 
We also choose to set $d_1 = \cdots = d_{L-1} \eqqcolon N_e$, i.e., the same \emph{number of neurons} 
for each layer, to simplify the set of parameters. 


\section{Main result}

In this section, we focus on proving the convergence of the Lagrangian scheme \eqref{eq:L-scheme}.
This algorithm only uses approximations of feedback controls.


Since our estimates will need Lipschitz continuous controls,
and since the exact optimal solution  in general does not involve such regular controls,
we first introduce $\eta$-weak approximations as follows.

\begin{defi}\label{def:eps-weak}
  For a given sequence $\eta=(\eta_n,\eta_{n+1},\dots,\eta_{N-1})\in (\R_+^*)^{N-n}$, 
  we say $\ha=(\ha_n,\dots,\ha_{N-1})\in \cA^{N-n}$ is
  an $\eta$-weak approximation of $(V_k)_{n\leq k\leq N}$ if \\
$(i)$ $(\ha_n,\dots,\ha_{N-1})$ are Lipschitz continuous controls, \\
  $(ii)$ 
    $\bigg|\E[V_k(X_k)] -  \E[ G^{\ha_k}(X_k) \vee V_{k+1} (F^{\ha_k}(X_k))]\bigg|\leq \eta_k$, for all 
     $k\in \llbracket n, N-1\rrbracket$.
  
\end{defi}

Notice that by using Prop.~\ref{prop:eps-weak-approx-refined}, 
it is possible to construct $\eta$-weak approximations. By recursion, it is furthermore
possible to construct  the controls such that, for all 
$k\in \llbracket n, N-1\rrbracket$
\be  \label{eq:weak-controls-precised-estimate}
  \bigg|\E[V_k(X_k)] -  \E[ G^{\ha_k}(X_k) \vee V_{k+1} (F^{\ha_k}(X_k))]\bigg|\leq 
    \eta_k([\ba_n^{\eps_n}], [\ba_{n+1}^{\eps_{n+1}}], \dots, [\ba_{k-1}^{\eps_{k-1}}]),
\ee
for given strictly positive functions $\eta_k$ (i.e., $\eta_k>0$, $\eta_{k+1}(x_k)>0$, $\dots$, $\eta_{N-1}(x_k,\dots,x_{N-2})>0$).

We now state our last assumption that will be needed on the approximate dynamics $F_h$.
First assumptions where already introduced in \hypfourA. The new assumption is the following.
In what follows, we denote $f^a(x):=f(x,a(x))$ and recall that $[f^a]$ is the Lipschitz constant of $f^a$.

\medskip


\begin{em}
\noindent
{\bf \hypfourB} 
  %
    There exist a constant $\delta>0$, 
    for any Lipschitz continuous function $a(\cdot)\in \cA$ and $h>0$ such that 
    $h [f^a]\leq \delta$,
    $x\converge F_h(x,a(x))$ is 
    one-to-one and onto on $\R^d$, Lipschitz continuous,
    with Lipschitz bound
    \be\label{eq:H4-majo} 
      %
      %
      [F_h(.,a(.)) - i_d]\leq c\, ([f]_1+ [f]_2[a])\, h
    \ee
    (where $i_d(x):=x$),
    where
    $c\geq 0 $ is a universal constant (independent of $a,h$).
\end{em} 

Assumption \hypfourB\ will needed in order to use a change of variables formula (for $y=F^a(x)$, corresponding
to $p$ iterates of $y\conv F_h(y,a(x))$ starting from $y=x$).

\begin{rem}
  Lemma \ref{lem:RKmajoration} shows that \hypfourA-\hypfourB\ are satisfied 
  for Runge Kutta schemes as defined in~\ref{def:RK}.
  For instance, dynamics $F_h$ 
  such as the Euler scheme $F_h(x,a) = x + h f(x,a)$, or
  the Heun scheme  $F_h(x,a) = x + \frac{h}{2} (f(x,a) + f(x+ h f(x,a),a))$
  satisfy \hypfourA-\hypfourB.
  Higher order schemes, or implicit schemes, could be used as well.
\end{rem}

\begin{rem}
Assumption \hypfourB\ is satisfied by the exact characteristics. 
  Indeed, for any regular control $a \in \mathcal{A}$, the map $x \to F^a_h(x):=y_{x}^a(h)$
  (where $y(s)=y_x^a(s)$ is the solution of $\dot y(s)=f^a(y(s))$ with $y(0)=x$)
  is one-to-one and onto on $\R^d$ and satisfies 
  $D F_h^ a(x) = \exp B_h$ with $B_h=\int_0^h D f^ a(y_{x}^a(s))ds$.
  Then denoting $\|D F_h^a\|_\infty = \sup_{x\in \R^d} \|D F_h^a(x)\|_\infty$, we have
$[F_h^a-i_d]=\| DF_h^a - I\|_\infty \leq e^{\|B_h\|_\infty}-1 \leq 2 \|B_h\|_\infty$ as soon as for instance
$\|B_h\|_\infty\leq \frac{1}{2}$, with also 
  $\|B_h\|_\infty\leq h \|Df^a\|_\infty  = h [f^a]$, hence \eqref{eq:H4-majo} holds true 
  with $c=2$ and $\delta=\frac{1}{2}$.
When $a$ is only Lipschitz regular, the same bound is obtained by a regularization argument.
\end{rem}

\begin{rem}
        Assumption \eqref{eq:H4-majo} implies the following bounds:
        \be
          & & \|DF_h^a(x)\|_\infty \leq [F_h^a] \leq 1 + c [f^a]h \leq e^{c[f^a]h},
             \quad \mbox{a.e. $x\in \R^d$} \label{eq:DF-bound} 
        \ee
from which we can deduce
        \be
          & & |\det(DF_h^a(x))| \leq e^{d c [f^a] h}, \mbox{ a.e. $x\in \R^d$.} \label{eq:det-bound}
        \ee
        This estimate will be used in the change of variable Lemma~\ref{lem:cv}.
Note also that \hypfourB\ and
$\max_{a\in A} |F(0,a)|\leq Ch$ for some constant $C$,
implies the bound \eqref{eq:F-bound-sec4} of \hypfourA.
\end{rem}

\medskip


\begin{rem} \label{rem:hur-not-for-deterministic}
  In~\cite{hur-pha-21-a}, 
  an assumption on the controlled transition probability of a stochastic process 
  (say $x=X_k \conv x'=X^a_{k+1}$ for a given control $a$) is made,
  which is to be measure of the form
  $$
       r(x,a; x') \mu(dx'),
  $$
  for some measure $\mu$ which has a finite first order moment,
  and assuming a uniform bound $\|r\|_\infty<\infty$.
  However this assumption cannot be satisfied in our deterministic context 
  (where, typically, $x' = x + \dt f(x,a)$). 
  Instead, assumption \hypfourB\ will be used
   (in the change of variable Lemma~\ref{lem:cv})
   in order to get a recursive error bound estimate.
\end{rem}

In the following, 
we recall that $(\hV_k)_{0\leq k\leq N}$ corresponds to the Lagrangian scheme \eqref{eq:L-scheme}.
Also we have $\hV_k(x)\geq V_k(x)$, and therefore 
we look for an upper bound of $\hV_k(x)-V_k(x)$.

\begin{thm}\label{thm:1}
  Assume (H0)-\hypfourB, $N \geq 1$, and let $n\in \llbracket 0, N \rrbracket$.
For $k=n,\dots,N-1$, let $\eta_k:\R^{n-k}\conv \R_+^*$ be given functions
($\eta_n>0$ is a constant, $\eta_{n+1}>0$ is a function of one variable, and so on).
  Let $\ba=(\ba_0,\dots,\ba_{N-1})$ be an $\eta$-weak approximation
  in the sense of~\eqref{eq:weak-controls-precised-estimate}.
Then
\beno
   & & \hspace{-1cm} \E\left[(\hV_n - V_n)(X_n)\right]  \\
   & \leq &
   \inf_{(a_n\dots,a_{N-1}) \in \bigotimes_{k=n}^{N-1}\numA_{k}}
   \bigg( (\eps_n^{a_n}+\eta_n) + C_n^{a_n} (\eps_{n+1}^{a_{n+1}}+\eta_{n+1}([a_n])
   \\
   &      & \hspace{4cm}
   + \dots + C_k^{a_k}\cdots C_{N-2}^{a_{N-2}} (\eps_{N-1}^{a_{N-1}}+\eta_{N-1}([a_n],\dots,[a_{N-2}])) \bigg)
\eeno
  where $C_k^a:= C_{k,\dt} e^{dc[f^a]\dt}$
  (with $c\geq 0$ is as in \hypfourB, and $C_{k,\dt}$ as in \hypfive), 
\begin{align}\label{eq:defEpsilonak}
  \eps^{a}_k :=C_F  ([g]+[V_{k+1}])\, \dt \, \E_k\Big[|a(X_{k}) - \ba(X_k)|\Big],
\end{align}
  and $[V_{k+1}]$ is bounded as in Lemma~\ref{lem:vn-lip-bound}, $C_F$ satisfies \eqref{eq:C_F}.
\end{thm}

Note that the consistency of the scheme (with respect to the dynamics $f$, as in \eqref{eq:F-consistent}) 
is not needed in the previous Theorem, because 
the result only focuses on the error between the semi-discrete problem and its approximation by a Lagrangian scheme.

\begin{corol}\label{cor:a-approx}
  In the same way, for the perturbed algorithm \eqref{eq:L-scheme-a-approx}
  the same error bound holds where each term  $(\eps_k^{a_k}+\eta_k)$ is replaced by 
  $(\eps_k^{a_k}+ \eta_k  + \mg_k)$. 
\end{corol}

The proof of Theorem~\ref{thm:1} is postponed to section~\ref{sec:thm-proof}
(the proof of Corollary~\ref{cor:a-approx} follows exactly the same lines).
We now give two corollaries of the previous theorem.

\begin{corol}\label{cor:1}
  Assume (H0)-\hypfourB, and $N\geq 1$. 
Let $\numA_n^{\Theta}$ denote the control approximation space at time $t_n$,
with explicit dependency over the size $\Theta$.
We denote by $\Theta\to \infty$ the limit when some parameters go to infinity
  (for instance the number of neurons of a neural network). 
We assume that for any $n=0,\dots,N-1$, any Lipschitz continuous function $\bar a\in \cA$ can be approximated 
  by some function of $a\in \numA_n^\Theta$ up to any arbitrary precision,
  which we write as 
  \begin{align}\label{eq:AM-approx}
    \lim_{\Theta\conv\infty} \inf_{a\in \numA_n^{\Theta}}
      \E[|a(X_{n})-\ba(X_{n})|]=0.
  \end{align}
  Let $(\hV^\Theta_n)$ be the corresponding $L$-scheme values associated with sets 
    $(\cA_n^{\Theta_n},\cA_{n+1}^{\mT_{n+1}},\dots,\cA_{N-1}^{\mT_{N-1}})$.
Then 
$$
    \lim_{\Theta\converge\infty} \max_{0\leq n\leq N} \E[\hV_n^\Theta(X_n)-V_n(X_n)] = 0.
$$
  (where $\mT\conv \infty$ means here that $\mT_k\conv \infty$ for all $k=n,\dots,N-1$).
\end{corol}

\begin{rem} \label{rem:groupsort}
  Notice that Group Sort neural networks satisfies furthermore 
  (see \cite{anil2019sorting,tanielian2021approximating}):
  \begin{align}\label{eq:AM-approx-groupsort}
    \lim_{\Theta\conv\infty} \inf_{a\in \numA_n^{\Theta},\ [a]\leq [\bar a]}
      \E[|a(X_{n})-\ba(X_{n})|]=0.
  \end{align}
\end{rem}

\begin{proof}[Proof of Corollary~\ref{cor:1}]
  Let $\eps>0$.
  Let $\eta_n:=\eps/(2N)$.
  By assumption \eqref{eq:AM-approx}, 
  $$\lim\limits_{\mT\conv\infty} \inf\limits_{a_k\in \numA_k^\mT} 
     \E_{k+1}[|a_k - \ha_k|] =0 \quad\quad \forall n \leq k \leq N-1 $$
  and therefore 
  $\lim\limits_{\mT\conv\infty} \inf\limits_{a_k\in \numA_k^\mT} \eps_k^{a_k} =0$, 
  where $\eps_k^{a}$ is defined in \eqref{eq:defEpsilonak}.
  Hence we can find $a_n\in \numA_{n}^{\mT_n}$ (for $\mT_n$, the size of $\numA_{n}$, large enough) such that  
  $\eps_n^{a_n}\leq \frac{\eps}{2N}$, and therefore
  $$
    \eps_n^{a_n}+\eta_n\leq \frac{\eps}{N}.
  $$

  Then let $\eta_{n+1}([a_n]):= \eps/(2N C_n^{a_n})$.
  There exists $a_{n+1}\in \numA_{n+1}^{\mT_{n+1}}$ 
  (for $\Theta_{n+1}$ large enough) 
  such that 
  $$
    C_n^{a_n}(\eps_{n+1}^{a_{n+1}}+\eta_{n+1}([a_n])) \leq \frac{\eps}{N},
  $$
  and so on, until we chose 
  $\eta_{N-1}([a_{n}],\dots,[a_{N-2}])
  := \eps/(2N C_n^{a_n}C_{n+1}^{a_{n+1}}\cdots C_{N-2}^{a_{N-2}})$
  and then find $a_{N-1}\in \numA_{N-1}^{\mT_{N-1}}$ such that 
  $$
    C_n^{a_n}C_{n+1}^{a_{n+1}}\cdots C_{N-2}^{a_{N-2}} \big(\eps_{N-1}^{a_{N-1}} 
    +   \eta_{N-1}([a_{n}],\dots,[a_{N-2}]) \big) \leq \frac{\eps}{N}
  $$
  (we have $N-n$ such bounds).
  By using the bound of Theorem~\ref{thm:1}, 
  the sum of all error terms is bounded by $ (N-n) \frac{\eps}{N}$, and
  therefore
\beno
   \E[\hV_n(X_n) - V_n(X_n)] 
   & \leq & 
   \eps.
\eeno
  This shows that $\lim\limits_{\mT\conv\infty} \E[\hV^{\mT}_n(X_n) - V_n(X_n)] =0$.
  The desired result follows since we have only a finite number $N$ of such terms.
\end{proof}

Notice that in the previous result $N\geq 1$ is given. 
This does not give in general a convergence result
as $N\conv\infty$, because of the uncontrolled Lipschitz constants that appear in the bounds of 
Theorem~\ref{thm:1}.

However, in the case the optimal controls $(\ba_n)_n$ can be shown to be Lipschitz continuous
with a uniform Lipschitz constant, we may improve the result. 
We suppose that the numerical feedback space $\numA$ can be restricted to Lipschitz functions
with a controlled Lipschitz constant. For instance, if $\numA$ is a neural network space, one could choose
the GroupSort activation function and bound the weights to obtain this estimate 
(see \cite{anil2019sorting,tanielian2021approximating}).

\begin{corol} \label{cor:2}
  Assume (H0)-\hypfourB, $N\geq 1$,
  and that there exists a sequence of optimal feedback control (denoted $\ba$)
  which are Lipschitz continuous: 
  $\exists L\geq 0$, $\forall 0\leq k\leq N-1$, $[\ba_k]\leq L$.
  Then 
  \begin{align*}
    \max_{0\leq n < N}\E[\hV_n(X_n) - V_n(X_n)] 
       \leq 
       K_N \ \inf_{(a_0\dots,a_{N-1}) \in \bigotimes_{k=0}^{N-1}\numA_{k},
       \ [a_k]\leq [\ba_k]\,\forall k}
       \bigg(\dt \sum_{k=0}^{N-1} 
         \E_k\big[|a_k(X_k) - \ba_k(X_k)|\big]
       \bigg),
  \end{align*}
where 
\be\label{eq:K_N}
  K_N := 2 C_F ([g]\vee [\varphi]) e^{(d+1) c ([f]_1 + L [f]_2) T} 
       \max_{0\leq n<N} \prod_{k=n}^{N-1} C_{k,\dt}.
\ee
  Furthermore, in the case of uniform densities, we can use the estimate \eqref{eq:useful-estim} 
  to deduce a bound for $K_N$ which is independent of $N$ (other situations could also lead 
  to a uniform bound for $K_N$).
\end{corol}

\begin{proof}[Proof of Corollary~\ref{cor:2}]
  We make use of the bound of Theorem~\ref{thm:1} with $\eta_k=0$, $\forall k$.
  Notice that $[f^{\ba_k}] \leq [f]_1 + L [f]_2$, and also, with $[a_k]\leq [\ba_k]\leq L$,
  we have $[f^{a_k}] \leq [f]_1 + L [f]_2$.
  Then 
  $$
   \prod_{n\leq k\leq N-1} C^{\ba_k} 
   \leq \big(\prod_{n\leq k\leq N-1} C_{k,\dt}\big) e^{d c \sum_{k=n}^{N-1} [f^{\ba_k}] \dt}
   \leq \big(\prod_{n\leq k\leq N-1} C_{k,\dt}\big) e^{d c ([f]_1+L[f]_2) T}
  $$
  as well as 
  $[V_{k+1}]\leq ([g] \vee [\varphi]) e^{c ([f]_1+L[f]_2) T}$
  and therefore (using also $[g]\leq [g]\vee [\varphi]$)
  $$ \eps_k^{a_k} \leq  2 C_F ([g] \vee [\varphi]) e^{c ([f]_1+L[f]_2) T } 
      \dt\ \E_k \big[|a_k(X_k) - \ba_k(X_k)|\big].
  $$ 
  The desired result follows.
\end{proof}



%

\section{Proof of Theorem~\ref{thm:1}}
\label{sec:thm-proof}
We first state a change of variable Lemma,
giving a statement for either an exact characteristic 
($x\conv y_{t,x}^a$) or for an approximate one ($x\conv F^a(x)$).
Only the second statement will be used in the convergence analysis.

\begin{lem}[Change of variable] \label{lem:cv}
  Let $a:\R^d\converge \R^\mk$ be a given Lipschitz continuous function.
  \begin{enumerate}[label=(\roman*)]
  \item 
  Let $t\converge y^a_{t,x}$ denotes the characteristic associated 
  with dynamics~$x\converge f^a(x)$ and such that $y^a_{0,x}=x$.    
  We assume the following analogue of \hypfive\ in the continuous case:
  $$
    \mbox{$y^a_{\dt,\mO_k}\subset \mO_{k+1}$, $\forall k=0,\dots,N-1$},
  $$
  and
  $$ 
      \widetilde{C}_{k,\dt}:=\max\limits_{0\leq k\leq N-1} 
     \sup_{x\in \mO_k} 
      \frac{\rho_k(x)}{\rho_{k+1}(y^a_{\dt,x})}<\infty .
  $$
  Then
  for any non-negative measurable function $\Phi:\mO_{k+1}\converge \R$,
  \begin{align}\label{eq:estim-L-1}
    \E_k[\Phi(y^a_{\dt,X_k})] \leq \widetilde{C}_{k,\dt}
      \, e^{d [f^a]\dt} \ \E_{k+1}[\Phi(X_{k+1}) ].
  \end{align}
  
  \item  Suppose \hypfourA, \hypfourB\ and \hypfive.
    Assume $p\geq 1$ is such that $h[f^a]\leq \delta$, where $h=\frac{\dt}{p}$.
    Then for any non-negative measurable function $\Phi:\mO_{k+1}\converge \R$,
  \be \label{eq:estim-L-2}
      \E_k [\Phi(F^a(X_k)) ] \leq C_{k,\dt} \, e^{d c [f^a] \dt} \ \E_{k+1} [\Phi(X_{k+1}) ]
  \ee 
      where $c\geq 0$ is as in \hypfourB, $C_{k,\dt}$ 
      is as in \hypfive. 
  \end{enumerate}
\end{lem} 

\begin{proof}[Proof of Lemma~\ref{lem:cv}]
  $(i)$
  Let $y_{t,x}$ be the solution at time $t$ of the differential equation $\dot y_{t,x} = f^a(y_{t,x})$ for $t\in \R$,
  with $y_{0,x}=x$. Then for any function $\Phi\geq 0$ and $t\geq 0$,
  \beno
    \E_k [\Phi(y_{t,X_k})] = \int_{\mO_k} \Phi(y_{t,x}) \rho_k(x)\, dx \leq e^{d [f^a] t} \int_{y_{t,\mO_k}} \Phi(x') \rho_k(y_{-t,x'}) dx'
  \eeno
  (the Jacobian of the change of variable, as well as its inverse, 
  is bounded by $e^{d [f^a] t}$ for $t\geq 0$).
  Since the r.v. $X_k$ has density law $\rho_k(x)dx$, we deduce
  \be\label{eq:b2}
    \E_k[ \Phi(y_{t,X_k}) ] \leq 
    e^{d [f^a]t}  \int_{y_{t,\mO_k}} \Phi(x') \frac{\rho_k(y_{-t,x'})}{\rho_{k+1}(x')} 
    \rho_{k+1}(x') dx'.
  \ee
  Let $t=\dt$,
  we have $y^a_{t,\mO_k}\subset \mO_{k+1}$ by assumption \eqref{eq:H5-b}. 
  Also, for any $x'=y^a_{\dt,x} \in y^a_{\dt,\mO_k}$,
  we have $\rho_k(y^a_{-\dt,x'})/\rho_{k+1}(x')= \rho_k(x)/\rho_{k+1}(y^a_{\dt,x})
   \leq \widetilde{C}_{k,\dt}$ by assumption~\eqref{eq:H5-c}.
  Together with \eqref{eq:b2} this allows to conclude to the desired bound.

  $(ii)$ The proof is completely similar to $(i)$.
\end{proof}

We are now in position to prove the main result.

\begin{proof}[Proof of Theorem~\ref{thm:1}]
Our aim is to bound recursively the quantity
\begin{align*}
  e_n : = \E [\hV_n(X_n)-V_n(X_n)].  
\end{align*}
Let $\eta_n>0$. By Prop.~\ref{prop:eps-weak-approx-refined}, there exists $\ba_n\in\cA_n$, Lipschitz continuous,
such that 
$$
    \big| \E[V_n(X_n)]-\E[G^{\ba_n}\vee V_{n+1}(F^{\ba_n}(X_n))] \big| \leq  \eta_n.
$$
Recall that $\hV_n$ satisfies
\begin{align*}
  \E\left[\hV_n(X_n)\right]
  = \inf_{a\in \numA_n} 
    \E\left[ G^a(X_n) \bigvee \hV_{n+1}(F^{a}(X_n)) \right],
\end{align*}
hence 
\beno
  & & \hspace{-2cm} \E\left[\hV_n(X_n) - V_n(X_n)\right]  \\
  & & 
  \leq \inf_{a\in \numA_n} 
    \E\left[ G^a(X_n) \bigvee \hV_{n+1}(F^{a}(X_n))  \ - \
               G^{\ba_n}(X_n)\bigvee V_{n+1}(F^{\ba_n}(X_n))\right]  +  \eta_n.
\eeno

Thus, using $\max(a,b) - \max(c,d) \leq \max(a-c,b-d)$, we have
\beno
  & & \hspace{-2cm} G^a(x) \bigvee \hV_{n+1}(F^{a}(x)) - G^{\ba_n}(x) \bigvee V_{n+1}(F^{\ba_n}(x)) \\
  & \leq & \max_{0\leq j <p} (g(Y_{j,x}^{a}) - g(Y_{j,x}^{\ba_n})) \bigvee
    \Big(\hV_{n+1}(F^{a}(x)) - V_{n+1}(F^{\ba_n}(x)) \Big) \\
  & \leq & \Big(\max_{0\leq j <p} [g]|Y_{j,x}^{a} - Y_{j,x}^{\ba_n}| \Big) \bigvee
    \Big(\hV_{n+1}(F^{a}(x)) - V_{n+1}(F^{\ba_n}(x)) \Big).
\eeno
We use the decomposition and following bounds
\beno
  &     & \hspace{-2cm} \hV_{n+1}(F^a(x)) - V_{n+1}(F^{\ba_n}(x)) \\
  & =   & \Big(\hV_{n+1}(F^a(x)) - V_{n+1}(F^a(x))\Big)
           + \Big(V_{n+1}(F^a(x))-V_{n+1}(F^{\ba_n}(x))\Big) \\
  &\leq & \Big(\hV_{n+1}(F^a(x)) - V_{n+1}(F^a(x))\Big) + [V_{n+1}]|F^a(x) - F^{\ba_n}(x)|.
\eeno

We deduce from the previous estimates
\beno
  & & \hspace{-2cm} \E\left[\hV_n(X_n) - V_n(X_n)\right] 
  \leq \inf_{a\in \numA} \bigg( 
   [g]\ \E\Big[\max_{0\leq j \leq p} \left|Y_{j,X_n}^{a} - \YcA_{j,X_n}^{\ba_n}\right| \Big]
   \\
  & &  +\ \E\big[(\hV_{n+1} - V_{n+1}){(F^a(X_n))}\big] 
       + [V_{n+1}]\ \E\big[|F^a(X_n) - F^{\ba_n}(X_n)|\big] 
    \bigg)+ \eta_n \\
  & \leq &  \inf_{a\in \numA} \bigg( 
      \E\big[(\hV_{n+1} - V_{n+1}){(F^a(X_n))}\big] 
      + C_F ([g]+[V_{n+1}]) \dt \E \big[|a(X_n)-\ba_n(X_n)|\big] \bigg) + \eta_n
\eeno
where the estimate of Lemma~\ref{lem:estim-a-abar}$(ii)$ has been used for the last inequality.

Then by using the change of variable Lemma~\ref{lem:cv},
we obtain
\beno
    e_n 
     & \leq & \inf_{a\in\numA_n}  
          \bigg( C_{n,\dt} e^{dc[f^a]\dt}\E\big[(\hV_{n+1} - V_{n+1})(X_{n+1})\big] 
          + C_F ([g]+[V_{n+1}]) \dt \E \big[|a(X_n)-\ba_n(X_n)|\big] \bigg) +\eta_n\\
     & \leq & \inf_{a\in\numA_n}  C_n^a e_{n+1} + \big(\eps^{a}_n +\eta_n\big)
\eeno
where $\eps^{a}_n :=C_F  ([g]+ [V_{n+1}])\, \dt \, \E\Big[|a(X_{n}) - \ba_n(X_n)|\Big]$
and $C_n^a := C_{n,\dt} e^{dc[f^a] \dt}$.
By induction, and using the fact that $e_N=0$ because
$\hV_N = V_N$, we obtain (for given coefficients $\eta_k>0$, $k=n,\dots,N-1$):
\beno
  e_n \leq \inf_{(a_n,\dots,a_{N-1})\in\bigotimes_{k=n}^{N-1}\numA_{k}}  \left[
        (\eps^{a_n}_n + \eta_n) 
        + C_n^{a_n} (\eps^{a_{n+1}}_{n+1} + \eta_{n+1}) 
        + \dots 
        + C_n^{a_n}\cdots C_{N-2}^{a_{N-2}}(\eps^{a_{N-1}}_{N-1} + \eta_{N-1})\right].
\eeno
However, we can improve this bound.
For a given $a_n \in \numA_n$, we have a constant $C_n^{a_n}$ which depends of the Lipschitz constant $[a_n]$. 
We can chose a coefficient $\eta_{n+1}=\eta_{n+1}([a_n])>0$ (which may have a dependency over $[a_n]$), 
and proceed in the same way.
By Prop.~\ref{prop:eps-weak-approx-refined}, there exists $\ba_{n+1}\in\cA_{n+1}$, Lipschitz continuous,
such that 
$$
    \big| \E[V_{n+1}(X_{n})]-\E[G^{\ba_{n+1}}\vee V_{n+2}(F^{\ba_{n+1}}(X_{n+1}))] \big| \leq  \eta_{n+1}([a_n]).
$$
Then we obtain the bound
\beno
    e_n 
     & \leq & \inf_{(a_n,a_{n+1}) \in\numA_n \times \numA_{n+1}}
        \bigg(\big(\eps_n^{a_n} +\eta_n\big)  + 
         C_n^{a_n} \big(\eps_{n+1}^{a_{n+1}} +\eta_{n+1}([a_n])\big)
         + C_n^{a_n} C_{n+1}^{a_{n+1}} e_{n+2} 
         \bigg).
\eeno
At the next step, we can chose a coefficient $\eta_{n+2}=\eta_{n+2}([a_n],[a_{n+1}])$, and so on.
By induction, we conclude to the desired bound.
\end{proof}

\newcommand{\cP}{{\mathcal P}}
\newcommand{\card}{{\mathrm Card}}

\newcommand{\VOID}[1]{#1} 
\newcommand{\VOIDBIS}[1]{#1} 

\newcommand{\SLTWOH}{H-scheme\xspace} 
\newcommand{\figw}{coucou1} 
\newcommand{\ldir}{coucou2} 
\newcommand{\rdir}{coucou3} 

\newcommand{\Example}[2]{\subsection{Example \arabic{subsection}: #1 \label{ex:#2}}} 



\section{Numerical results}\label{sec:num}

In the following $d$-dimensional examples (where $d\geq 2$),
two-dimensional "local" and "global" errors are computed in the following way.
Depending on the example, a two-dimensional plane of reference $\cP={\rm Vect}(w_1,w_2)$ is set (passing through the origin),
where $w_1,w_2$ are chosen vectors of $\R^d$, and
a uniform grid mesh $x_k=a_{i} w_1 + b_j w_2$ (for $k=(i,j)$, $|a_i|,|b_j|\leq R_{\max}$
for a given $R_{\max}>0$) is chosen in the plane $\cP$ in order to compute the
exact solution and to compare with the numerical solution. 

Given a threshold $\eta>0$, the errors are computed at the last iteration by
\begin{align*}
	e^{\eta}_{L^1_{\eta}} \coloneqq
	 \frac
         {\sum\limits_{\{x_i \in \mO_\eta\}}  |v(0,x_i)- \hV_0(x_i)|}
         {\sum\limits_{\{x_i \in \mO_\eta\}}  1}
	\quad\text{and}\quad
        e^{\eta}_{L^\infty_{loc}} \coloneqq \max_{\{x_i \in \mO_\eta\}}  |v(0,x_i)- \hV_0(x_i)|
\end{align*}
where $v(0,\cdot)$ is the analytical solution at time 
$t=0$, $\hV_0$ is its approximation by the scheme used, $\mO_\eta\coloneqq\{x \in \mO, |v(0,x)|\leq \eta\}$, and $\Omega$ is the bounded computational domain. Notice that 
$$
  e_{L^1_{\eta}}\simeq \frac{\|v(0,.)\|_{L^1(\mO_\eta)}}{\|1\|_{L^1(\mO_\eta)}},
  \quad 
  e_{L^\infty_{\eta}}\simeq \max_{x\in \mO_\eta} |v(0,x)|.
$$

The global errors, corresponding to the case $\eta=+\infty$,  will be denoted $e_{L^1}$ and $e_{L^\infty}$ 
(i.e.,  $\mO_\eta=\mO$).
Unless otherwise stated, the local errors are computed with $\eta=0.1$, and denoted $e_{L^1_{loc}}$ and $e_{L^\infty_{loc}}$.


We use feedforward neural networks with ReLu activation function on the inner layers. 
Some other activation functions were also tested, including the sigmoid and the tanh functions. 
We found that ReLu was performing better on our cases, and we report only these results.
If not otherwise stated, the output activation function is the identity,
and the Heun scheme is used for $F_\dt$, with $p=5$
substeps (excepted for Example 1 where $p=1$ and $p=5$ are compared).

Implementation of neural networks uses python TensorFlow 2, with Adam optimizer (see \cite{tensorflow2015-whitepaper}),
and the architecture is an Intel Xeon Gold 6140 Processor with 2 CPUs and a total of 36 cores.

\subsection{Example \arabic{subsection} : Rotation with obstacle}
\label{ex:rot} 

\begin{figure}[!hb]
\VOID{
\centering
\renewcommand{\figw}{0.32\linewidth}
\renewcommand{\ldir}{numrot}
\includegraphics[width=\figw]{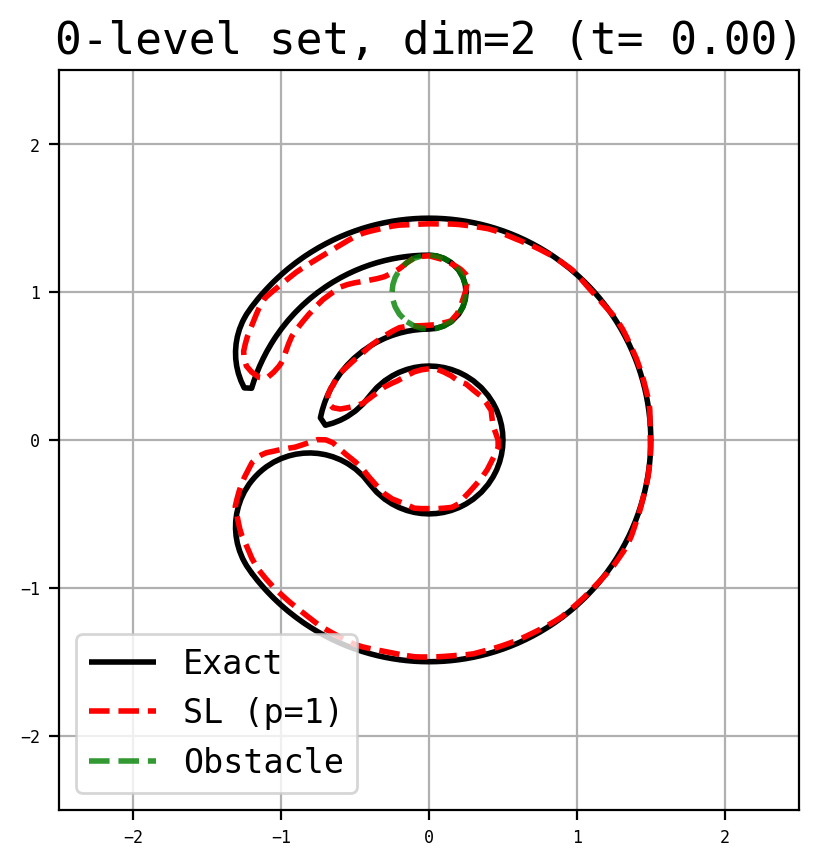}
\includegraphics[width=\figw]{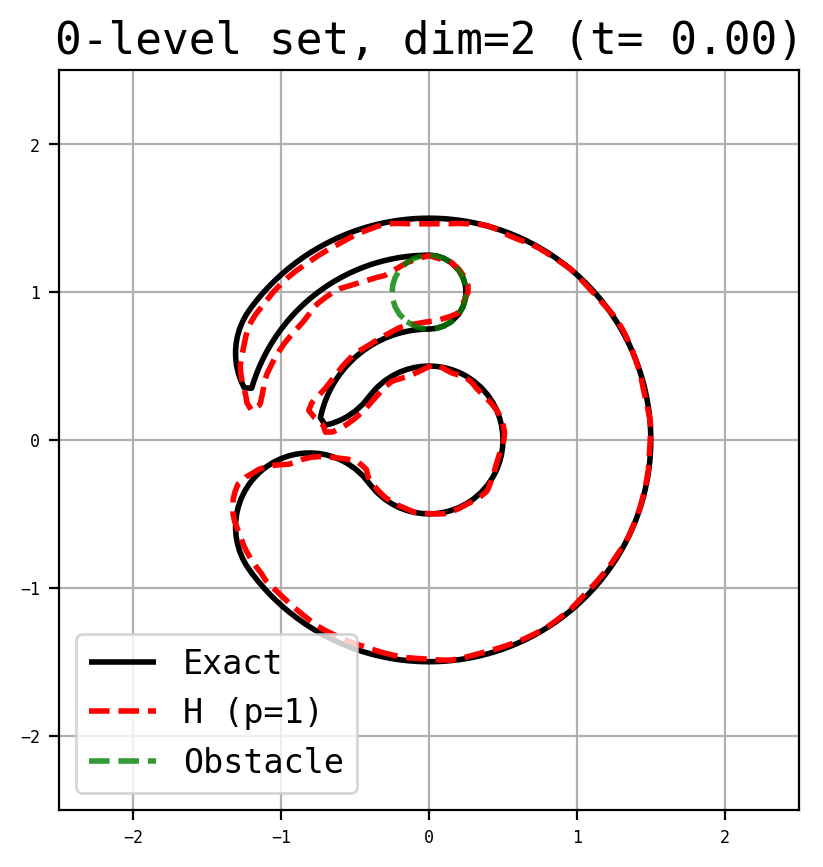}
\includegraphics[width=\figw]{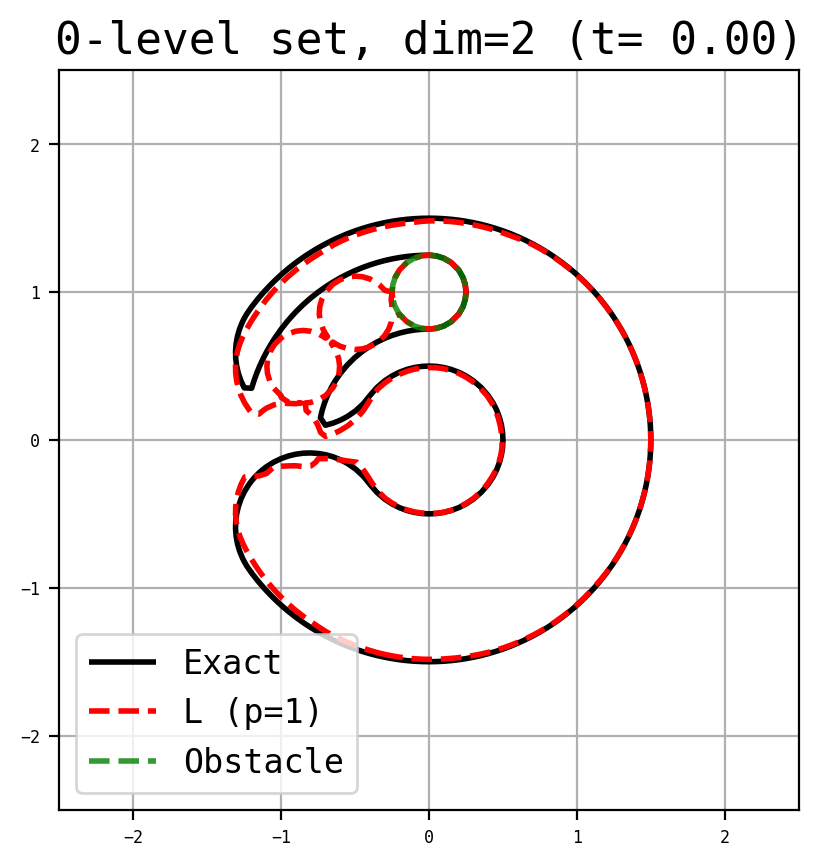}
\\
\includegraphics[width=\figw]{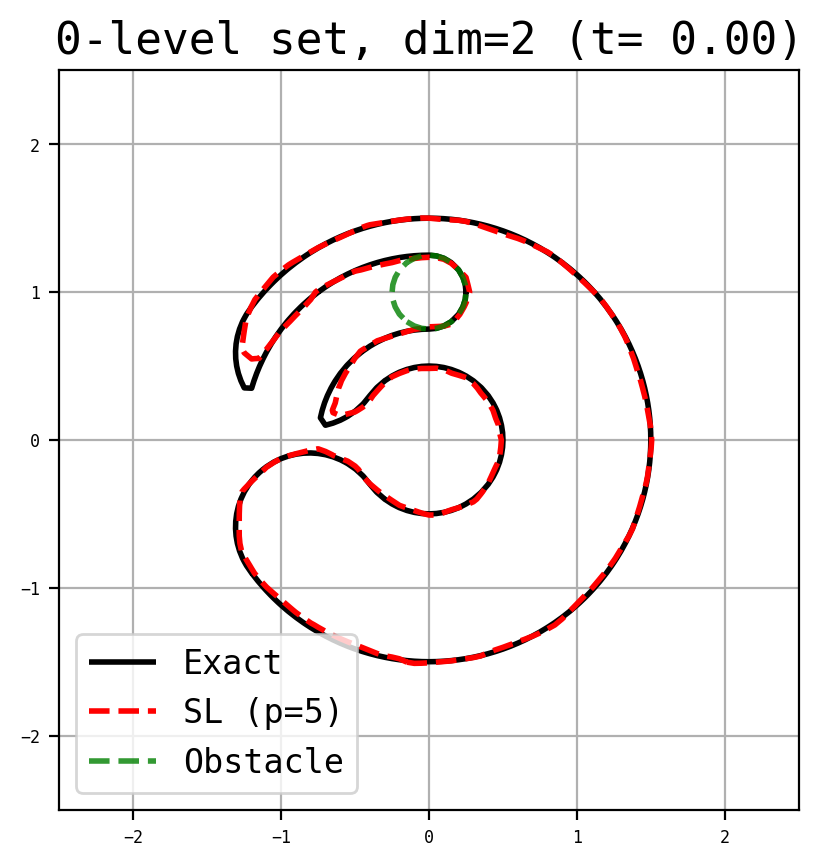}
\includegraphics[width=\figw]{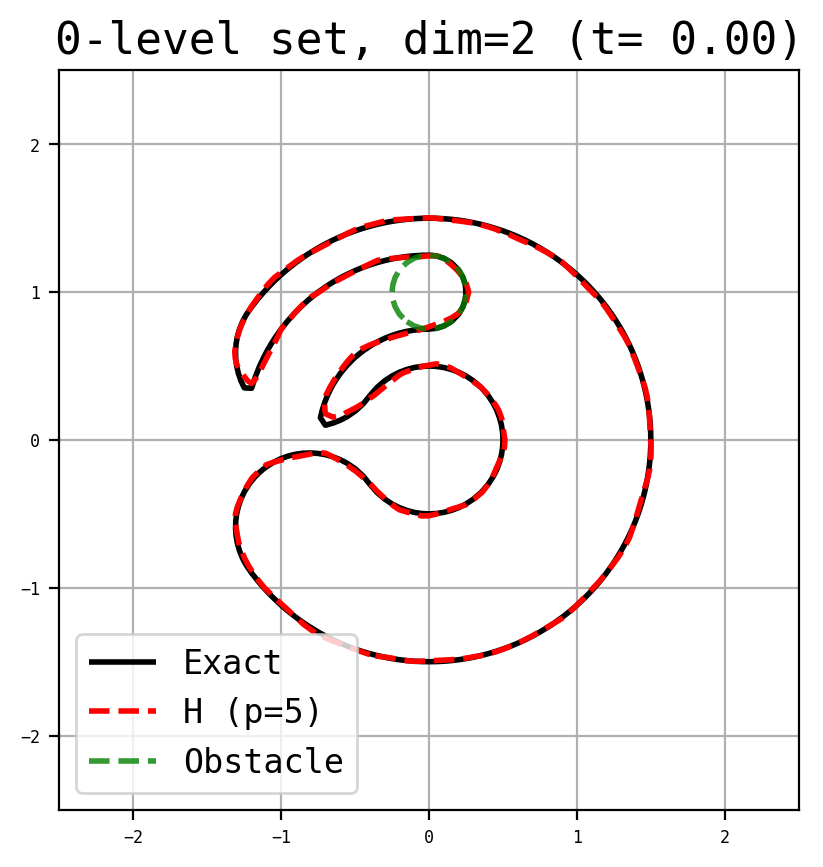}
\includegraphics[width=\figw]{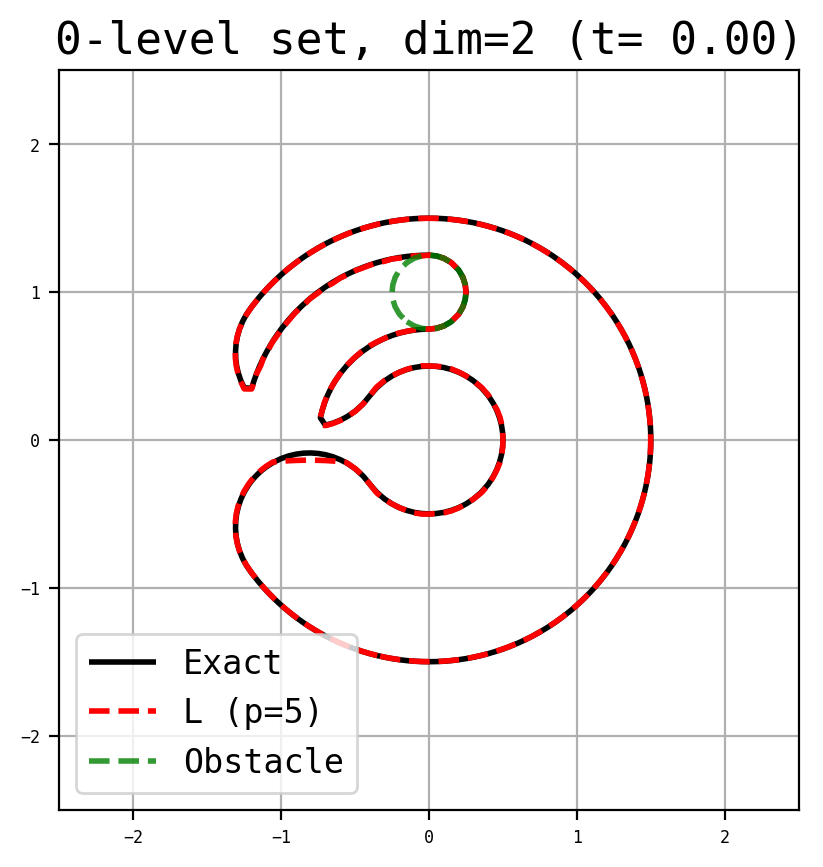}
\\
\includegraphics[width=\figw]{\ldir/SL2sO20/%
fig_NNsurf5_SL2_D2_Nit5_Arelu_Neu40_Lay3_Next10_Order20_Batch1000_MH0_VbistruncLmax_nbOpt1_FeedForward01}
\includegraphics[width=\figw]{\ldir/SL2fO20/%
fig_NNsurf5_SL2_D2_Nit5_Arelu_Neu40_Lay3_Next10_Order20_Batch1000_MH0_VbistruncLmax_nbOpt1_FeedForward01}
\includegraphics[width=\figw]{\ldir/SL3O20/%
fig_NNsurf5_SL3_D2_Nit5_Arelu_Neu40_Lay3_Next10_Order20_Batch1000_MH0_VbistruncLmax_nbOpt1_FeedForward01}
} 
\caption{
 (Example \arabic{subsection})
 The \SLTWO (left), the \SLTWOH (middle) and the \SLTHREE (right) are 
 tested with Euler scheme with $p=1$ (top) and Heun scheme with $p=5$ (middle/bottom).
  The bottom figures corresponds to the surface plots of $z=v(0,x,y)$ (blue), the plot of the  obstacle function (orange), 
  and the 0-level set (red line).
 Networks uses 3 hidden layers, 40 neurons, with $N=5$ time steps. 
\label{fig:exrot}
} 
\end{figure}

This first problem is a two-dimensional example. 
We aim at computing the backward reachable set of a target disk $\mathscr{D}(x_A,r_0)$ before time $T$,
while avoiding the region $\mathscr{D}(x_B,r_1)$ 
with the following parameters
$$
  x_A=(1,0),\quad x_B=(0,1), \quad r_0=0.5, \quad  r_1=0.25, \quad \mbox{and}\ T=0.4
$$
(see Fig.~\ref{fig:exrot}).
The dynamics $f(x,a)$ with controls $a\in[-1,1]$ 
is given by  
$$ \mbox{$f((x_1,x_2),a):= 2\pi a (-x_2,  x_1)$ with $a\in A:=[-1,1]$}
$$
and corresponds to a clockwise to counter-clockwise rotation.
We set 
$$ \varphi(x) :=  \|x-x_A\|_2 - r_0 \quad \mbox{and} \quad 
         g(x) := r_1 - \|x-x_B\|_2.
$$
The value $v(t,x)$ of this problem problem (as defined in \eqref{pb:1}, with $t\in [0,T]$ and $x\in \R^2$)
is also solution of the following HJB equation with an obstacle term
\be \label{eq:exrot-pde}
  & & \min(- v_t  + \max_{a\in[-1,1]} f(x,a)\cdot \nabla_x v, \ v-g(x)) = 0, \quad t\in [0,T], \\
  & & v(T,x)=\max(\varphi(x),g(x)).
\ee
Here, the control networks use the sigmoid output activation function, with value in $[0,1]$, and is converted to $[-1,1]$ by a linear transformation. 


In Fig~\ref{fig:exrot}, we compare the \SLTWO, the \SLTWOH and the \SLTHREE. 
Errors are given in Table~\ref{tab:exrot-errors}. 

We first investigate the influence of the substeps ($p\geq 1$). We choose 
$F_h$ as the Heun scheme, with $N=5$ time steps ($\dt=T/N$), 
and compare the results using $p=1$ or $p=5$ (recall that $p$ is the number of substeps in order to approximate the caracteristic
with a constant control $a$ on a given time interval $[t_k,t_k+\dt]$).

The results, for all schemes, are clearly in favor of using $p=5$ (better characteristic approximation) which 
benefit from the regions of regularity of the control. Hence, for the forthcoming examples, we will always use the Heun scheme with $p=5$.

Notice that for this low-dimensional example ($d=2$), only a small number of stochastic gradient iterations is enough to obtain reasonable results, and
in particular to observe the contribution of $p$. 

We also compare the three schemes for $p=5$, looking at the relative errors. We observe that the L-scheme gives the best results,
the H-scheme gives intermediate results and the SL-scheme is less precise.
 Here we observe that a local $L^1$ relative error less or equal to $10^{-2}$ 
 corresponds to an almost perfect result to the eye.



\begin{table}[H]
	\centering
	\begin{tabular}{|c|ccccc|cc|cc|r|} \hline
	\multirow{2}*{Scheme} & \multicolumn{5}{c|}{Parameters} & \multicolumn{2}{c|}{Global errors} & \multicolumn{2}{c|}{Local errors} & \multirow{2}*{$\begin{matrix}\text{CPU time} \\\text{(s.)}\end{matrix}$} \\ 
	\cline{2-10} & $N$ & lay. & neur. & $M$ & S.G. it. & $L_{\infty}$ & $L_{1}$ rel. & $L_{\infty}$ & $L_{1}$ rel. & \\ 
	\hline \hline 
	SL ($p=1$) &  5 &  3 &  40 & 1000 & 1000 & 2.23e-01 & 5.23e-02 & 1.07e-01 & 3.14e-02 &   133.00 \\ \hline 
	SL ($p=5$) &  5 &  3 &  40 & 1000 & 1000 & 1.22e-01 & 1.72e-02 & 1.02e-01 & 1.19e-02 &   182.94 \\ \hline \hline
	 H ($p=1$) &  5 &  3 &  40 & 1000 & 1000 & 2.12e-01 & 5.39e-02 & 1.13e-01 & 2.58e-02 &   180.18 \\ \hline 
	H  ($p=5$) &  5 &  3 &  40 & 1000 & 1000 & 1.20e-01 & 7.96e-03 & 7.83e-02 & 7.93e-03 &   285.84 \\ \hline \hline
	 L ($p=1$) &  5 &  3 &  40 & 1000 & 1000 & 5.99e-01 & 4.74e-02 & 5.01e-01 & 2.55e-02 &    54.42 \\ \hline 
	L  ($p=5$) &  5 &  3 &  40 & 1000 & 1000 & 2.10e-01 & 3.22e-03 & 2.00e-01 & 4.08e-03 &   106.46 \\ \hline 
	\end{tabular}
	\caption{(Example \arabic{subsection}) Comparison of schemes
	\label{tab:exrot-errors}
	}
\end{table}

\begin{figure}[!hb]
\VOID{
\centering
\renewcommand{\figw}{0.32\linewidth}
\renewcommand{\ldir}{numrot/DGM}
\includegraphics[width=\figw]{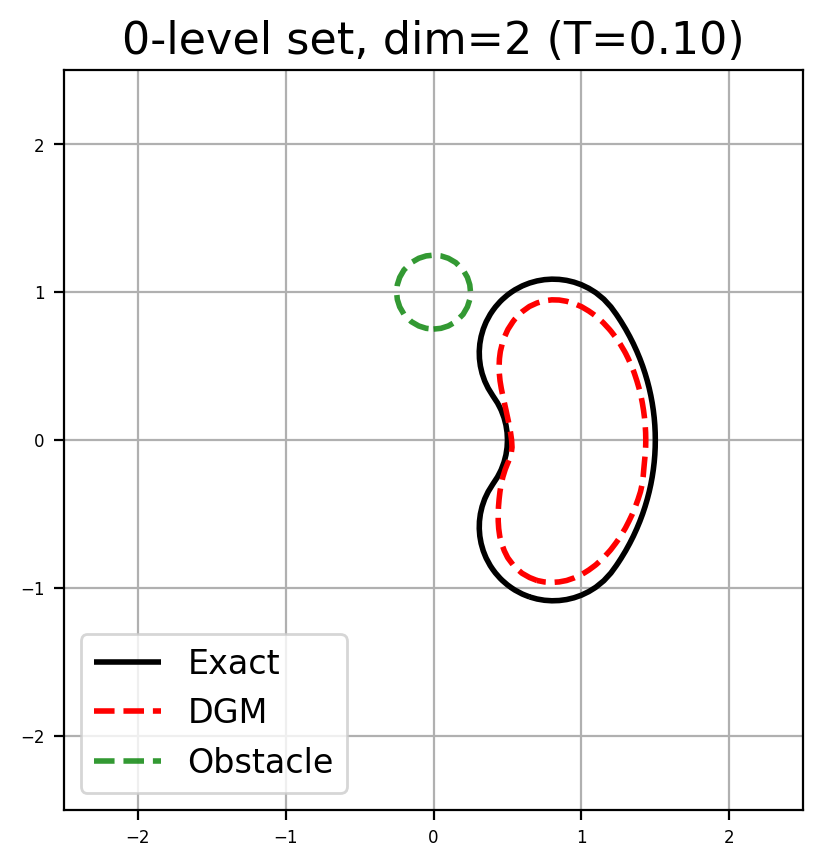}
\includegraphics[width=\figw]{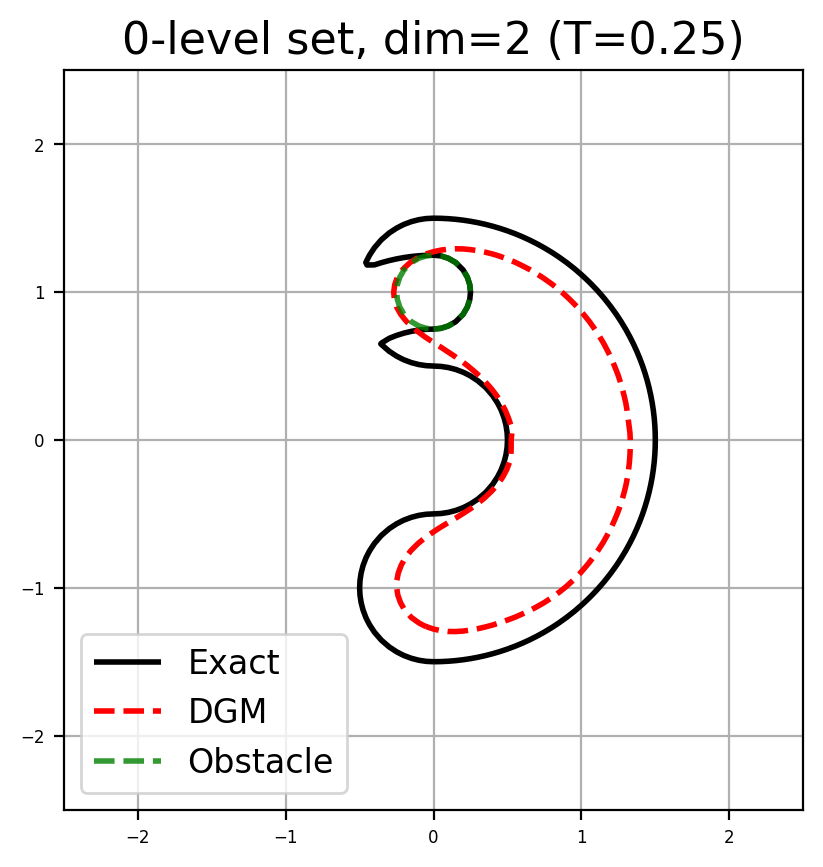}
\includegraphics[width=\figw]{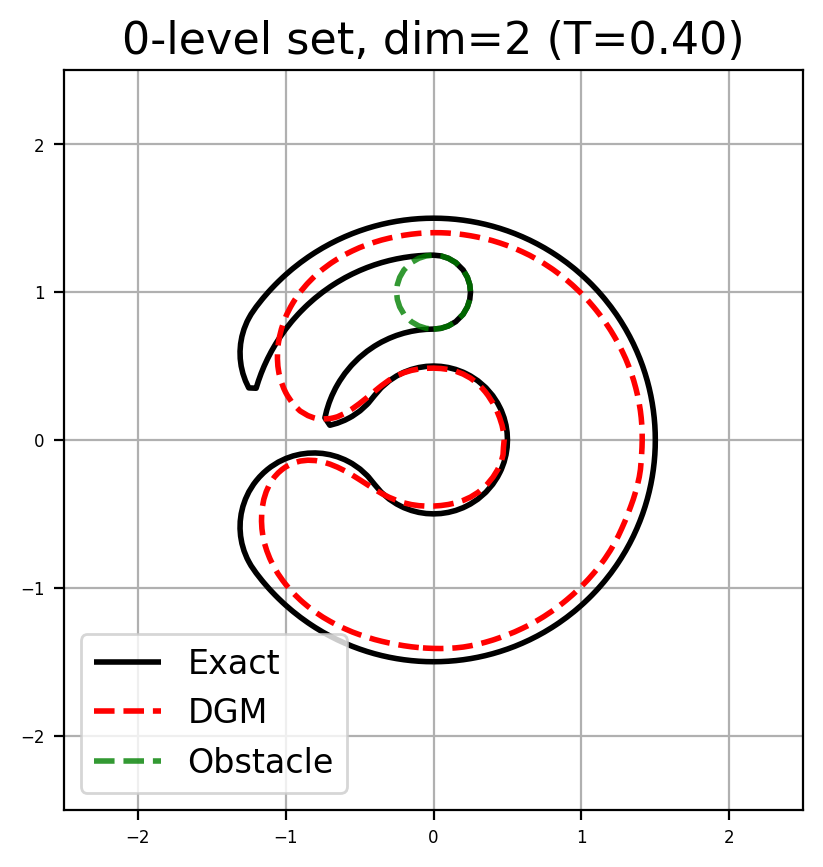}
} 
\caption{
 (Example \arabic{subsection})
 DGM direct method in dimension $d=2$ (the computation is done in dimension $d+1$ to include time). 
 Results are given at $t=0$.
 Final time is set to $T=0.1$ (left), $T=0.25$ (middle) and $T=0.4$ (right).
\label{fig:exrot-DGM}
} 
\end{figure}

Finally, on this example, we have also tested a direct method (the DGM approach of \cite{sirignano2018dgm}),
where a global space-time DNN is used in order to approximate
the value $(t,x)\conv v(t,x)$ solution of the PDE~\eqref{eq:exrot-pde}. 
However, in our experiments, we found that the DNN in general
fails to see the obstacle part of the solution. 
A typical illustration is given in Figure~\ref{fig:exrot-DGM}, where 3 simulations with increasing 
final time $T$ are presented.
We considered neural networks with $tanh$ activation function, both in the inner and output layers. 
In the presented results, the network uses 3 inner layers of 40 neurons. At each iteration of the minimization, 
the stochastic gradient draws 10,000 points in the space-time domain and 1000 points on the border $t=T$
{($100,000$ iterations of stochastic gradient used).}

\subsection{Example \arabic{subsection} : eikonal equation}
\label{ex:eik} 

Next we consider a $d$-dimensional problem, with no obstacle term, for various dimensions $d=6,7,8$
(the next examples will consider obstacles).

More precisely the dynamics is $f(x,a):= a$ with $a\in A\coloneqq\overline{\mathscr{B}}(0,1)$, the closed unit ball of $\R^d$ (for the Euclidean norm).
The function $\varphi$ is 
  $$ \varphi(x):=\min\bigg( \|x-x_A\|_2 - r_0, \  \|x-x_B\|_2 - r_0\bigg)  $$
  with $x_A=(1,0,\dots,0)$ and  $x_B=(-1,0,\dots,0)$,
and parameters $T=1.0$ and $r_0=0.5$.
Hence the value is defined as the solution of \eqref{pb:1} (with $g:=-\infty$). 
The analytical solution is known and given by $v(t,x)=\min\big( (\|x-x_A\|_2-(T-t))_+ - r_0, \  (\|x-x_B\|_2-(T-t))_+ - r_0\big)$.

The corresponding HJB equation (for $x\in \R^d$), using 
$\max_{a\in A} f(x,a)\cdot \nabla_x v = \|\nabla_x v\|$, is the following eikonal equation
\be \label{eq:exeik-pde}
  & & - v_t  + \|\nabla_x v\| = 0, \quad t\in [0,T] \\
  & & v(T,x)=\varphi(x).
\ee

\begin{figure}[H]
\VOID{
\centering
\renewcommand{\figw}{0.32\linewidth}
\renewcommand{\rdir}{numeik/}
\renewcommand{\ldir}{2022-06-24-07h17-save-eik-SL3-D08-N04-NEXT4000-4000-o01-B1000}
\includegraphics[width=\figw]{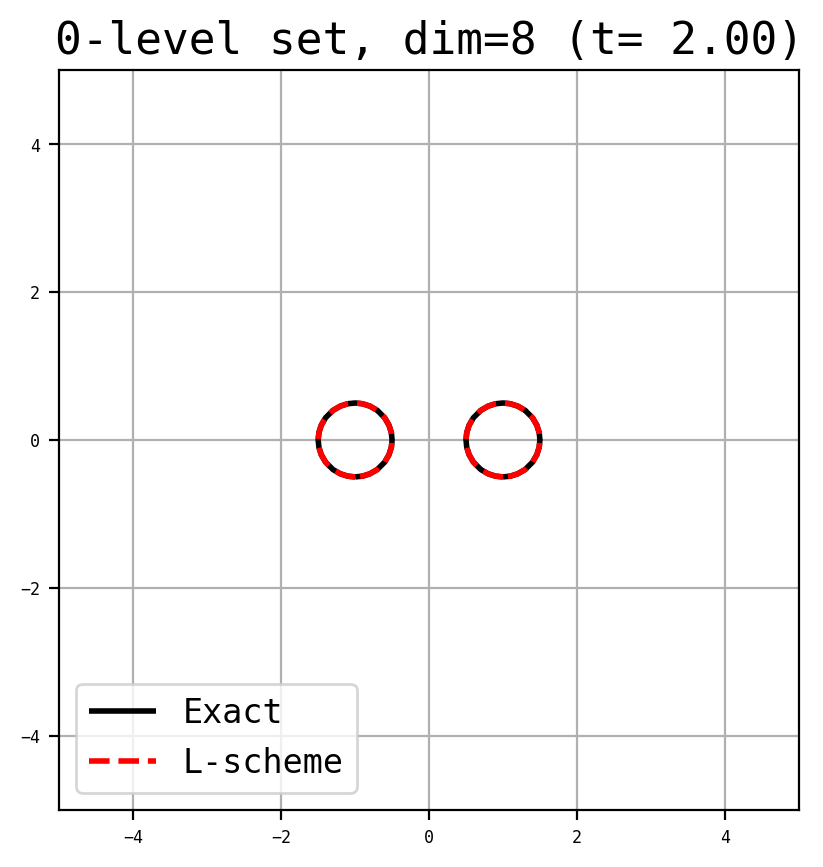}
\includegraphics[width=\figw]{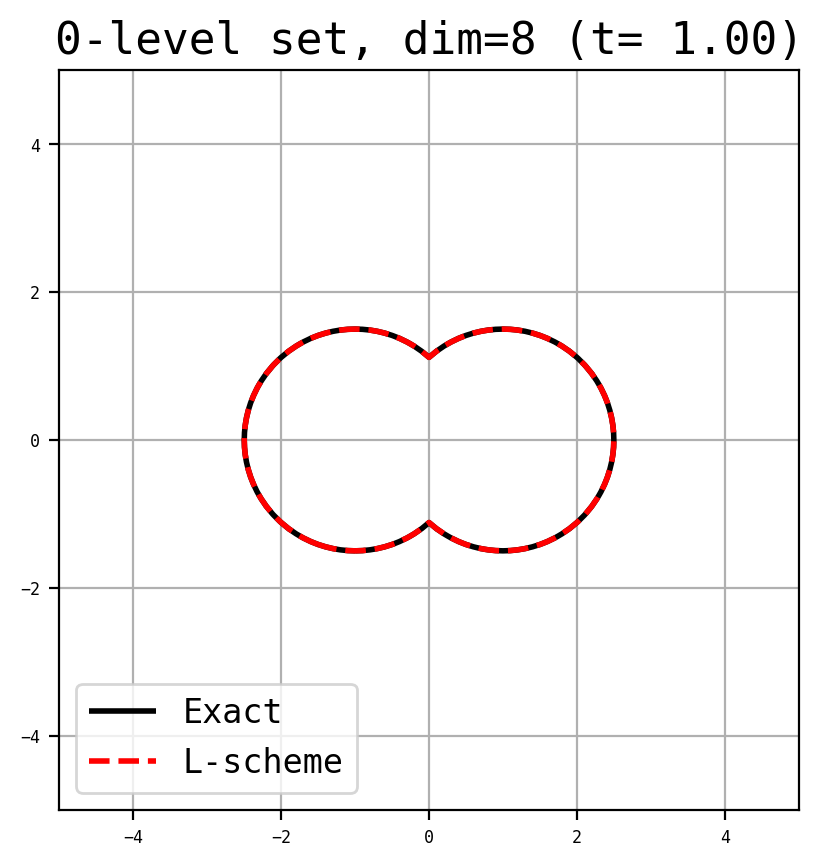}
\includegraphics[width=\figw]{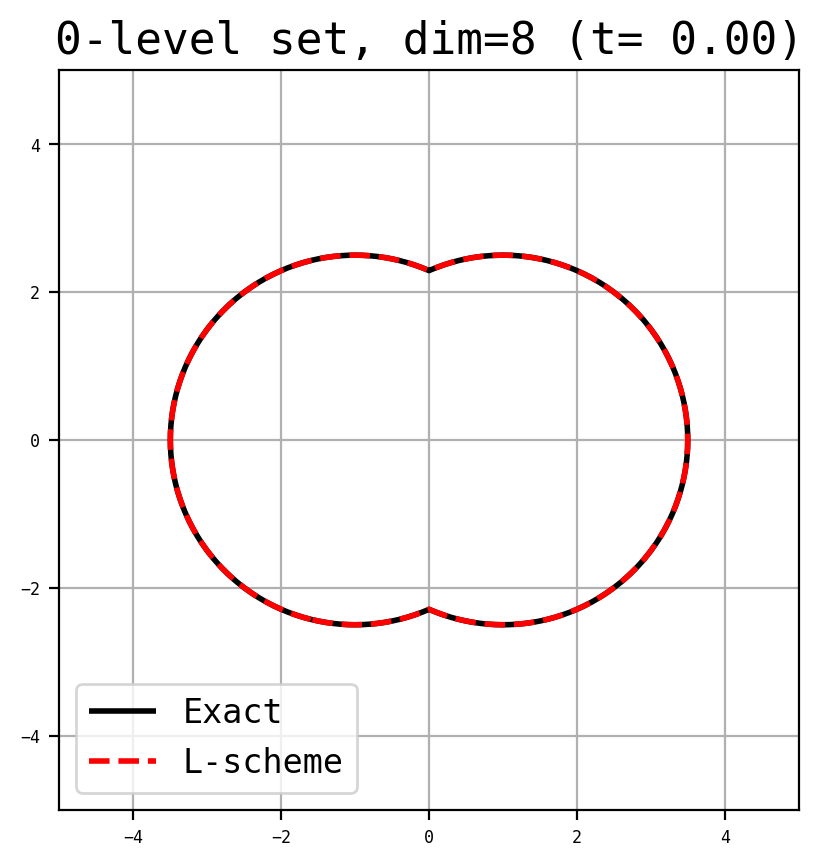}
\\
} 
\caption{
  (Example \arabic{subsection}) 
  Eikonal equation, \SLTHREE,
  dimension $d=8$,
  at time $t=T=2.0$ (left, terminal condition), $t=1.0$ (center), $t=0.0$ (right).
  Networks of 3 hidden layers and 40 neurons; $N=4$ time steps.
  \label{fig:exeik}
}
\end{figure}

Here, we choose the control networks to take their values in $\mathbb{R}^d$. The results are then converted from $\mathbb{R}^d$ to the unit ball 
$\overline{\mathscr{B}}(0,1)$  of $\R^d$
by using the map $p \mapsto \frac{p}{\max(1,\|p\|)}$. 
(Numerical tests showed that the choice of the map may affect the results, and the results may deteriorate in particular when 
using an anisotropic map.)
Errors are given in Table~\ref{tab:exeik-errors} for dimensions $d=6,7,8$,
and some illustrations are given in Fig.~\ref{fig:exeik} for dimension $d=8$
(results for $d\in\left\{6,7\right\}$ are indistinguishable to the eye from the case $d=8$, and they are not included).
Errors and figures are computed in the plane $\cP$ generated by the first two vectors $e_1,e_2$ of the canonical basis of $\R^d$.

In particular we observe that the $L$-scheme performs well (numerical and exact $0$-level sets are indistinguisable to the eye),
as long as a sufficient number of SG iterations is used, and that the control map from $\R^d$ to $\overline{\mathscr{B}}(0,1)$ is well chosen.

\begin{table}[H]
	\centering
	\begin{tabular}{|c|cccccc|cc|cc|c|} \hline
	\multirow{2}*{Scheme} & \multicolumn{6}{c|}{Parameters} & \multicolumn{2}{c|}{Global errors} & \multicolumn{2}{c|}{Local errors} & \multirow{2}*{$\begin{matrix}\text{CPU} \\ \text{time} \end{matrix}$} \\ 
	\cline{2-11} & $d$ & $N$ & layers & neurons & $M$ & S.G it. & $L_{\infty}$ & $L_{1}$ rel. & $L_{\infty}$ & $L_{1}$ rel. & \\ 
	\hline \hline 
	         L & 6 &  4 & 3 &  40 &  1000 & 100000 & 2.16e-02 & 1.96e-03 & 4.06e-04 & 1.58e-04 &  1h26 \\ \hline 
	         L & 7 &  4 & 3 &  40 &  1000 & 200000 & 5.00e-02 & 3.41e-03 & 1.51e-02 & 1.26e-04 &  3h55 \\ \hline 
	         L & 8 &  4 & 3 &  40 &  1000 & 400000 & 1.99e-01 & 1.81e-02 & 4.39e-04 & 2.19e-04 & 10h31 \\ \hline 
	\end{tabular}
	\caption{(Example \arabic{subsection}) $L$-scheme, dimensions $d=6,7,8$
	\label{tab:exeik-errors}
	}
\end{table}


\Example{$d$-dimensional advection with obstacle}{adv}

We now consider an elementary $d$-dimensional advection problem with an obstacle term,
and compare the \SLTWO, the \SLTWOH and the \SLTHREE.
The problem is to reach the target $\{\varphi(x)\leq 0\}$,
while avoiding an obstacle $\left\{g(x) \leq 0\right\}$, 
with linear dynamics $f(x,a) \coloneqq - a e$ where $e \in \R^d$ and the control  $a$ lies in $A\coloneqq[0,1]$.
The corresponding HJB equation is
\be \label{eq:exadv-pde}
  & & \min\bigg(- v_t  + \max_{a\in [0,1]} ae \cdot \nabla_x v, \ v- g(x)\bigg) = 0, \quad t\in [0,T] \\
  & & v(T,x)=\max(\varphi(x),g(x)).
\ee
Equivalently, $\max_{a \in [0,1]} (a e\cdot \nabla v) = \max(0, e \cdot \nabla v)$.
The reachable set at time $t$ is given by $\left\{v(t,\cdot) \leq 0\right\}$ (corresponding to the set of points that can reach the target before time $t$).
The target function $\varphi$ and the obstacle function $g$ are defined by
$$
   \varphi(x) \coloneqq  \|x-A_0\|_2 - r_0 \quad \text{and} \quad g (x) \coloneqq r_1 - \|x-A_1\|_2
$$
so that $\{\varphi(x)\leq 0\} = \overline{\mathscr{B}}(A_0,r_0)$, and $\{g(x)\geq 0\} = \overline{\mathscr{B}}(A_1,r_1)$.
The following parameters are considered:
$$
   e=(1,1,...,1)/\sqrt{d},\quad
   A_0=-(1,1,...,1)/\sqrt{d},\quad
   A_1=(0,0,\dots,0),\quad
   r_0 = 0.5,\quad
   r_1 = 0.25.
$$
Here, the exact solution can be computed as $v(t,x) = \varphi(p(x,t)) \vee g(q(x))$, where 
\begin{align*}
	p(x,t) \coloneqq x - \max\left(0,\min\left(\left<x-A_0, e\right>,T-t\right)\right) e\quad
	\text{and}\quad
	q(x) = x - \max\left(\left<x - A_1, e\right>,0\right) e.
\end{align*}

\begin{figure}[H]
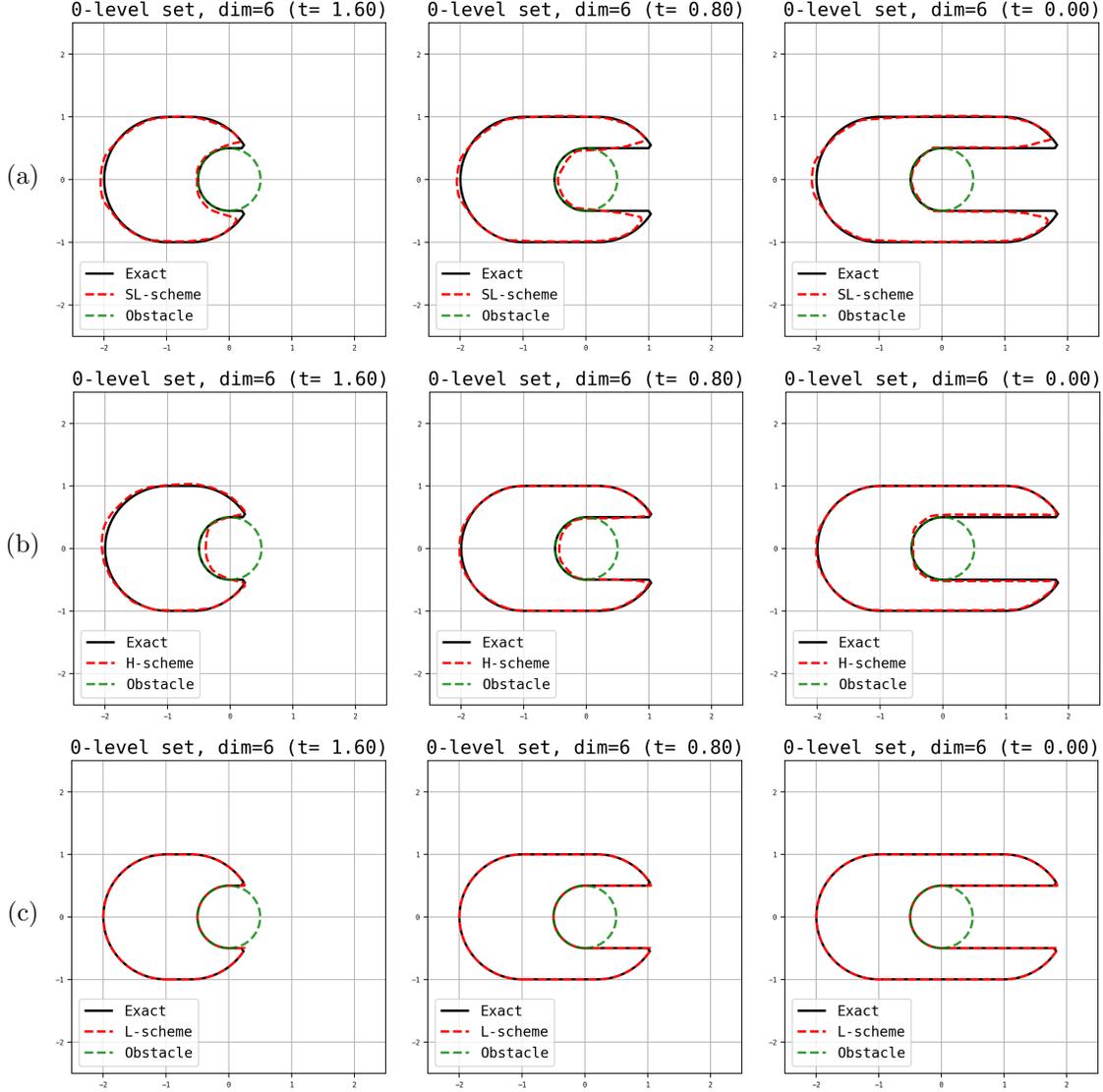

\VOID{
  \centering
  \renewcommand{\figw}{0.30\linewidth}
\renewcommand{\rdir}{numadv/}
  \renewcommand{\ldir}{\VOIDBIS{SL2SEMI_D6}}
  (a)
  \includegraphics[height=\figw,align=c]{\rdir\ldir/show_fig_zerolevel1_SL2_D6_Nit5_Arelu_Neu40_Lay3_Next2000_Order20_Batch2000_MH0_VbistruncLmax_nbOpt1_FeedForward01}
  \includegraphics[height=\figw,align=c]{\rdir\ldir/show_fig_zerolevel3_SL2_D6_Nit5_Arelu_Neu40_Lay3_Next2000_Order20_Batch2000_MH0_VbistruncLmax_nbOpt1_FeedForward01}
  \includegraphics[height=\figw,align=c]{\rdir\ldir/show_fig_zerolevel5_SL2_D6_Nit5_Arelu_Neu40_Lay3_Next2000_Order20_Batch2000_MH0_VbistruncLmax_nbOpt1_FeedForward01}
  \\
  \renewcommand{\ldir}{\VOIDBIS{SL2FULL_D6}}
  (b)
  \includegraphics[height=\figw,align=c]{\rdir\ldir/show_fig_zerolevel1_SL2_D6_Nit5_Arelu_Neu40_Lay3_Next2000_Order20_Batch2000_MH0_VbistruncLmax_nbOpt1_FeedForward01}
  \includegraphics[height=\figw,align=c]{\rdir\ldir/show_fig_zerolevel3_SL2_D6_Nit5_Arelu_Neu40_Lay3_Next2000_Order20_Batch2000_MH0_VbistruncLmax_nbOpt1_FeedForward01}
  \includegraphics[height=\figw,align=c]{\rdir\ldir/show_fig_zerolevel5_SL2_D6_Nit5_Arelu_Neu40_Lay3_Next2000_Order20_Batch2000_MH0_VbistruncLmax_nbOpt1_FeedForward01}
  \\
  \renewcommand{\ldir}{\VOIDBIS{SL3_D6}}
  (c)
  \includegraphics[height=\figw,align=c]{\rdir\ldir/show_fig_zerolevel1_SL3_D6_Nit5_Arelu_Neu40_Lay3_Next2000_Order20_Batch2000_MH0_VbistruncLmax_nbOpt1_FeedForward01}
  \includegraphics[height=\figw,align=c]{\rdir\ldir/show_fig_zerolevel3_SL3_D6_Nit5_Arelu_Neu40_Lay3_Next2000_Order20_Batch2000_MH0_VbistruncLmax_nbOpt1_FeedForward01}
  \includegraphics[height=\figw,align=c]{\rdir\ldir/show_fig_zerolevel5_SL3_D6_Nit5_Arelu_Neu40_Lay3_Next2000_Order20_Batch2000_MH0_VbistruncLmax_nbOpt1_FeedForward01}
  \\
}       
  \caption{
    (Example \refexadv) 
    Results obtained with
    \SLTWOH (row (a)), \SLTWO (row (b)), and \SLTHREE (row (c)) respectively. 
    Dimension $d=6$, $N=5$ time steps, 
    neural networks of 3 layers and 40 neurons.
  }
  \label{fig:exadv}
\end{figure}

For the control networks we use the sigmoid as the output activation function (output in $[0,1]$).
For the figure and error computations, we have chosen a grid in the 2-dimensional plane $\cP=Vect(u,v)$ where
$$
  u=e \equiv (1,1,\dots,1)/\sqrt{d}, \qquad  v=(1,-1,0,\dots,0)/\sqrt{2}.
$$
(Notice that for such parameters the exact 0-level set is the same independently of the dimension $d$).
In order to perform the SG iterations, the size of the random batch points is set to $M=2000$ (as well as for the value approximation by neural networks,
step $(ii)$ of \SLTWO).
Results are given in Table~\ref{tab:exadv-errors} and in Figure~\ref{fig:exadv}, for dimension $d=6$
(the difference between the schemes is more clear when the dimension is not too small).


The CPU time reflects the computational cost of the projection of the value function that is present in the \SLTWO and the \SLTWOH.
Both the \SLTWO and the \SLTWOH need to optimize two networks per time step (one for the control and one for the value),
whereas the \SLTHREE needs only one (for the control). Additionnally, the \SLTWOH  computes
the whole characteristics, leading to a higher CPU time than the \SLTWO. 
(However, if the number of time steps $N$ grows, the \SLTHREE may become 
more expensive than the \SLTWO.)

Looking in particular at the figures in Figure~\ref{fig:exadv}, 
this example shows some kind of numerical diffusion that we may encounter with the~\SLTWO (and with the~\SLTWOH, to a lesser extent).

In this example, we have also numerically observed that an increasing number of stochastic gradient iterations were 
needed as the dimension increases. 



\begin{table}[H]
	\centering
	\begin{tabular}{|c|cccccc|cc|cc|c|} \hline
	\multirow{2}*{Scheme} & \multicolumn{6}{c|}{Parameters} & \multicolumn{2}{c|}{Global errors} & \multicolumn{2}{c|}{Local errors} & \multirow{2}*{$\begin{matrix}\text{CPU} \\ \text{time} \end{matrix}$} \\ 
	\cline{2-11} & $d$ & $N$ & lay. & neur. & $M$ & S.G. it. & $L_{\infty}$ & $L_{1}$ rel. & $L_{\infty}$ & $L_{1}$ rel. & \\ 
	\hline \hline 
	        SL & 6 &  5 &  3 &  40 & 2000 & 200000 & 9.08e-02 & 1.84e-02 & 5.77e-02 & 1.29e-02 & 6h41 \\ \hline 
	         H & 6 &  5 &  3 &  40 & 2000 & 200000 & 9.47e-02 & 1.50e-02 & 6.42e-02 & 1.08e-02 & 8h31 \\ \hline 
	         L & 6 &  5 &  3 &  40 & 2000 & 200000 & 2.14e-03 & 9.58e-05 & 1.79e-03 & 9.54e-05 & 4h59 \\ \hline 
	\end{tabular}
	\caption{(Example \arabic{subsection}) Advection with obstacle, comparison of schemes
	\label{tab:exadv-errors}
	}
\end{table}

\subsection{Example \refeiknc : eikonal advection equation with obstacle, large drift} \label{sec:eiknc}

We consider now a mixed $d$-dimensional eikonal/advection equation with an obstacle term:
\be \label{eq:exeiknc-pde}
  & & \min\bigg(- v_t  + \max_{a\in A} f(x,a)\cdot \nabla_x v, \ v- g(x)\bigg) = 0, \quad t\in [0,T] \\
  & & v(T,x)=\max(\varphi(x),g(x)).
\ee
with $f(x,a) = b e_1 + c a$, where $e_1=(1,0,...,0)^t\in \R^d$, the  control $a$ belongs to $A \coloneqq \mathbb{S}^{d-1}$ the unit ball of $\mathbb{R}^d$,
$b\in \mathbb{R}$ is a coefficient corresponding to the "drift",  
and $c\geq 0$ is a speed coefficient for the eikonal part of the equation.
Equivalently, 
    $$\max_{a \in A} (f(x,a)\cdot \nabla v) = 
       b \cdot \frac{\partial v}{\partial x_1} + c \|\nabla_x v\|.$$

The obstacle term and terminal condition are defined as 
\begin{align*}
        g (x) := \min\left(g_{\max} - c_{e} \left|x_1 - g_c\right|, c_x \left|x_\perp\right| + g_{\min}\right),
        \quad 
        \varphi(x) := \|x\| + \alpha_{\min},
\end{align*}
where $c_e$, $c_x$, $g_c$, $g_{\max}$ $g_{\min}$ and $\alpha_{\min}$ are coefficients,
and, for a given $x=(x_1,\dots,x_d)^t\in \R^d$,  $x_\perp:=(0,x_2,\dots,x_d)^t$ 
(the orthogonal projection of $x$ on vect$(e_2,\dots,e_d)$).
Note that this obstacle term correspond to a wall obstacle with a tube opening centered around the $e_1$ axis 
(see for instance the green dotted line in Fig.~\ref{fig:exeiknc}).

The exact solution can be computed. Details are given in Appendix~\ref{app:A}.
In this example, more precisely, the following parameters are considered
$$ 
   g_{\max} = 2, \quad g_{\min} = -2, \quad c_e = 1, \quad c_x = 1.5, \quad g_c = 4, \quad b = 1, \quad c = 0.5,
   \quad \alpha_{\min}=-1.
$$
Here in particular $|b|>c$: the drift is dominant, which corresponds to a non-controllable situation.

\begin{figure}[H]
\VOID{
	\centering
	\renewcommand{\figw}{.32\linewidth}
	\renewcommand{\rdir}{numeiknc/}
	\renewcommand{\ldir}{EIKO_COMP_D4}
	\includegraphics[width=\figw]{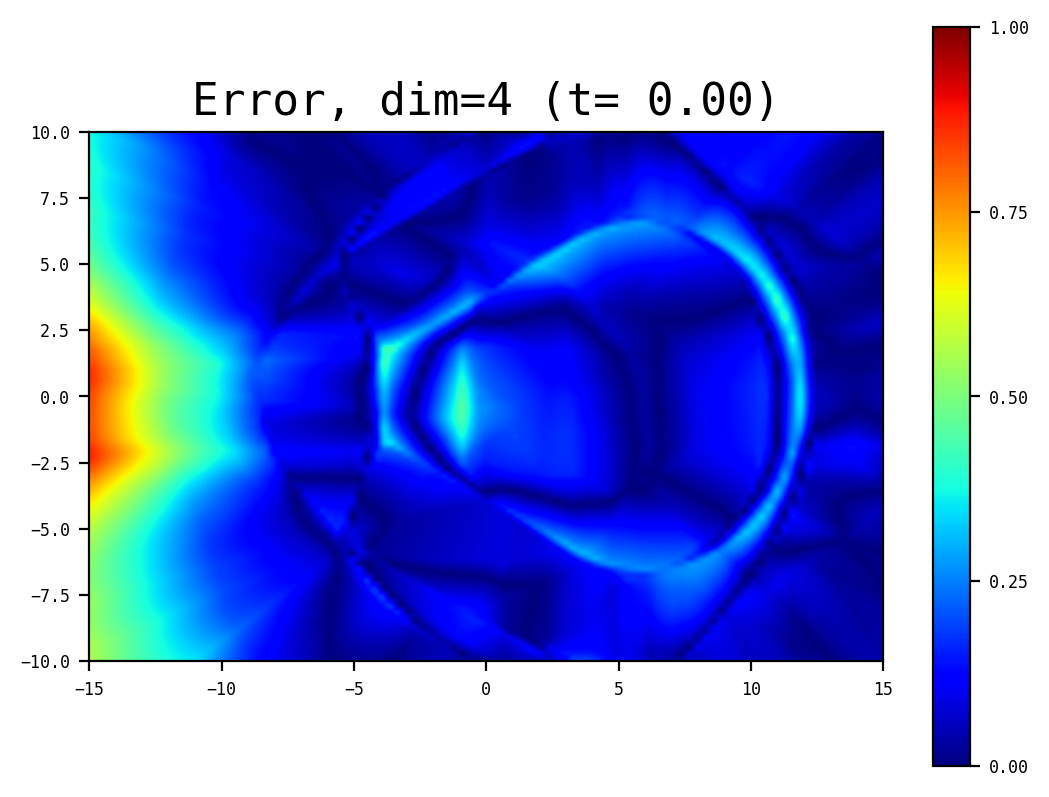}%
	\includegraphics[width=\figw]{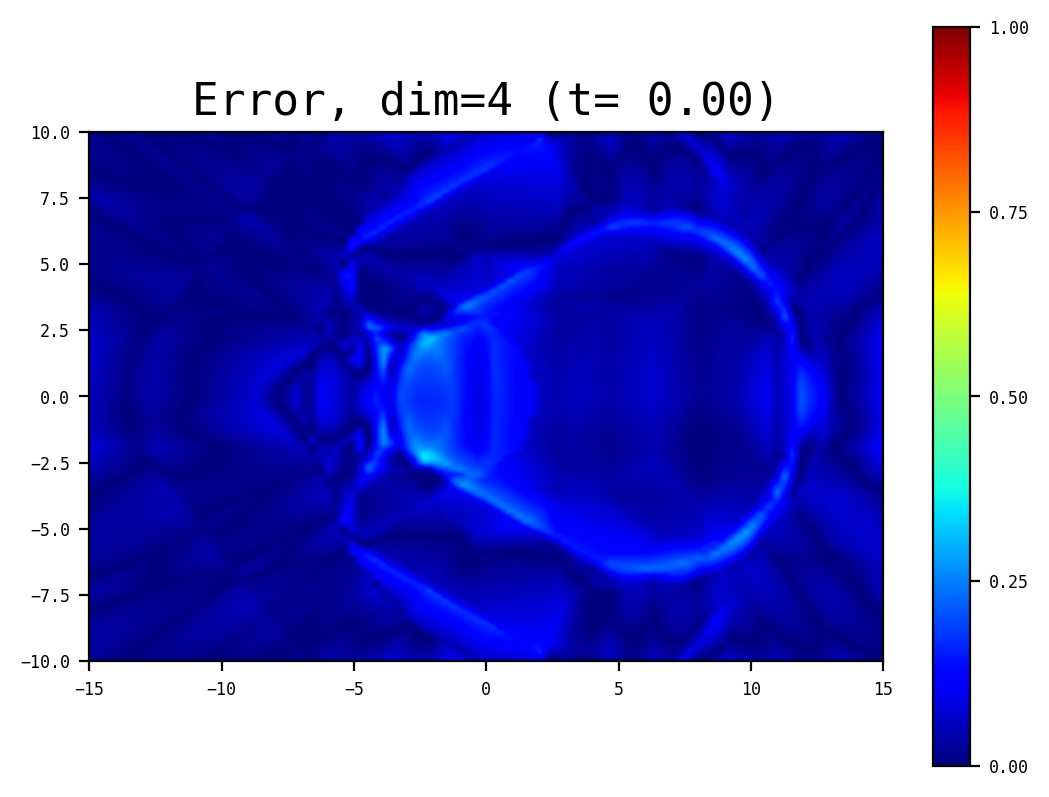}%
	\includegraphics[width=\figw]{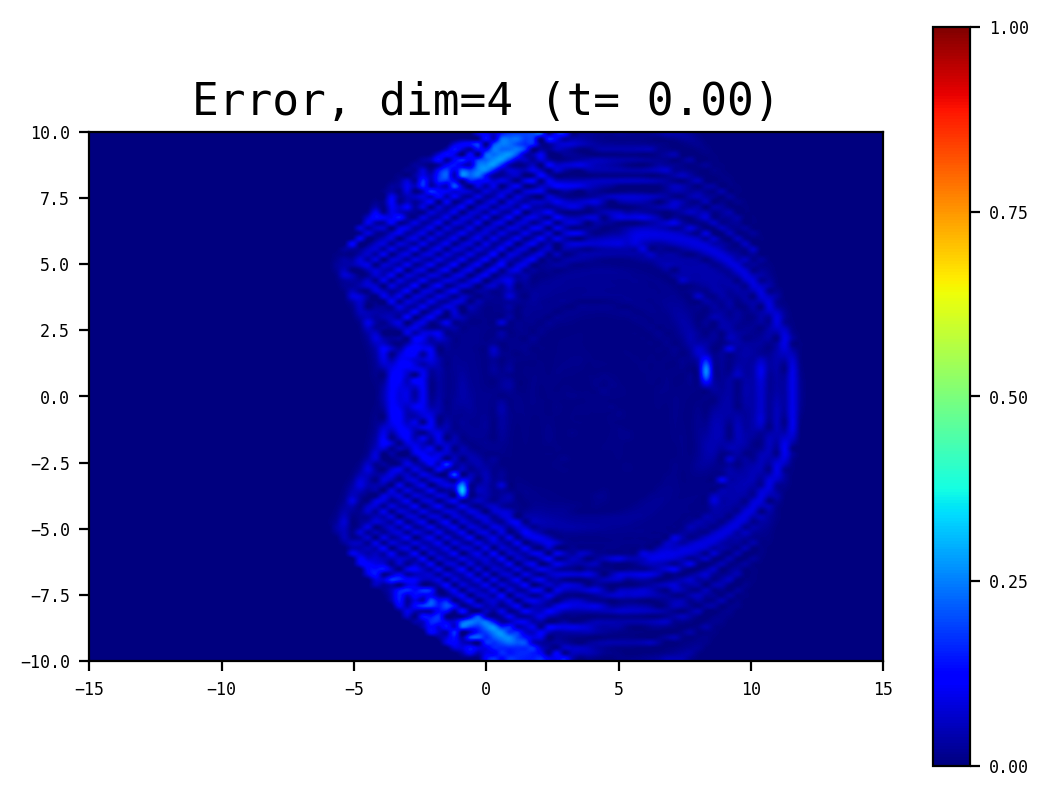}\\
	\renewcommand{\ldir}{EIKO_COMP_D4}
	\includegraphics[width=\figw]{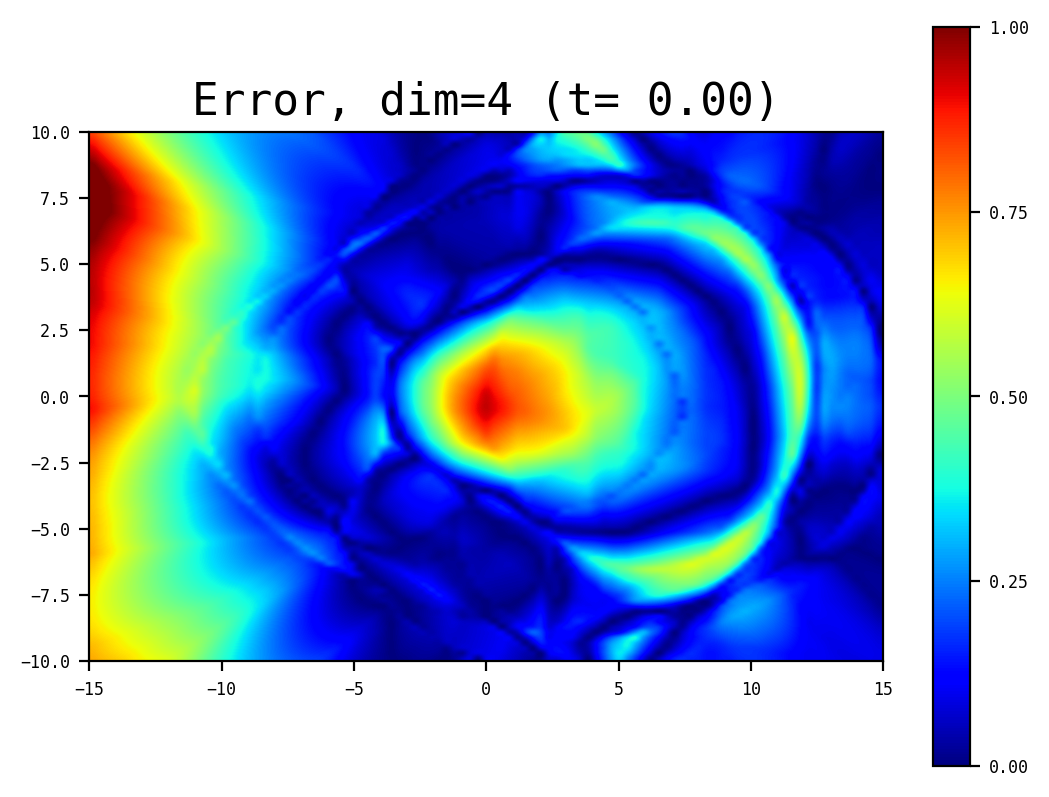}%
	\includegraphics[width=\figw]{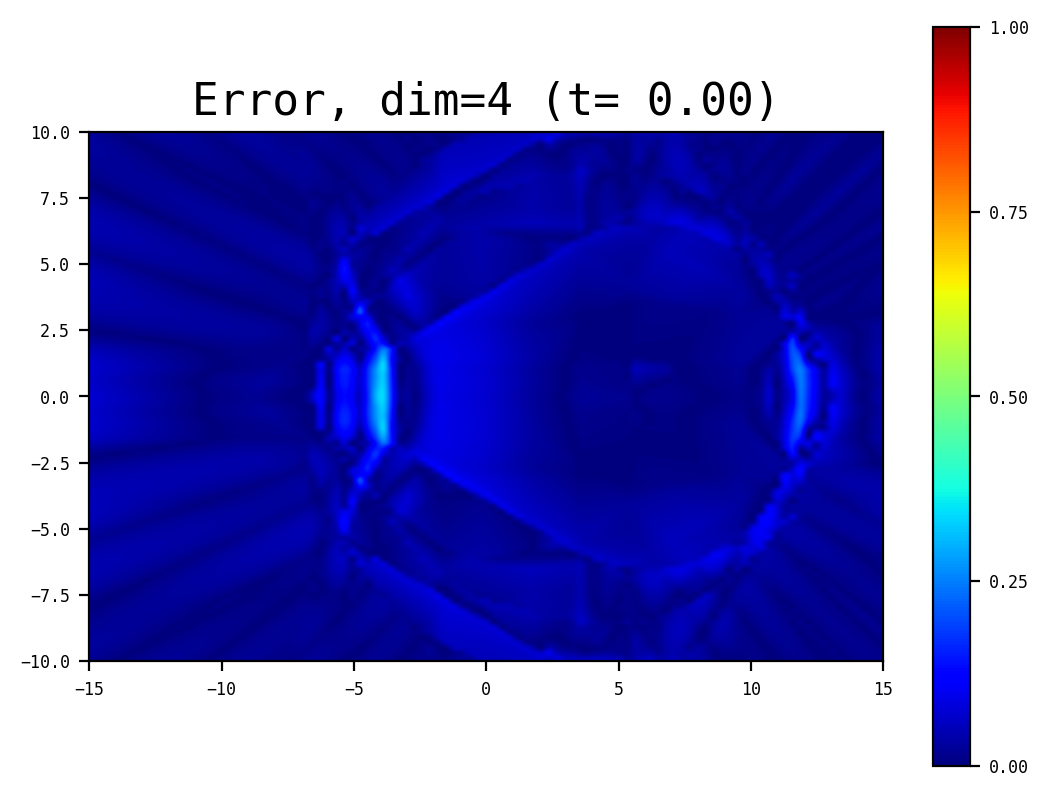}%
	\includegraphics[width=\figw]{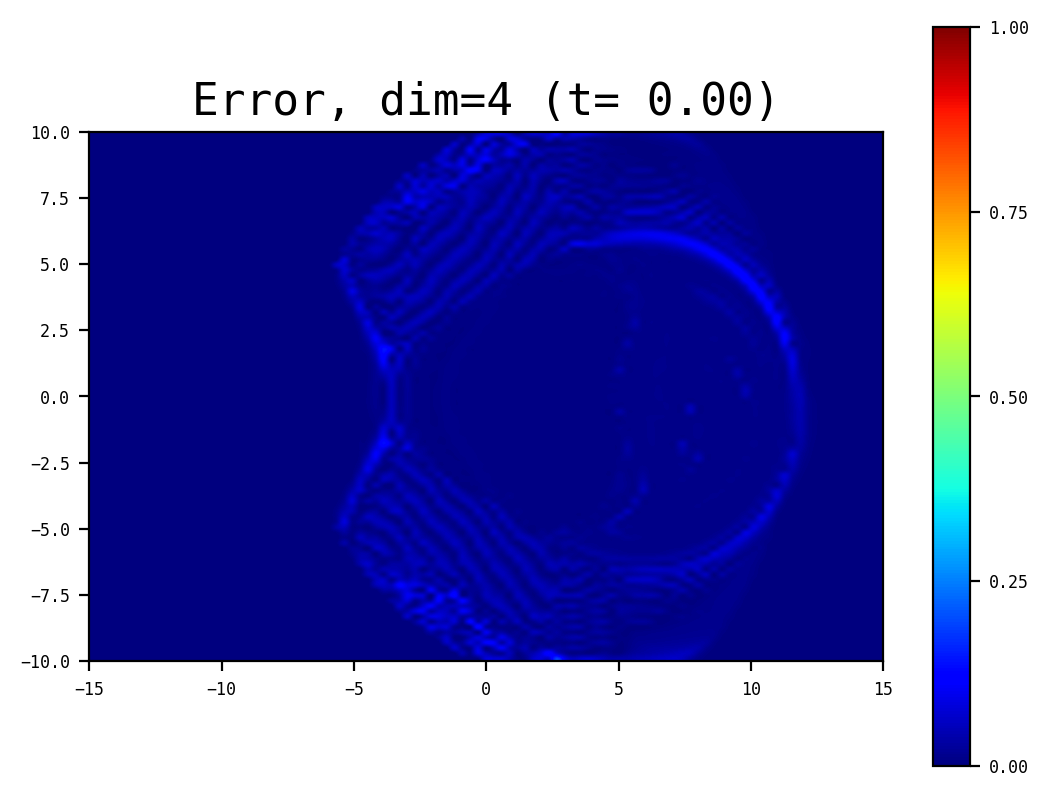}\\
}       
        \caption {(Example 4) 
        Comparison between \SLTWO (left), \SLTWOH (middle) and \SLTHREE (right),
        for $N=8$ (top) and $N=16$ (bottom).}
	\label{fig:exeiknc-comp}
\end{figure}

{\em Comparison of schemes in dimension $d=4$.}
First, the SL-, H- and L-schemes are compared.
Neural networks with 3 layers of 60 neurons are chosen, and each simulation uses 100,000 stochastic gradient step. 
Figure (\ref{fig:exeiknc-comp}) displays the error $|\hat{V}_0(x) - v(0,x)|$ in fonction of space, with $N\in\left\{8,16\right\}$ number of time steps. 
Results are shown in Table~\ref{tab:exeiknccomp-errors}.

The \SLTWO approximates both the control and the value function 
by neural networks. The projection of the latter is a source of errors, that accumulates
during the simulation. This drawback is avoided with the \SLTWOH and \SLTHREE, 
where the value function is computed as a composition of the (exact) target function $\varphi$
and the approximated controls. Again, for the error, the \SLTWOH and \SLTHREE behave better than the \SLTWO.

Furthermore, for the H-scheme and the L-scheme, when $N$ varies from $8$ to $16$, 
we observe very roughly that the $L^1$ (global and local) errors are divided by a factor two (this is less clear for the $L^\infty$ errors). This is not the case for the SL-scheme, for which errors have a tendency to accumulate more with time iterations.

\begin{table}[H]
	\centering
	\begin{tabular}{|c|cccccc|cc|cc|c|} \hline
	\multirow{2}*{Scheme} & \multicolumn{6}{c|}{Parameters} & \multicolumn{2}{c|}{Global errors} & \multicolumn{2}{c|}{Local errors} & \multirow{2}*{$\begin{matrix}\text{CPU} \\ \text{time} \end{matrix}$} \\ 
	\cline{2-11} & $d$ & $N$ & lay. & neur. & $M$ & S.G. it. & $L_{\infty}$ & $L_{1}$ rel. & $L_{\infty}$ & $L_{1}$ rel. & \\ 
	\hline \hline 
	SL & 4 &  8 & 3 &  60 & 4000 & 100000 & 8.60e-01 & 4.91e-02 & 1.98e-01 & 7.48e-02 & 6h14 \\ \hline 
	H & 4 &  8 & 3 &  60 & 4000 & 100000 & 3.34e-01 & 1.58e-02 & 2.01e-01 & 4.61e-02 & 11h26 \\ \hline 
	L & 4 &  8 & 3 &  60 & 4000 & 100000 & 3.24e-01 & 6.68e-03 & 1.09e-01 & 2.56e-02 & 8h16 \\ \hline \hline
	SL & 4 & 16 & 3 &  60 & 4000 & 100000 & 1.07e+00 & 9.21e-02 & 2.92e-01 & 1.06e-01 & 12h29 \\ \hline 
	H & 4 & 16 & 3 &  60 & 4000 & 100000 & 3.32e-01 & 9.62e-03 & 1.54e-01 & 2.88e-02 & 34h26 \\ \hline 
	L & 4 & 16 & 3 &  60 & 4000 & 100000 & 1.85e-01 & 3.57e-03 & 8.85e-02 & 1.64e-02 & 28h08 \\ \hline 
	\end{tabular}
	\caption{(Example 4) comparison between schemes
	\label{tab:exeiknccomp-errors}
	}
\end{table}

{\em Test of the L-scheme for increasing dimensions.} 
Next, the \SLTHREE is tested for several dimensions $d\in \left\{2,4,6,8\right\}$, and
results are given in Table~\ref{tab:exeiknc-errors}.
The neural network size is kept constant, with 3 layers of 60 neurons, 
as for the number of time iterations ($N=8$). 
In order to reach comparable precision, we have observed that the number of stochastic gradient iterations
has to grow with $d$ (as the dimension increases, more iterations are needed to explore the whole region of interest).
Otherwise, the scheme is relatively robust with respect to 
the physical dimension of the problem (see Fig.~\ref{fig:exeiknc}).

We observe for dimension $d=8$ some defects in the numerical solution (some oscillations appears). 
Because of CPU time limitations, we did not attempt using more S.G. iterations, although in principle (as observed for lower dimensions)
this should enable a better optimization and solve the problem.

\begin{table}[H]
	\centering
	\begin{tabular}{|c|cccccc|cc|cc|c|} \hline
	\multirow{2}*{Scheme} & \multicolumn{6}{c|}{Parameters} & \multicolumn{2}{c|}{Global errors} & \multicolumn{2}{c|}{Local errors} & \multirow{2}*{$\begin{matrix}\text{CPU} \\ \text{time} \end{matrix}$} \\ 
	\cline{2-11} & $d$ & $N$ & lay. & neur. & $M$ & S.G. it. & $L_{\infty}$ & $L_{1}$ rel. & $L_{\infty}$ & $L_{1}$ rel. & \\ 
	\hline \hline 
	L & 2 &  8 & 3 &  60 & 4000 & 50000 & 2.66e-01 & 5.99e-03 & 1.19e-01 & 4.61e-02 & 3h02 \\ \hline 
	L & 4 &  8 & 3 &  60 & 4000 & 100000 & 3.90e-01 & 6.77e-03 & 1.16e-01 & 2.69e-02 & 8h13 \\ \hline 
	L & 6 &  8 & 3 &  60 & 4000 & 400000 & 9.69e-01 & 1.09e-02 & 1.78e-01 & 2.88e-02 & 35h20 \\ \hline 
	L & 8 &  8 & 3 &  60 & 4000 & 600000 & 1.05e+00 & 3.75e-02 & 1.71e-01 & 2.95e-02 & 45h27 \\ \hline 
	\end{tabular}
	\caption{(Example 4) $L$-scheme, dimensions $d=2,4,6,8$
	\label{tab:exeiknc-errors}
	}
\end{table}

\begin{figure}[!hbtp]
\VOID{
	\centering
	\renewcommand{\figw}{.32\linewidth}
	\renewcommand{\rdir}{numeiknc/}
	\renewcommand{\ldir}{EIKNCD2NEXT500}
	\includegraphics[width=\figw]{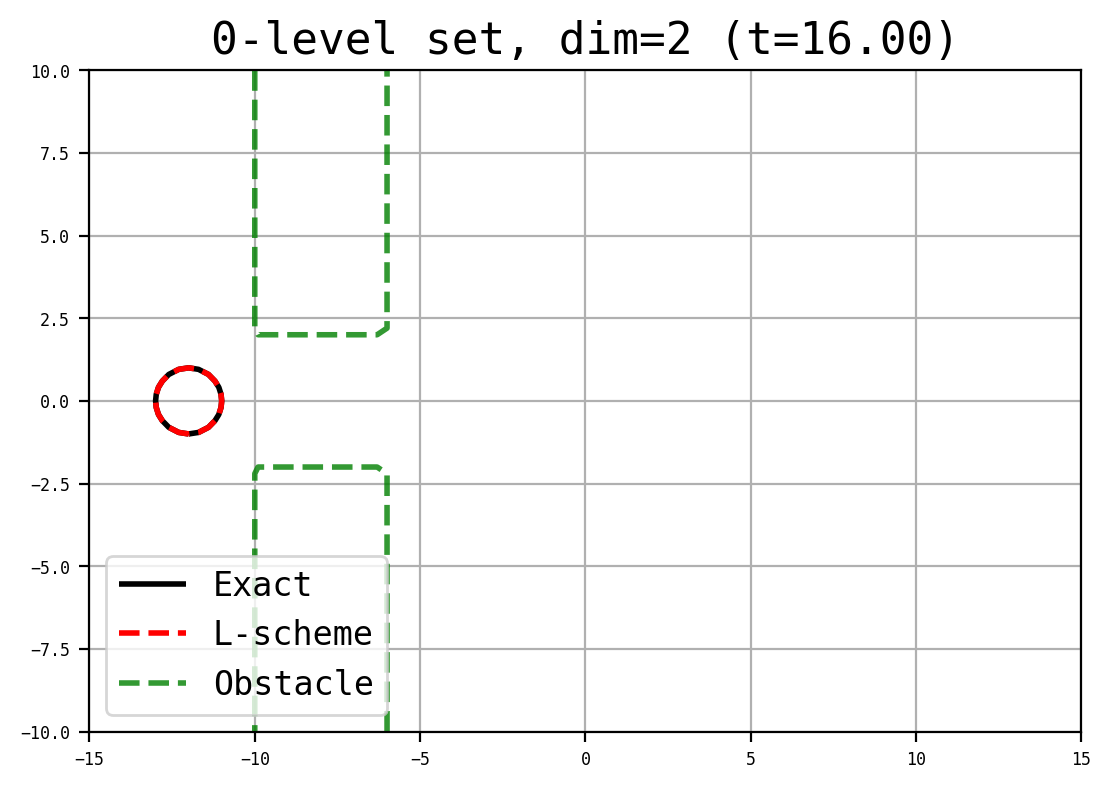}%
	\includegraphics[width=\figw]{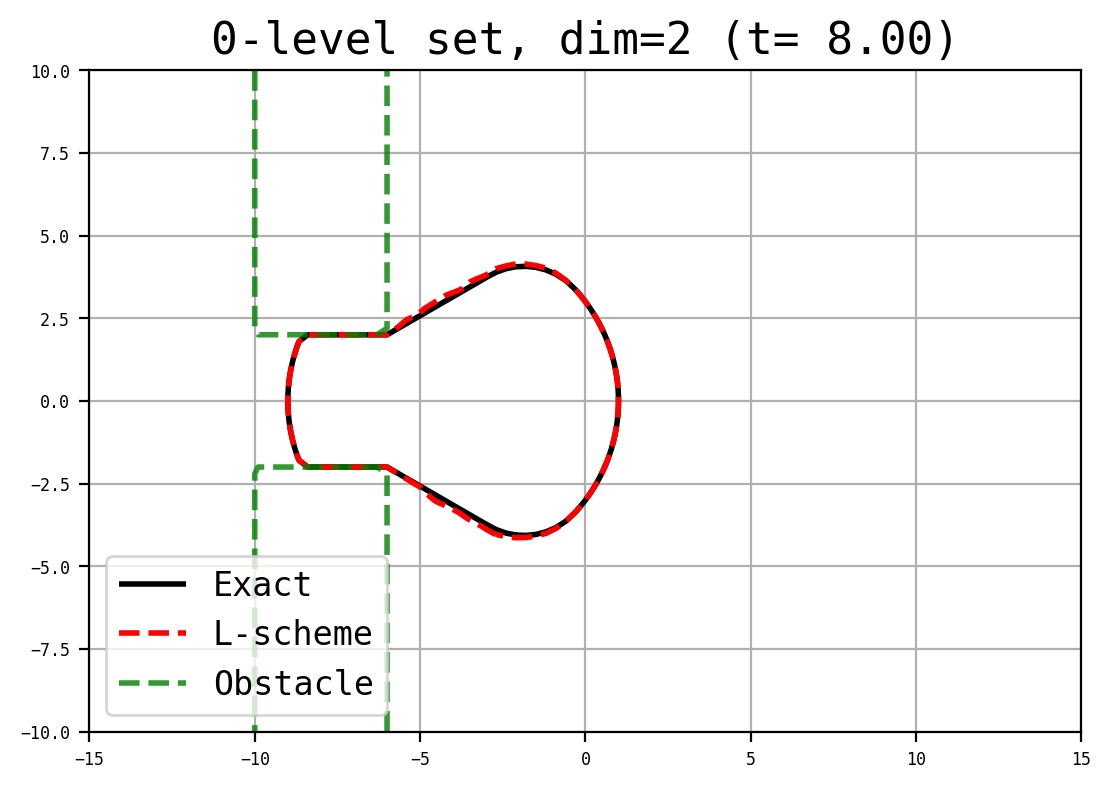}%
	\includegraphics[width=\figw]{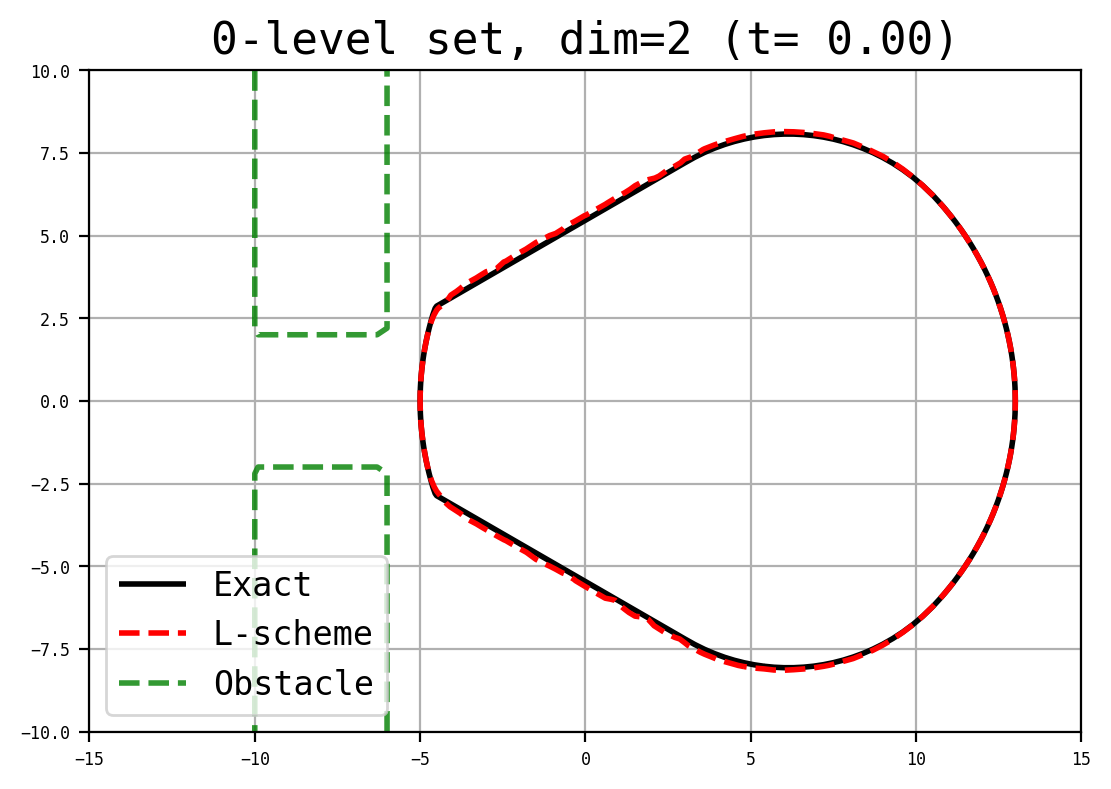}\\
	\renewcommand{\ldir}{EIKNCD4NEXT1000}
	\includegraphics[width=\figw]{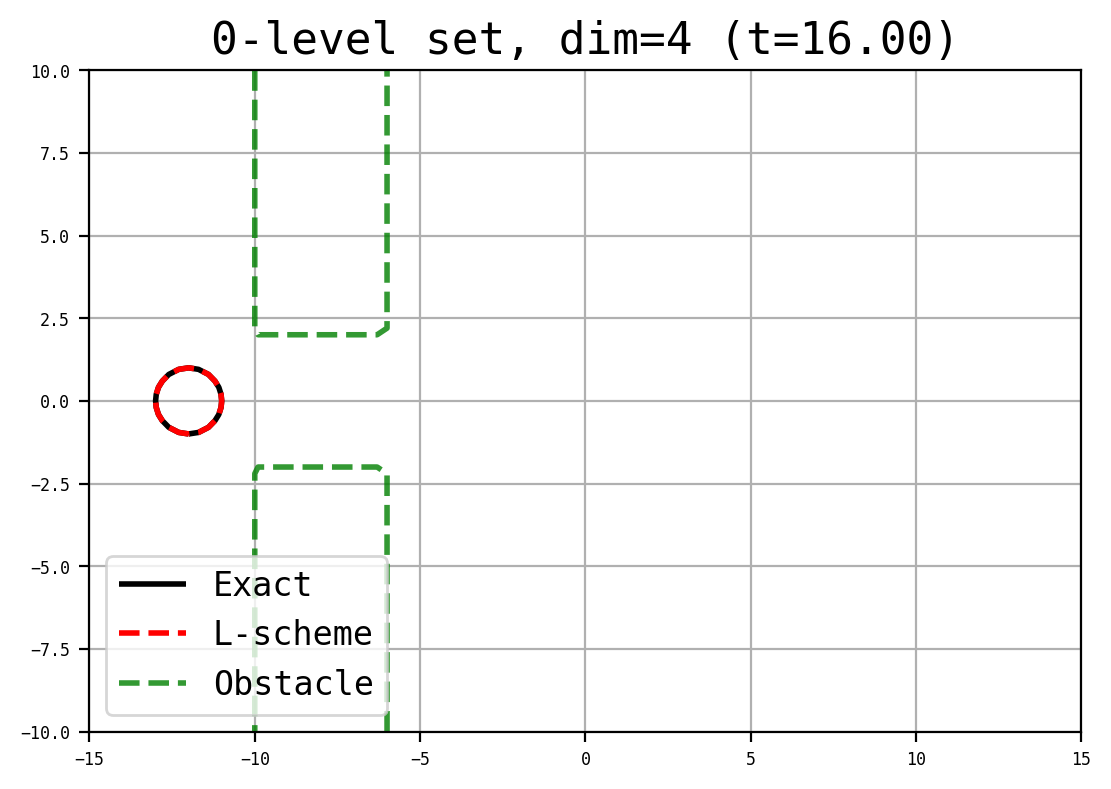}%
	\includegraphics[width=\figw]{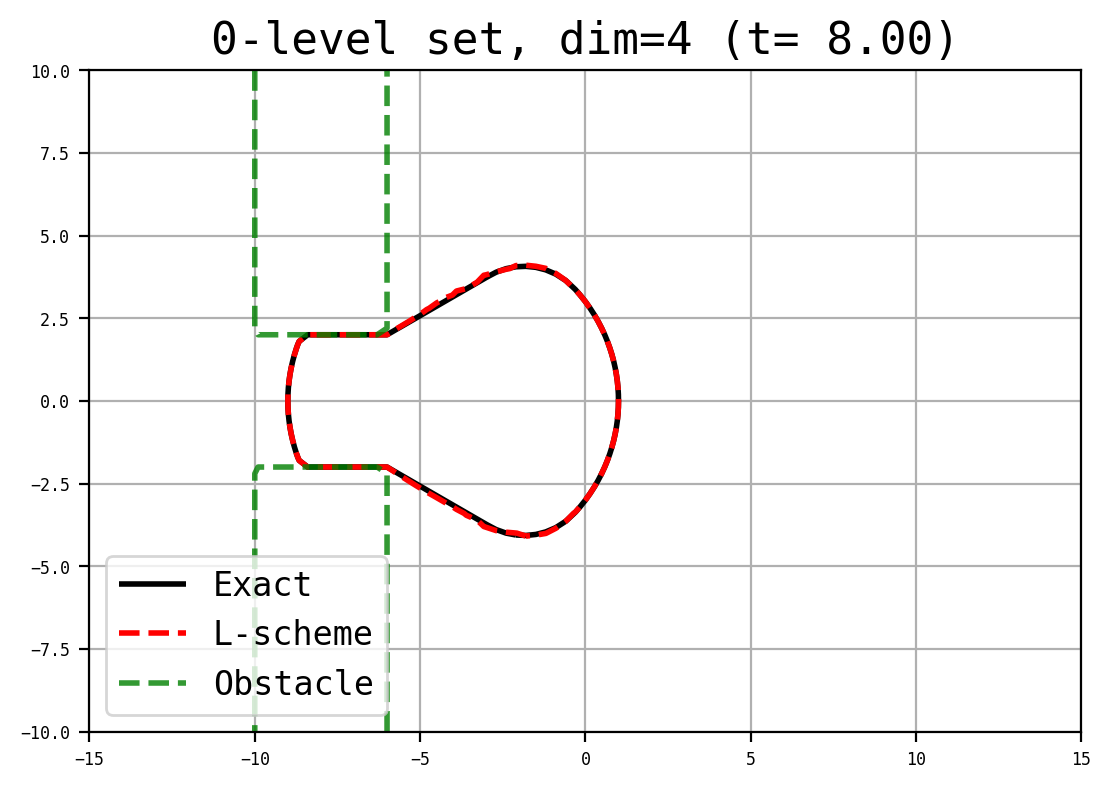}%
	\includegraphics[width=\figw]{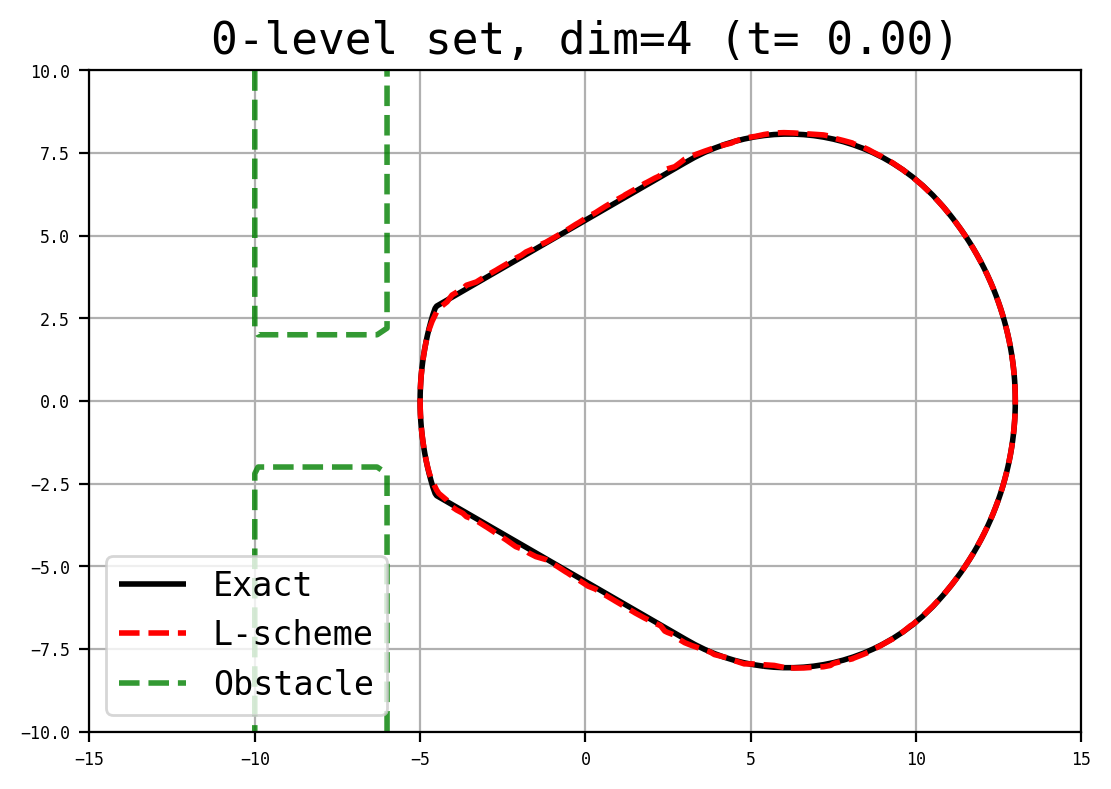}\\
	\renewcommand{\ldir}{EIKNCD6NEXT4000}
	\includegraphics[width=\figw]{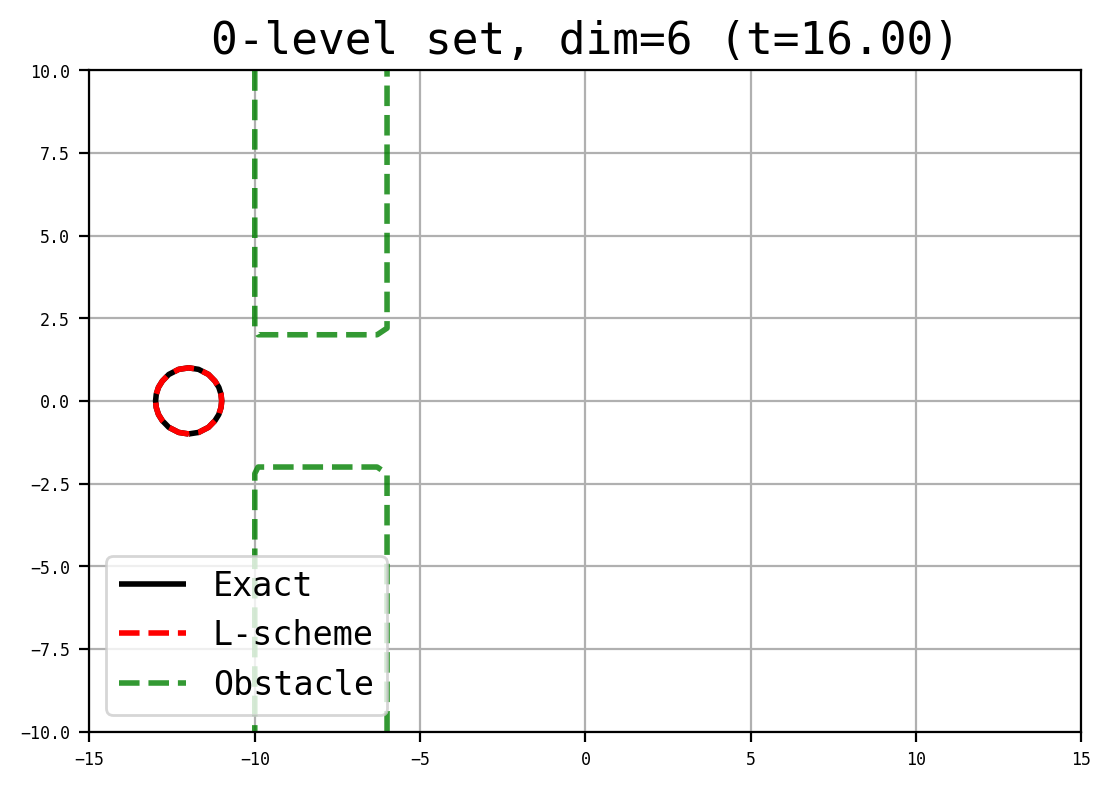}%
	\includegraphics[width=\figw]{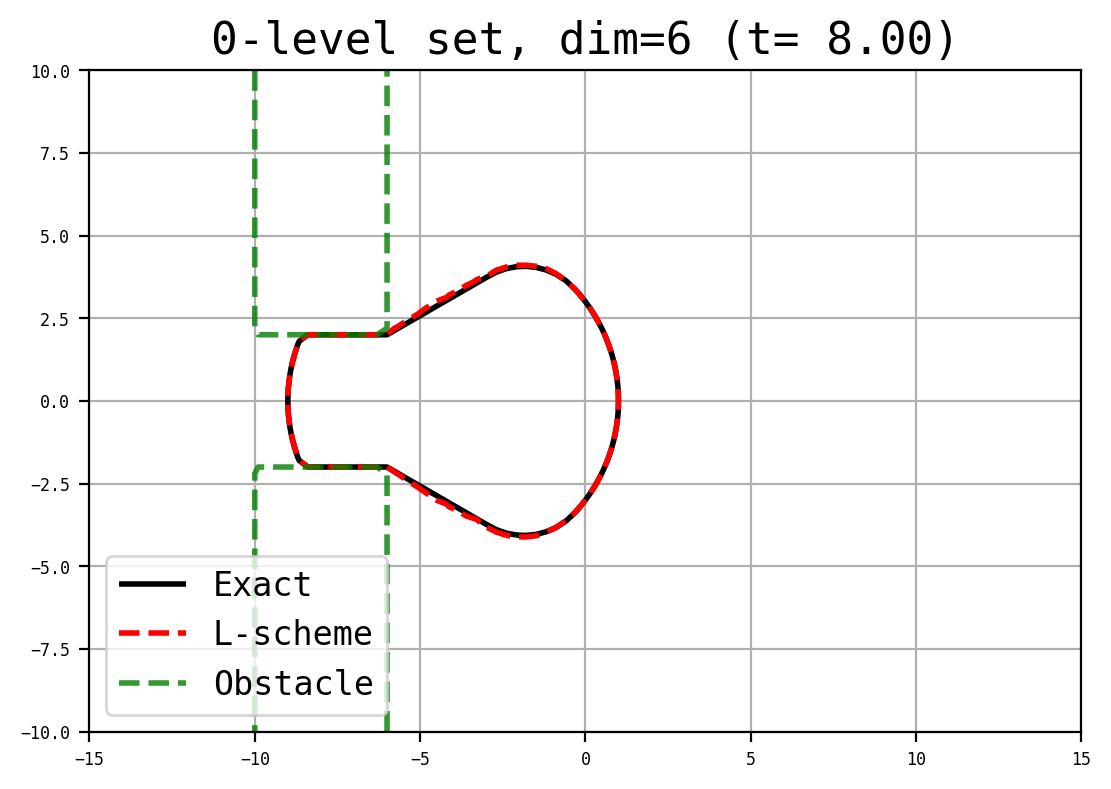}%
	\includegraphics[width=\figw]{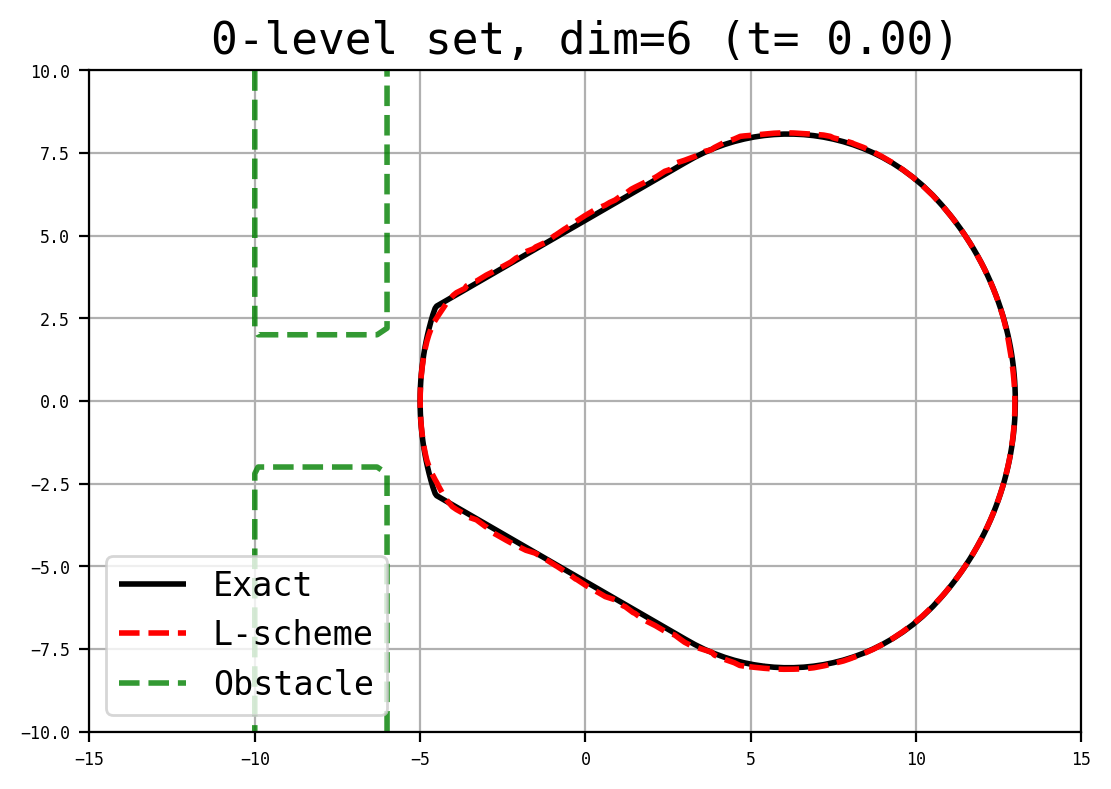}%
	\renewcommand{\ldir}{EIKNCD8NEXT6000}
	\includegraphics[width=\figw]{\rdir\ldir/%
        show_fig_zerolevel0_SL3_D8_Nit8_Arelu_Neu60_Lay3_Next6000_Order20_Batch4000_MH0_VbistruncLmax_nbOpt1_FeedForward01_nIRclip}
	\includegraphics[width=\figw]{\rdir\ldir/%
        show_fig_zerolevel4_SL3_D8_Nit8_Arelu_Neu60_Lay3_Next6000_Order20_Batch4000_MH0_VbistruncLmax_nbOpt1_FeedForward01_nIRclip}
	\includegraphics[width=\figw]{\rdir\ldir/%
        show_fig_zerolevel8_SL3_D8_Nit8_Arelu_Neu60_Lay3_Next6000_Order20_Batch4000_MH0_VbistruncLmax_nbOpt1_FeedForward01_nIRclip}
}       
	\caption {
        (Example \refeiknc)
	\SLTHREE, dimensions 2,4,6,8, $N=8$ time steps, networks : 3 layers of 60 neurons.}
	\label{fig:exeiknc}
\end{figure}


\subsection{Example \refeikc : eikonal advection equation with obstacle, small drift} \label{sec:eikc}

We now turn on a similar example as in example~\refeiknc, excepted for the coefficients which are now 
$$
  c=1 , \quad b= 0.5. 
$$

Results obtained with the \SLTHREE are given in Table~\ref{tab:exeikc-errors} and Fig.~\ref{fig:exeikc}.
Here $|b|<c$, the drift is small (this corresponds to a controllable situation). 
We observe that the front has to negociate a sharper angle 
near the boundary of the tube (see Fig.~\ref{fig:exeikc}).

As in example~\refeiknc, the results are rather robust with respect to dimension, 
provided the number of stochastic gradient iterations is large enough.

\begin{table}[H]
	\centering
	\begin{tabular}{|c|cccccc|cc|cc|c|} \hline
	\multirow{2}*{Scheme} & \multicolumn{6}{c|}{Parameters} & \multicolumn{2}{c|}{Global errors} & \multicolumn{2}{c|}{Local errors} & \multirow{2}*{Time} \\ 
	\cline{2-11} & $d$ & $N$ & lay. & neur. & $M$ & S.G it. & $L_{\infty}$ & $L_{1}$ rel. & $L_{\infty}$ & $L_{1}$ rel. & \\ 
	\hline \hline 
	L & 2 &  8 & 3 &  60 &  4000 & 50000 & 1.24e-01 & 2.94e-03 & 8.81e-02 & 5.21e-03 & 3h09 \\ \hline 
	L & 4 &  8 & 3 &  60 &  4000 & 100000 & 2.49e-01 & 4.70e-03 & 8.67e-02 & 5.45e-03 & 6h51 \\ \hline 
	L & 6 &  8 & 3 &  60 &  4000 & 400000 & 8.74e-01 & 3.70e-02 & 1.09e-01 & 1.01e-02 & 35h07 \\ \hline 
	\end{tabular}
	\caption{Errors for example \arabic{subsection}
	\label{tab:exeikc-errors}
	}
\end{table}

\newcommand{\figh}{0.33}

\begin{figure}[H]
\VOID{
  \centering
  \renewcommand{\figw}{.25\linewidth}
  \renewcommand{\figh}{.25\linewidth}
  \renewcommand{\rdir}{numeikc/}
  \renewcommand{\ldir}{EIKOD2500EXT}
  \includegraphics[width=\figw,height=\figh]{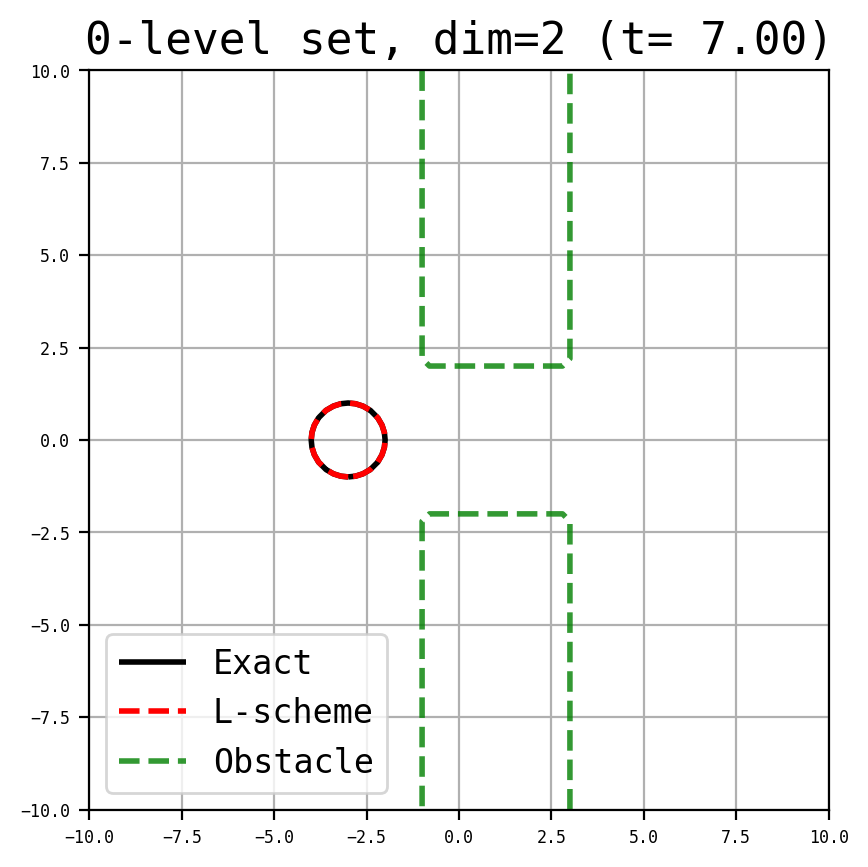}%
  \includegraphics[width=\figw,height=\figh]{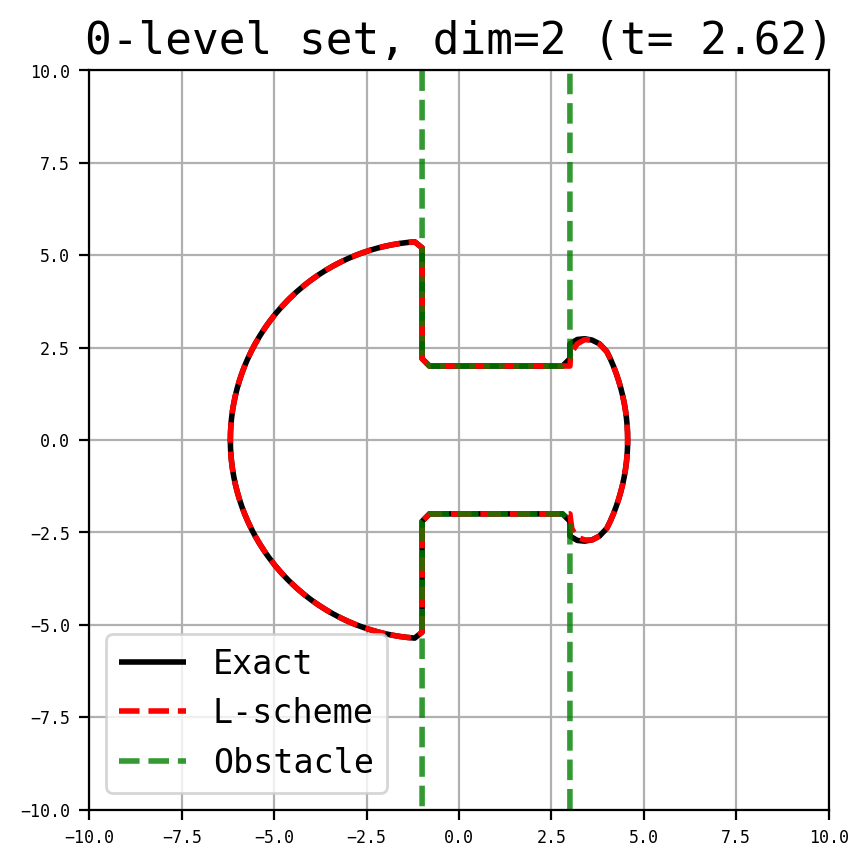}%
  \includegraphics[width=\figw,height=\figh]{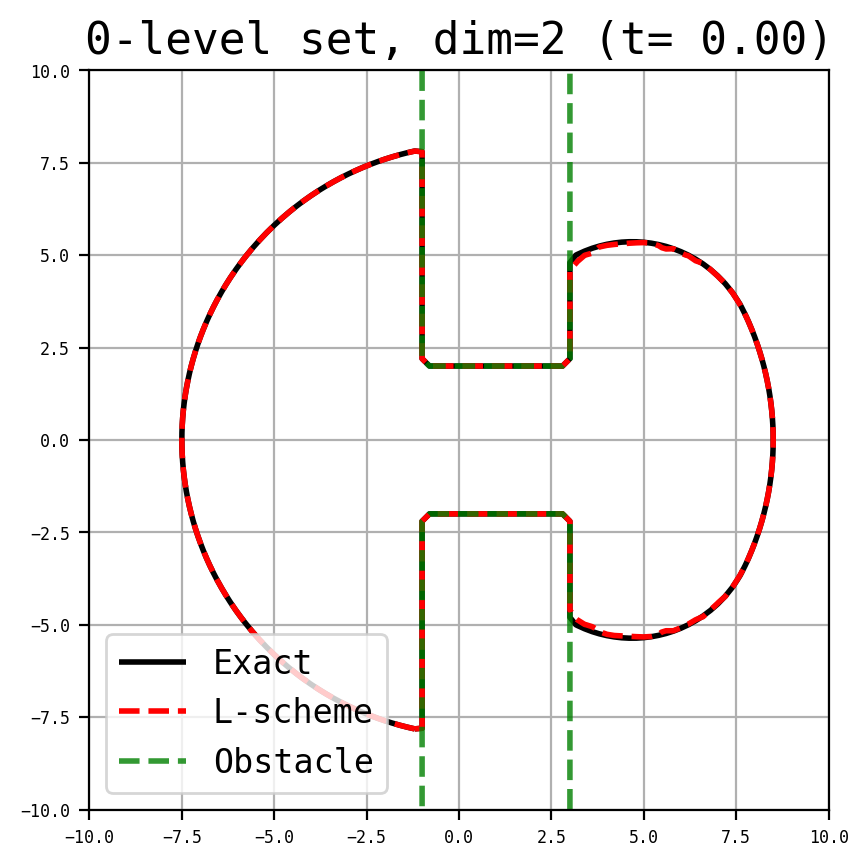}\\
	
  \renewcommand{\ldir}{EIKOD41000EXT}
  \includegraphics[width=\figw,height=\figh]{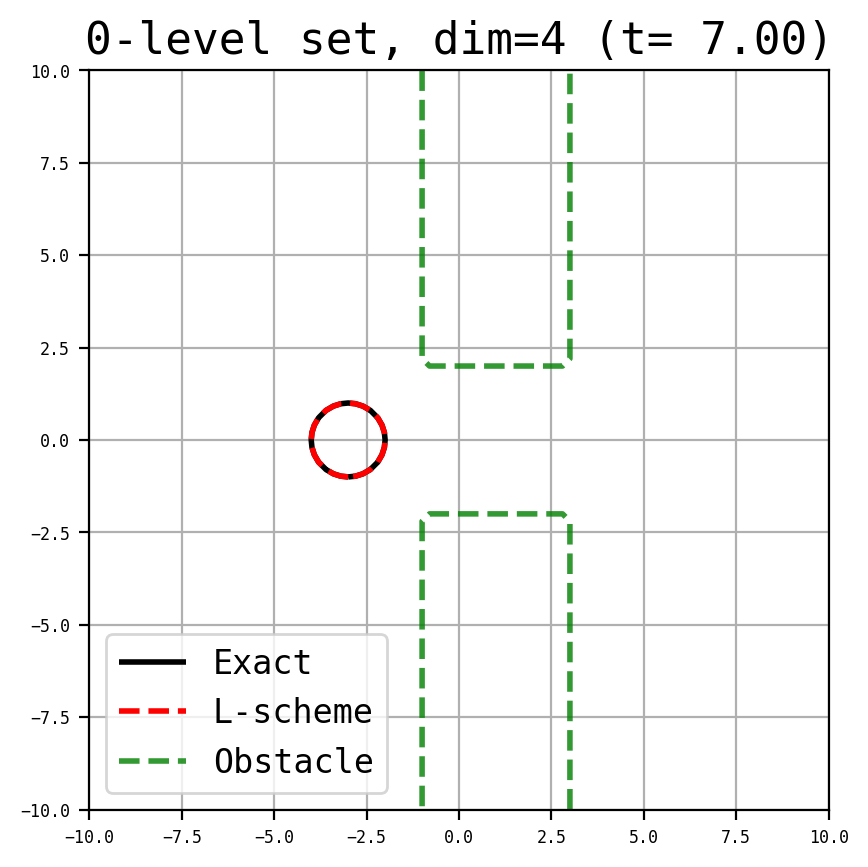}%
  \includegraphics[width=\figw,height=\figh]{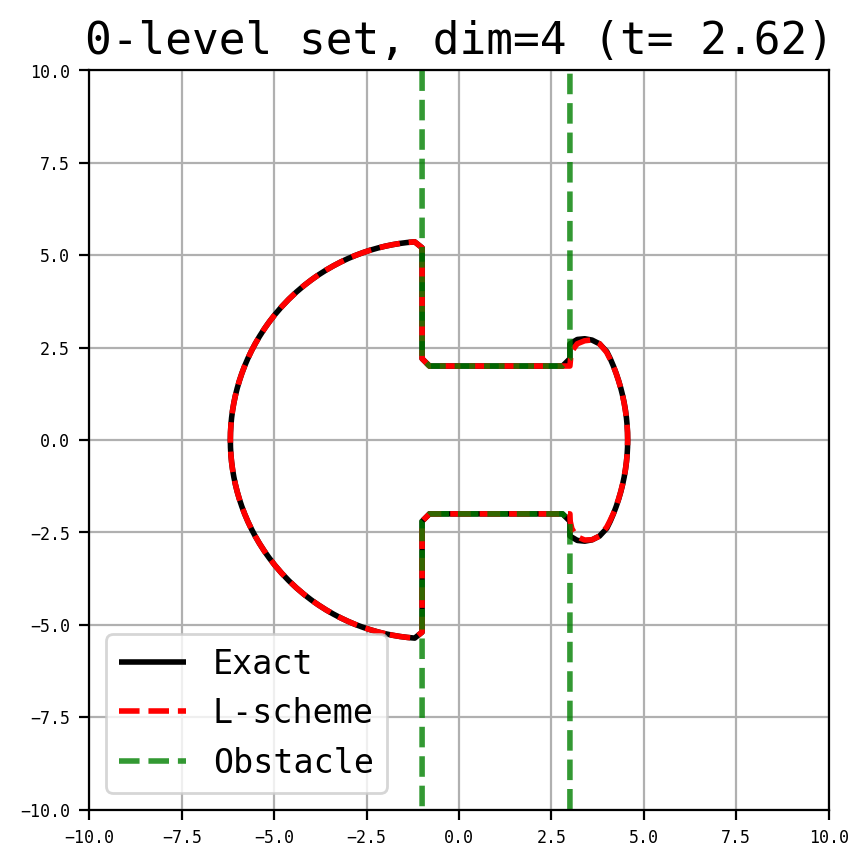}%
  \includegraphics[width=\figw,height=\figh]{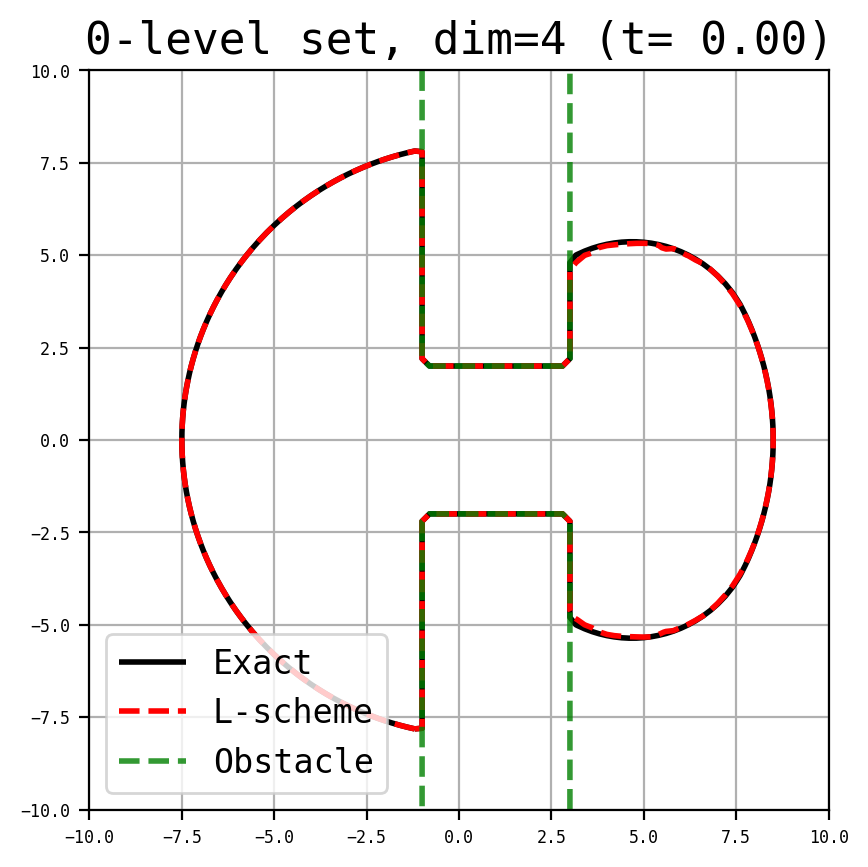}\\
	
  \renewcommand{\ldir}{EIKOD64000EXT}
  \includegraphics[width=\figw,height=\figh]{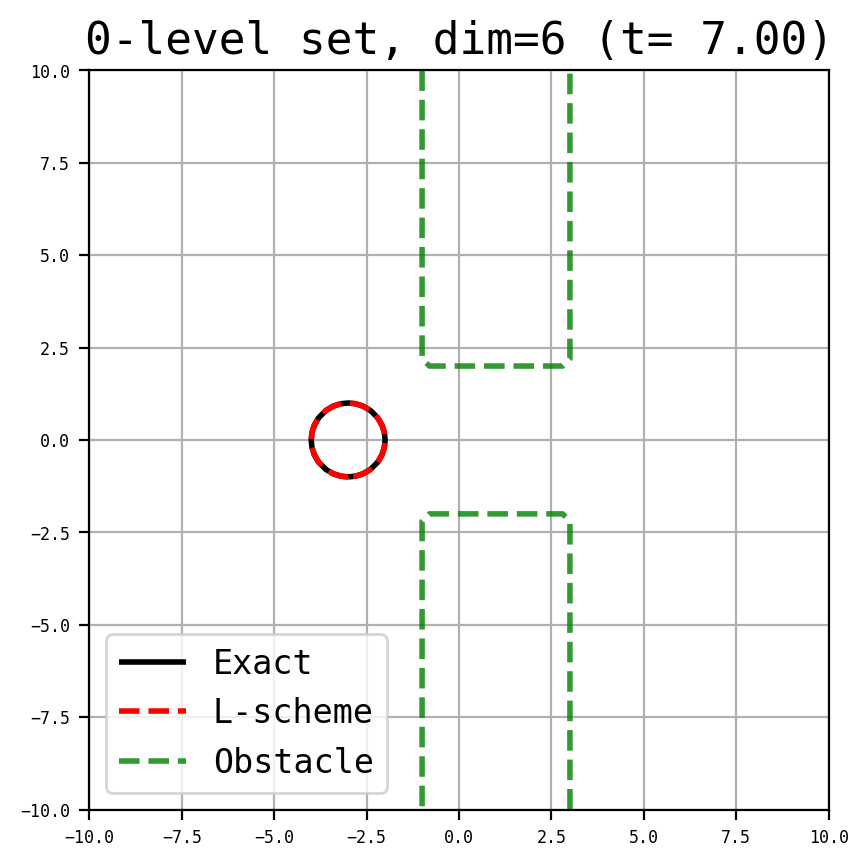}%
  \includegraphics[width=\figw,height=\figh]{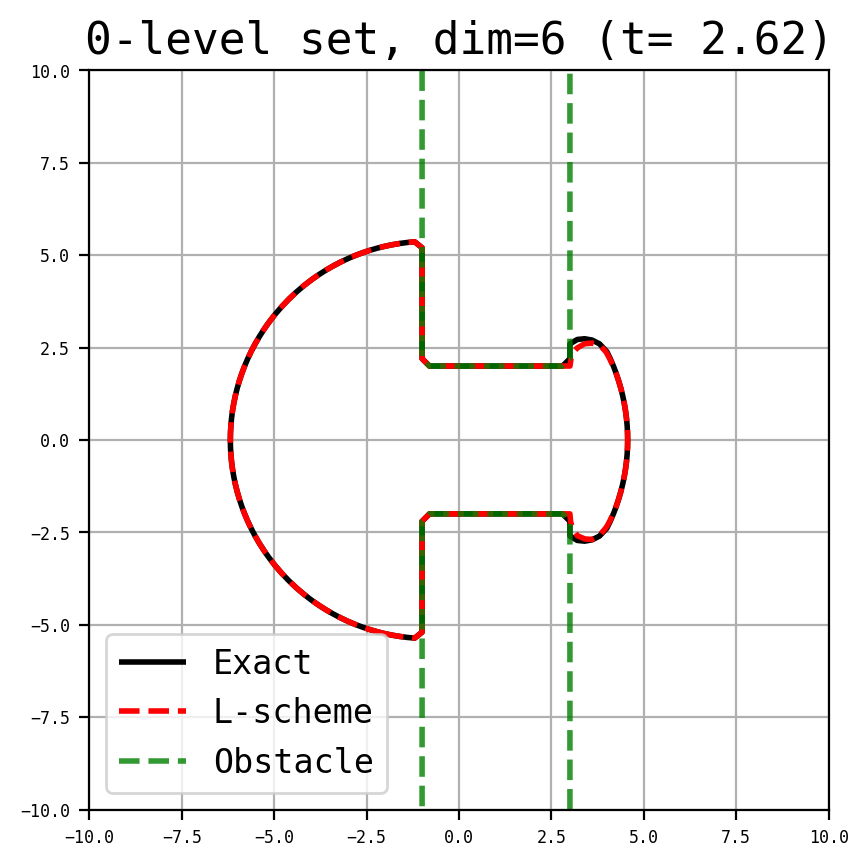}%
  \includegraphics[width=\figw,height=\figh]{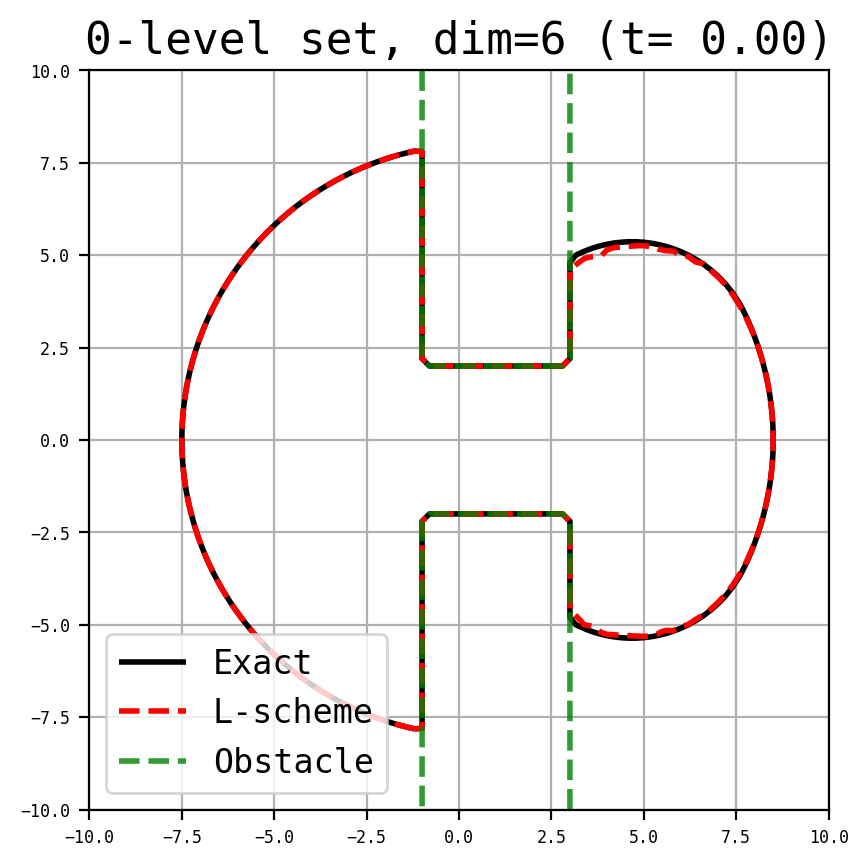}%
} 
\caption {(Example \arabic{subsection})
\SLTHREE, dimensions 2, 4 and 6. Networks with 3 layers of 60 neurons, $N=8$ time steps.
\label{fig:exeikc}
}
\end{figure}

\appendix               
\section{Semi-analytical solution for examples \refeiknc\ and \refeikc}\label{app:A}

We briefly 
describe how to compute the exact values for examples~\refeiknc\ and~\refeikc,
that is, in order  to compute $v=v(t,x)$ for given values $t$ and $x$.
We consider the case of data $\varphi(x) = \|x\| + \alpha_{\min}$, for a $\alpha_{\min} \in \mathbb{R}$.  

\begin{figure}[H]
  \centering
  \includegraphics[width=0.32\linewidth]{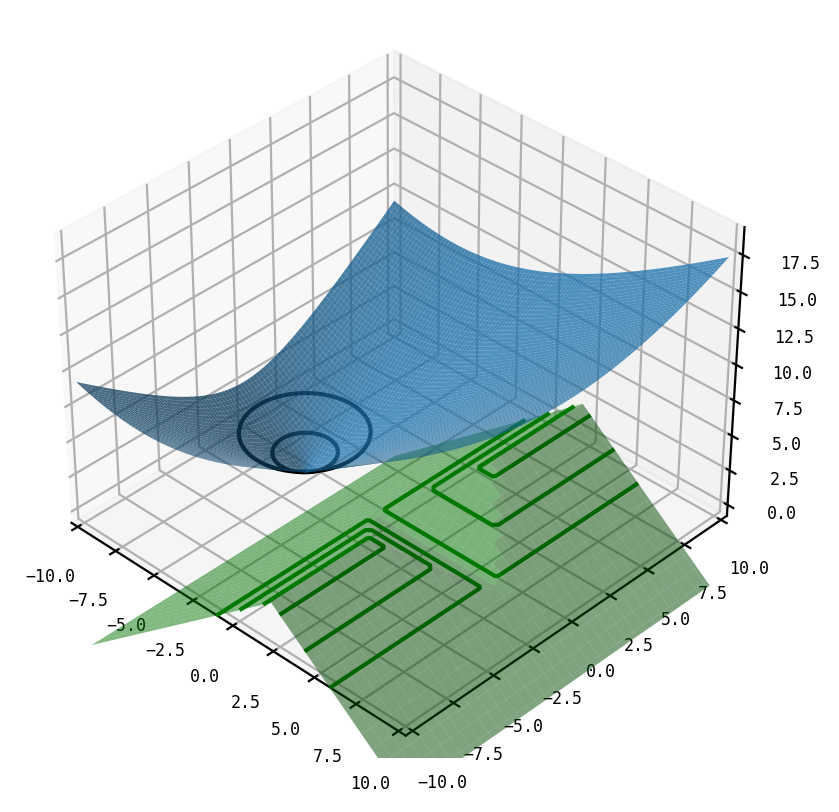}
  \includegraphics[width=0.32\linewidth]{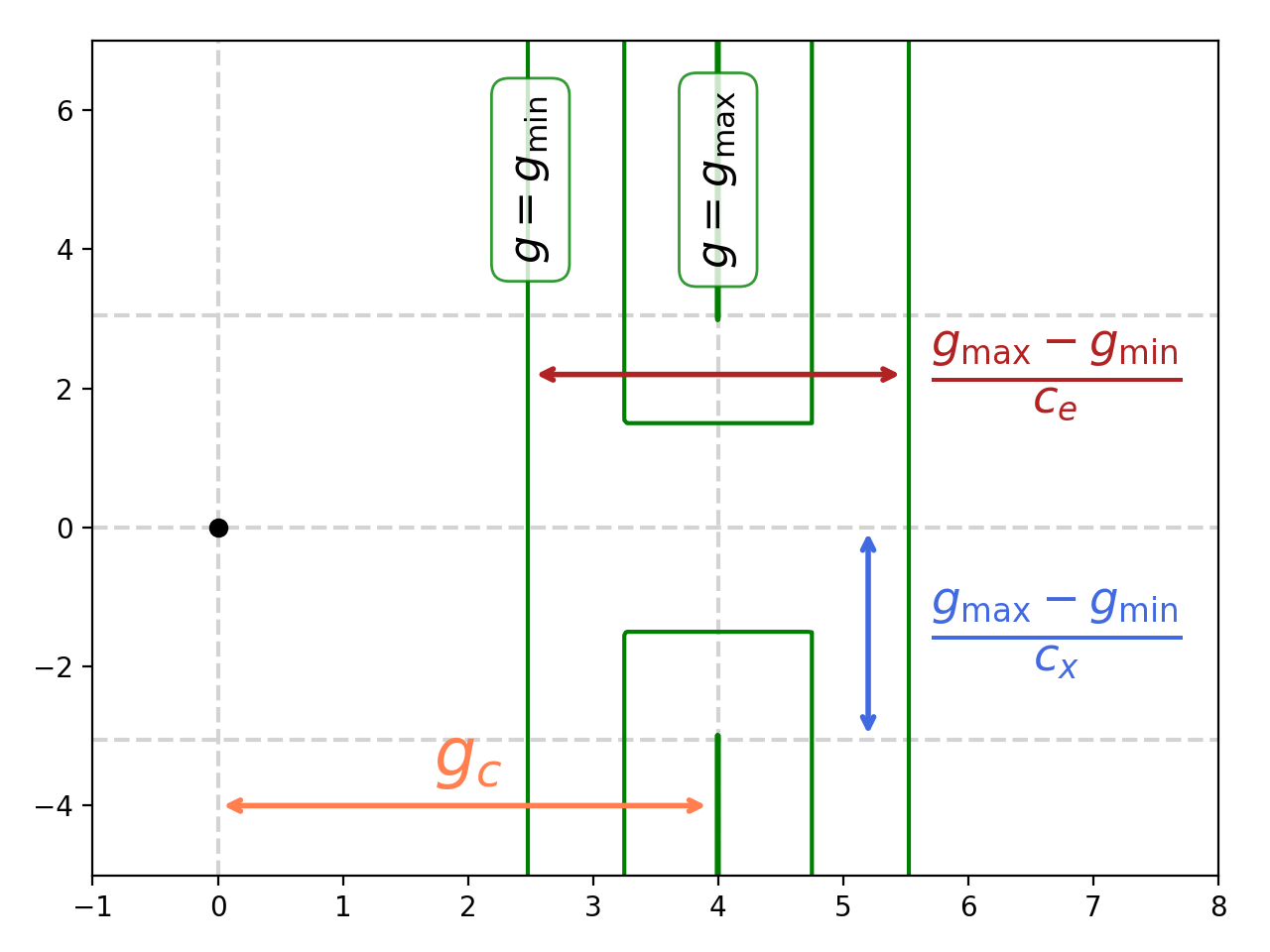}
  \includegraphics[width=0.32\linewidth]{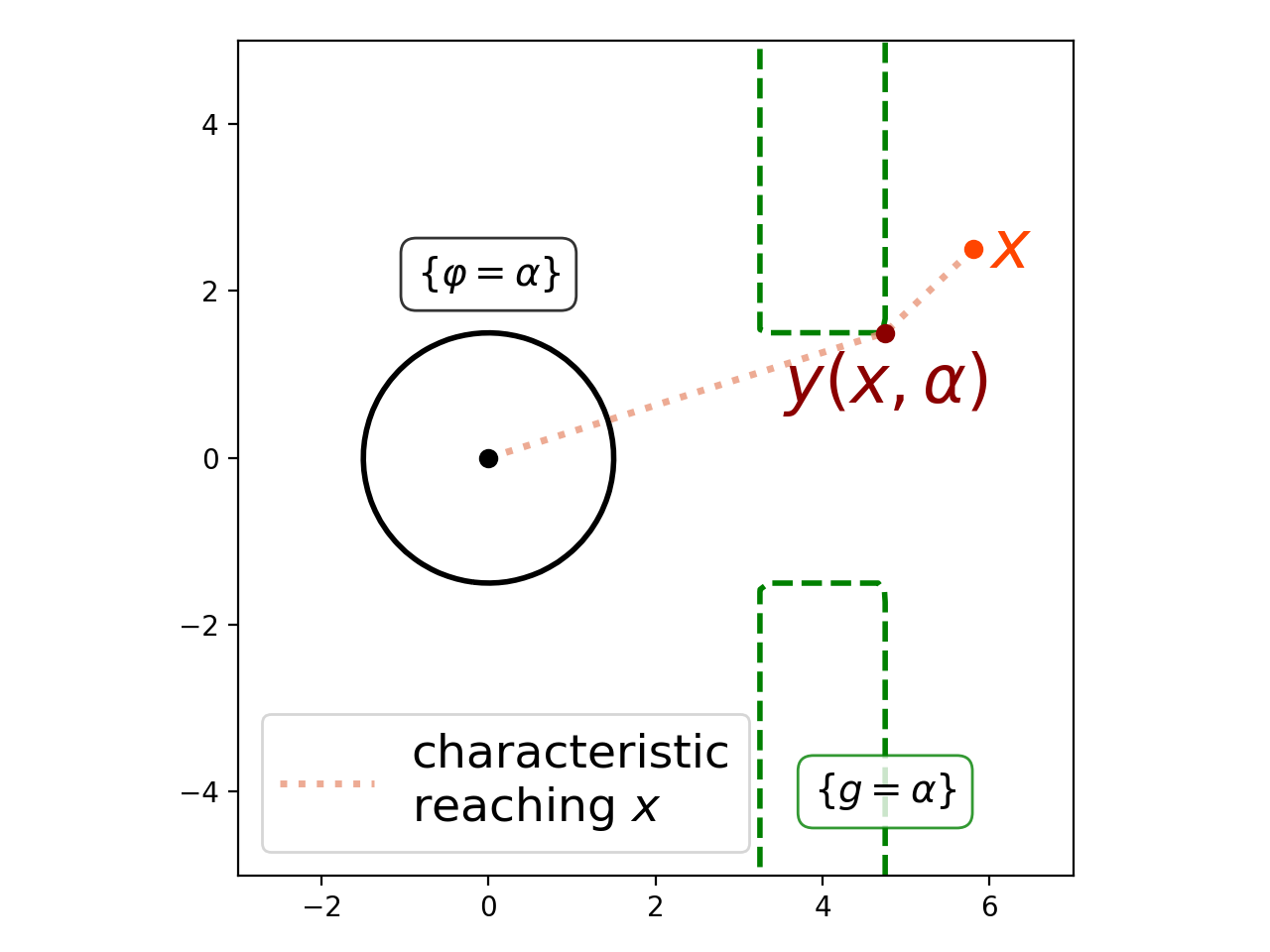}
  \caption {Initial condition (in blue) and obstacle function (in green).
  Illustration of parameters. Exemple of a characteristics with obstacle in the case $b=0$.}
  \label{fig:zermelo_obstacle_anatomy}
\end{figure}

For $\alpha \in [g_{\min}, g_{\max}]$, notice that the level set $\left\{g=\alpha\right\}$ 
corresponds to a wall pierced by a square door (see figure \ref{fig:zermelo_obstacle_anatomy}).

First, by using the symmetries of the problem around the $e_1$ axis, and an orthogonal axis along $x-x_1 e_1$,
it is possible to set back the problem into a 2-dimensional problem.

If no obstacle term is present (or if it does not modify the value $v$),
for a given point $x\in \R^d$, and for a given (level-set) value $v$, 
it is possible to compute the minimal time to reach $x$ from the initial level set front (corresponding to some point 
on the level set $\{\varphi(x)=v\}$).
More precisely, the optimal trajectory is a straight segment and the value satisfies
$$  
   v(t,x)= (\|x - b e_1 t \| -ct)_+ + \ma_{\min},
$$
and $t$ is also the time for the front $\{\varphi(.)=v\}$ to reach the point $x$.

In the more complex situation when the optimal trajectory from $\{\varphi(.)=v\}$ to point $x$ is not a straight segment,
it will be composed of two segments (one starting from some point on $\{\varphi(.)=v\}$ to reach some point $y$ on the level set $\{g(y)=v\}$,
the other one starting from $y$ to $x$ - as depicted in Fig.~\ref{fig:zermelo_obstacle_anatomy}-right).
Then the minimal time $t$ (associated to some value $v=v(t,x)$) is the sum of two minimal times $t_1$, $t_2$
such that 
\begin{subequations} \label{eq:NEWT}
\be
  & &   t=t_1 + t_2  \label{eq:NEWT-1}\\
  & &   \| y - b e_1 t_1\|^2 = (c t_1  + v - \ma_{\min})^2 \label{eq:NEWT-2}\\
  & &   \| x - (y   + b e_1 t_2)\|^2 = (c t_2)^2. \label{eq:NEWT-3} 
\ee
\end{subequations}
Note that $t_1$ is the time for the level set $\{\varphi(.)=v\}$ to reach $y$ on $\{g(.)=v\}$, $t_2$ is the time to reach
$x$ from $y$. Then, for a given value $v$,  $y=y(v)$ is known (it has an affine analytic expression in term of $v$). 
Times $t_1=t_1(v)$ and $t_2=t_2(v)$ are obtained as root solutions of 
\eqref{eq:NEWT-2} and \eqref{eq:NEWT-3}. 
Then, for given $(t,x)$, the value $v$ is obtained through a Newton algorithm for solving system~\eqref{eq:NEWT}. 
Note that on regions not attained by the front, we consider the complex roots $t_1(v)$ or $t_2(v)$ in the Newton algorithm
(in order to always have well-defined and continuous reaching times).




\bibliographystyle{abbrv}
\hypersetup{
    linkcolor=black, 
    citecolor=black, 
    urlcolor=black  }
\bibliography{mtns,comp,other}


\end{document}